\documentclass[12pt,a4paper]{article}
\usepackage[latin1]{inputenc} 
\usepackage[T1]{fontenc} 
\usepackage[english]{babel}
\usepackage{libertine}
\usepackage{libertinust1math}
\usepackage{amsmath}
\usepackage{amsfonts}
\usepackage{amssymb}
\usepackage{amsthm}
\usepackage{mathrsfs}
\usepackage{geometry}
\usepackage{hyperref}
\usepackage{accents}
\usepackage{etoolbox}
\usepackage{stmaryrd}
\usepackage{mathtools}
\usepackage[scr=boondox]{mathalfa}
\usepackage{pgf,tikz-cd,pgfplots,xcolor}
\usepackage[all]{xy}
\usetikzlibrary{arrows}
\pgfplotsset{compat=1.15}
\title{Unobstructedness of deformations for a d-semistable central corank one boundary point}
\author{Emeryck Marie}
\date{}

\newtheoremstyle{emi}{\topsep}{\topsep}{\itshape}{0pt}{}{\hspace{1ex}---}{1ex}{\thmname{{\itshape#1}}\thmnumber{ {\bfseries(#2.)}}\thmnote{ (#3.)}}
\theoremstyle{emi}
\newtheorem{theoreme}{Theorem}[subsection]
\newtheorem{lemme}[theoreme]{Lemma}
\newtheorem{proposition}[theoreme]{Proposition}
\newtheorem{definition}[theoreme]{Definition}
\newtheorem{corollaire}[theoreme]{Corollary}

\newtheoremstyle{emiparagraphe}{\topsep}{\topsep}{}{0pt}{}{\hspace{1ex}---}{1ex}{\thmnumber{{\bfseries#2.}}\thmnote{ (#3.)}}
\theoremstyle{emiparagraphe}
\newtheorem{paragraphe}[theoreme]{}

\newtheoremstyle{emiremarque}{\topsep}{\topsep}{}{0pt}{}{\hspace{1ex}---}{1ex}{\thmname{{\itshape#1}}\thmnote{ (#3.)}}
\theoremstyle{emiremarque}
\newtheorem*{remarque}{Remark.}
\newtheorem*{preuve}{Proof.}
\AtEndEnvironment{preuve}{\qed}

\begin{document}

\maketitle

\begin{center}
    \small{Ruhr Universit\"{a}t Bochum, Fakult\"{a}t f\"{u}r Mathematik; Universit\"{a}tsstra\ss e 150, 44801 Bochum, Germany \\  
    \textsc{E-Mail}: emeryck.marie@ruhr-uni-bochum.de}
\end{center}

\begin{center}
    \textsc{Abstract:}
\end{center}
In this paper, we prove the unobstructedness of the functor of log smooth deformations for an irreducible type II degeneration of complex abelian surfaces --- shifted gluing of a $\mathbb{P}^1$-bundle over an elliptic curve along two disjoint sections --- using the $T^1$-lifting technique and the degeneration on the first page of a Hodge--de Rham like spectral sequence for the sheaf of torsion-free differentials. From this result and by comparing the log smooth and ordinary flat deformations at first order of such surfaces, we deduce that the functors of locally trivial and flat deformations of such surfaces are unobstructed; we deduce a smoothability result for such surfaces as well. 


\vspace{0.8cm}

\section{Introduction.}

\subsection{Description and context of the article.}

In this paper, we are interested in studying the deformation theory and the logarithmic deformation theory of the varieties obtained as a gluing $X$ of a ruled surface over an elliptic curve $E$ along two (disjoint) sections along a shift. These gluings are non-normal (in particular, singular) proper (not always projective) varieties; nevertheless, by \cite[Proposition 2]{HS94}, we can suppose that the split rank two line bundle defining the ruled surface is of the form $\mathscr{O}_E \oplus \mathscr{L}$ where $\mathscr{L}$ is a line bundle \emph{of degree zero} on the elliptic curve $E$.

We prove several results about these glued surfaces, the main one is that their functor of logarithmic deformations, of locally trivial deformations and flat deformations are smooth (respectively: Theorem \ref{functor_log_def_smooth_gluing}, Corollary \ref{smooth_locally_trivial_kuranishi_space} and Corollary \ref{smooth_kuranishi_space}); we also prove that such varieties are always smoothable (Corollary \ref{smoothability_result}). Along the way, we prove that they are always d-semistable in the sense of \cite{Fri83} (Proposition \ref{result_d-semistability}), we compute the tangent space of the functor of flat deformations of $X$ (Corollary \ref{comp_coh_tangent_sheaf_gluing}) and compute (Corollary \ref{computation_tgt_space_log_smooth_deformations}) that of log smooth deformations of $X$ when the latter is equipped with a logarithmic structure of semistable type.

Such varieties have already been studied in the context of degenerations of complex abelian surfaces. More precisely, in \cite{HKW93}, the three authors construct a toroidal compactification of the moduli space of complex abelian surface equipped with a polarization of type $(1,p)$ where $p \geq 3$ is a prime number that they call the \emph{Igusa compactification}. In this book, they show that four types of points arise in the boundary of this compactified moduli space and the varieties that we are considering are part of the boundary points that correspond to type II degenerations --- in their language, they are called \textit{corank one central boundary points}. Their approach is nevertheless quite different from that of this paper since they construct the boundary points using Mumford's construction of degenerate abelian surfaces but do not study their deformation theory.

In the spirit of \cite{HS94}, where the authors study these type II degenerate abelian surfaces from a Hodge-theoretic perspective, we would like to investigate in a forthcoming paper a Torelli theorem for Deligne's mixed Hodge structure on these degenerate varieties.

\subsection{Acknowledgments and keywords.}

The author would like to thank Christian Lehn and Simon Felten for useful discussions, he was supported by the Deutsche Forschungsgemeinschaft research grant Le 3093/3-2.

%\textsc{MSC classification:} 14D06, 14B12, 14K10.

\subsection{Conventions and notations.}

In this document, a \emph{monoid} will always mean a commutative monoid and a \emph{ring} will always mean a commutative ring with unit. A \emph{variety} means a separated, reduced and locally of finite type scheme over a field. We will systematically work over the field $\mathbb{C}$ of complex numbers unless otherwise stated. The set of non-negative integers is denoted by $\mathbb{N}$. We denote by $\mathbb{C}[\varepsilon]$ the $\mathbb{C}$-algebra of dual numbers defined as $\mathbb{C}[\varepsilon]:=\mathbb{C}[t]/(t^2)$ and we denote by $\varepsilon$ the class of $t$ in that quotient.

We say that $f : \mathscr{X} \to \Delta$ is a \emph{semistable degeneration} if $\mathscr{X}$ is smooth, $X_0:=f^{-1}(0)$ is normal crossing (its irreducible components might not be smooth) and if $f$ is a proper and flat morphism that is smooth over the punctured disk $\Delta^\star$.

If $M$ is a monoid, $R$ is a ring and if $\varphi : M \to (R, \cdot)$ is a monoid homomorphism, the affine logarithmic scheme with underlying scheme $\mathrm{Spec}(R)$ and whose logarithmic structure is induced by $\varphi$ is denoted by $\mathrm{Spec}(M \xrightarrow{\varphi} R)$. We denote by $S_0$ the standard logarithmic point, see \ref{example_log_point} in the appendix for the definition.

\subsection{The glued varieties $X$.}

We start with an elliptic curve $E$ and a line bundle $\mathscr{L} \in \mathrm{Pic}(E)$ of degree zero. We consider the ruled surface $S:=\mathbb{P}_E(\mathscr{O}_E \oplus \mathscr{L})$, we have the structural morphism $p : S \to E$ which is a fiber bundle with fiber $\mathbb{P}^1$ and we have the universal rank one quotient $\mathscr{O}_S(1)$ of $\mathscr{O}_E \oplus \mathscr{L}$; the latter depends on $\mathscr{L}$ contrary to what the notation suggests but this will not cause confusion. By the universal property of the projectivization, the two surjections
\begin{center}
    $\mathscr{O}_E \oplus \mathscr{L} \twoheadrightarrow \mathscr{O}_E$ and $\mathscr{O}_E \oplus \mathscr{L} \twoheadrightarrow \mathscr{L}$
\end{center}
correspond to two sections
\begin{center}
    $\sigma_0 : E \to S$ and $\sigma_\infty : E \to S$
\end{center}
of $p$ such that
\begin{center}
    $\sigma_0^*(\mathscr{O}_S(1)|_{D_0}) \cong \mathscr{O}_E$ and $\sigma_\infty^*(\mathscr{O}_S(1)|_{D_\infty}) \cong \mathscr{L}$
\end{center}
where, for all $i \in \lbrace 0, \infty \rbrace$, we consider
\begin{center}
    $D_i:=\sigma_i(E) \subseteq S$ and $\tilde{D}:=D_0+D_\infty$;
\end{center} 
these are smooth divisors in $S$. For each $s \in E$, we define a new surface $X$ by the following cocartesian square:
\begin{center}
    \begin{tikzcd}
	{\tilde{D}} & S \\
	D & X
	\arrow["\nu", from=1-2, to=2-2]
	\arrow["\iota", from=1-1, to=1-2]
	\arrow["q"', from=1-1, to=2-1]
	\arrow["j", from=2-1, to=2-2]
\end{tikzcd}
\end{center}
where $\iota$ is the inclusion $\tilde{D} \hookrightarrow S$ and $q$ is defined by the identification of $D_0$ with $E$ and the identification of $D_\infty$ with $E$ post-composed with the translation by $s$. This surface $X$ exists by \cite[Th\'{e}or\`{e}me 5.4.]{Fer03} since $q$ is a finite morphism and since any finite set of points of $E$ or $S$ is contained in an affine open subset. The variety $X$ is a proper (local) normal crossing variety whose normalization is equal to $S$ and whose singular locus is equal to $D$. We denote by $\nu : S \to X$ the normalization morphism.

% The goal of that document is to compute the tangent space of the classical deformations of $X$ and the (log smooth) logarithmic deformation theory of $X$ when it makes sense; for that, we will be able to deduce smoothability results on $X$.

\section{Understanding $X$ from its normalization.}

In this section, we start by explaining the link between line bundles on $X$ and line bundles on its normalization and we deduce a d-semistability result for $X$. We pick a logarithmic structure of a specific type for our gluing $X$ and we state several short exact sequences relating logarithmic forms on $X$ and, ultimately, logarithmic forms on $S$ with at most logarithmic poles along $\tilde{D}$. In the last part of the section, we deduce from these short exact sequences that the variety $X$ has a trivial dualizing sheaf. 

\subsection{Line bundles.}

By \cite[V., Proposition 2.3.]{Har77}, we can formulate the link between line bundles on $X$ and line bundles on $S$.

\begin{proposition}
	A line bundle $\mathscr{E}:=p^*\mathscr{M} \otimes_{\mathscr{O}_S} \mathscr{O}_S(n)$ on $S$ descends to $X$ if and only if
	\begin{center}
		$\mathscr{M} \otimes_{\mathscr{O}_E} t_s^*\mathscr{M}^\vee \cong t_s^*\mathscr{L}^{\otimes n}$.
	\end{center}
\end{proposition}

\begin{preuve}
	By \cite[Th\'{e}or\`{e}me 2.2.]{Fer03}, the condition to descend is given by
	\begin{center}
		$\sigma_0^*\mathscr{E}|_{D_0} \otimes_{\mathscr{O}_E} t_s^*\sigma_\infty^* \mathscr{E}|_{D_\infty}^\vee \cong \mathscr{O}_E$.
	\end{center}
	Since we have $\sigma_0^*\mathscr{O}_S(1)|_{D_0} \cong \mathscr{O}_E$ and $\sigma_\infty^* \mathscr{O}_S(1)|_{D_\infty} \cong \mathscr{L}$, we can rewrite it as
	\begin{center}
		$\mathscr{M}  \otimes_{\mathscr{O}_E} t_s^*\mathscr{M}^\vee \otimes t_s^*\mathscr{L}^{\otimes -n} \cong \mathscr{O}_E$,
	\end{center}
	which can also be rewritten as
	\begin{center}
		$\mathscr{M} \otimes_{\mathscr{O}_E} t_s^* \mathscr{M}^\vee \cong t_s^*\mathscr{L}^{\otimes n}$.
	\end{center}
\end{preuve}

From this description, we can deduce d-semistability results for $X$, let us first recall this notion.

\begin{definition}
    Let $Y$ be a normal crossing complex variety. \\
    We say that $Y$ is \emph{d-semistable} if the sheaf $\mathscr{T}^1_Y:=\mathscr{E}xt^1_{\mathscr{O}_Y}(\Omega^1_Y, \mathscr{O}_Y)$ on $Y_{\mathrm{sing}}$ is isomorphic to $\mathscr{O}_{Y_{\mathrm{sing}}}$ where $Y_{\mathrm{sing}}$ is seen as a reduced subvariety of $Y$.
\end{definition}

\begin{remarque}
    By \cite[Lemma 12.1.]{FKat96}, the sheaf $\mathscr{T}^1_Y$ is in fact a line bundle on $Y_{\mathrm{sing}}$.
\end{remarque}

By \cite[Proposition 11.7.]{FKat96}, this has a translation in terms of logarithmic geometry; we refer the interested reader to the appendix concerning all notions of logarithmic geometry and especially to Definition \ref{definition_log_structure_semistable_type} for the definition of a logarithmic structure of semistable type.

\begin{proposition}
    Let $Y$ be a normal crossing complex variety. The following two assertions are equivalent:
    \begin{enumerate}
        \item $Y$ is d-semistable.
        \item $Y$ has a logarithmic structure of semistable type.
    \end{enumerate}
\end{proposition}

\begin{remarque}
    This logarithmic structure of semistable type is not unique: such logarithmic structures on $X$ are exactly the global sections of the sheaf of sets $\mathrm{LS}_X$, see Definition \ref{definition_sheaf_LS} in the appendix.
\end{remarque}

We now study the d-semistability of $X$.

\begin{proposition} \label{result_d-semistability}
    The variety $X$ is d-semistable.
\end{proposition}

\begin{preuve}
    Since the sheaf $\mathscr{T}^1_X$ is supported on the singular locus of $X$, we can work in an \'{e}tale neighborhood of it and therefore suppose that $X$ has \emph{global} normal crossings, i.e. $X$ is the union of two irreducible components $X_0$ and $X_\infty$ and $D=X_0 \cap X_\infty$. By \cite[Lemma 12.1.]{FKat96}, we have
 \begin{center}
     $\mathscr{T}^1_X \cong \sigma_0^*\mathscr{N}_{D_0/S} \otimes_{\mathscr{O}_E} t_s^*\sigma_\infty^*\mathscr{N}_{D_\infty/S}$.
 \end{center}
 By the adjunction formula, we have for all $i \in \lbrace 0 , \infty \rbrace$
 \begin{center}
     $\mathscr{N}_{D_i/S} \cong \omega_{D_i} \otimes \omega_S^\vee|_{D_i} \cong \omega_S^\vee|_{D_i}$
 \end{center}
 since $D_i$ is an elliptic curve for all $i \in \lbrace 0, \infty \rbrace$. Now, we know that
 \begin{center}
     $\omega_S \cong \omega_{S/E} \otimes p^*\omega_E \cong \omega_{S/E} \cong p^*\mathscr{L} \otimes_{\mathscr{O}_S} \mathscr{O}_S(-2)$
 \end{center}
 where the second isomorphism comes from the fact that $\omega_E \cong \mathscr{O}_E$ and the third comes from \cite[III., Exercise 8.4.]{Har77}. From that, we therefore deduce that
\begin{center}
	$\mathscr{T}^1_X \cong \sigma_0^*(p^*\mathscr{L}^\vee \otimes_{\mathscr{O}_S} \mathscr{O}_S(2)) \otimes_{\mathscr{O}_E} t_s^*\sigma_\infty^*(p^*\mathscr{L}^\vee \otimes_{\mathscr{O}_S} \mathscr{O}_S(2)) \cong \mathscr{L}^\vee \otimes t_s^*\mathscr{L}^\vee \otimes t_s^*\mathscr{L}^{\otimes 2} \cong \mathscr{L}^\vee \otimes_{\mathscr{O}_E} t_s^*\mathscr{L}$.
\end{center}
Then, we indeed get that
\begin{center}
	$\mathscr{T}^1_X=\mathscr{O}_D \Leftrightarrow t_s^*\mathscr{L}=\mathscr{L}$
\end{center}
but the right-hand side is true since $\mathscr{L}$ is a degree zero line bundle on an elliptic curve and is therefore invariant under pullback along any translation.
\end{preuve}

We now fix once and for all a logarithmic structure of semistable type $\mathscr{M}_X$ on $X$. In other words, we fix a global section of the sheaf of sets $\mathrm{LS}_X$ and by \ref{global_sections_LS}, it corresponds to a non-zero complex number that we denote by $\lambda$. For this precise logarithmic structure of semistable type on $X$, a chart for the morphism $X \to S_0$ is given \'{e}tale locally by the triple $(\mathbb{N}^2 \xrightarrow{a} \mathbb{C}[z_0,z_1,z_2]/(z_0z_1), \mathbb{N} \xrightarrow{b} \mathbb{C}, \mathbb{N} \xrightarrow{\delta} \mathbb{N}^2)$ where 
    \begin{itemize}
        \item $a$ is defined by $a(p,q):=\lambda^pz_0^pz_1^q$
        \item $b$ is defined by $b(1):=0$
        \item $\delta$ is defined by $\delta(1):=(1,1)$.
    \end{itemize}

From now on, the symbol $X$ will refer to the \emph{logarithmic scheme} $X$ equipped with the logarithmic structure of semistable type fixed above and $\underline{X}$ will refer to the underlying \emph{scheme}. When no confusion can arise (for example, writing $\mathscr{O}_X$ instead of $\mathscr{O}_{\underline{X}})$, we will drop the underline notation for the underlying scheme. Essentially, the only cases where this notation will be used are when we will deal with sheaf of differentials (logarithmic or classical) forms and objects that involve the logarithmic structure of $X$.

\subsection{Logarithmic differential forms.}

In this subsection, if $Y \to Z$ is a morphism of logarithmic schemes, the sheaf of \emph{logarithmic} $p$-forms of $Y$ over $Z$ will be denoted by $\Omega^p_{Y/Z}$ and the sheaf of \emph{ordinary} $p$-forms will be written $\Omega^p_{\underline{Y}/\underline{Z}}$, see Definition \ref{definition_log_cotangent_sheaf} in the appendix for a definition.

\subsubsection{The case of functions and 1-forms.}

\begin{proposition} \label{ses_sections_gluing}
	We have the following short exact sequence of $\mathscr{O}_X$-modules
	\begin{center}
		\begin{tikzcd}
0 \arrow[r] & \mathscr{O}_X \arrow[r, "\nu^\sharp"] & \nu_*\mathscr{O}_S \arrow[r, "\delta"] & j_*\mathscr{O}_D \arrow[r] & 0
\end{tikzcd}
	\end{center}
	where $\delta(f):=\sigma_0^*f|_{D_0} - t_{-s}^* \sigma_\infty^*f|_{D_\infty}$ for any local section $f$ of $\nu_*\mathscr{O}_S$ and $j : D \hookrightarrow X$ is the inclusion.
\end{proposition}

We continue with a similar result for logarithmic 1-forms.

\begin{proposition} \label{ses_log1-forms}
	We have an exact sequence of $\mathscr{O}_X$-modules:
	\begin{center}
	\begin{tikzcd}
0 \arrow[r] & \Omega^1_{X/S_0} \arrow[r, "\nu^*"] & \nu_*\Omega^1_{\overline{S}/S_0} \arrow[r, "\varphi"] & j_*\Omega^1_{D/S_0}  \arrow[r] & 0
\end{tikzcd}
\end{center}
where
\begin{itemize}
    \item $\overline{S}$ is the logarithmic scheme with underlying scheme $S$ and the logarithmic structure pulled back from the one on $X$
    \item $D$ is equipped with the restriction of the logarithmic structure of $X$
    \item $\varphi$ maps a section $\alpha$ of $\nu_*\Omega^1_{\overline{S}/S_0}$ to the section $\sigma_0^*\alpha|_{D_0}-t_{-s}^*\sigma_{\infty}^* \alpha|_{D_\infty}$ of $j_*\Omega^1_{D/S_0}$.
    \end{itemize}
\end{proposition}

\begin{preuve}
	If we start with a logarithmic form $\alpha \in \Omega^1_{X/S_0}$, then, because of the commutativity of the square
\begin{center}
	\begin{tikzcd}
	{D_0 \coprod D_\infty} & {\overline {S}} \\
	D & X
	\arrow["\nu", from=1-2, to=2-2]
	\arrow["j"', from=2-1, to=2-2]
	\arrow[from=1-1, to=1-2]
	\arrow["{\sigma_0 \coprod (t_s \circ \sigma_\infty)}"', from=1-1, to=2-1]
\end{tikzcd},
\end{center}
we get that the image of $\alpha$ in $j_*\Omega^1_{D/S_0}$ is equal to zero. To prove that $\Omega^1_{X/S_0}$ is indeed the kernel of $\varphi$, we can work locally for the \'{e}tale topology. In that setting, we can suppose that
\begin{center}
	$X=\mathrm{Spec}(\mathbb{C}[z_0,z_1,z_2]/(z_0z_1))$
\end{center}
	with logarithmic structure given by the logarithmic structure associated to the prelogarithmic structure given by
	\begin{center}
		$(p,q) \in \mathbb{N}^2 \mapsto \lambda^p z_0^pz_1^q \in \mathbb{C}[z_0,z_1,z_2]/(z_0z_1)$.
	\end{center}
	In that case, the sheaves of logarithmic forms can be computed explicitly, we have:
	\begin{center}
		$\Omega^1_{X/S_0} \cong \mathbb{C}[z_0,z_1,z_2]/(z_0z_1) \cdot \dfrac{\mathrm{d}z_0}{z_0} \oplus \mathbb{C}[z_0,z_1,z_2]/(z_0z_1) \cdot \mathrm{d}z_2$ \\
		$\nu_*\Omega^1_{\overline{S}/S_0} \cong \mathbb{C}[z_0, z_2] \dfrac{\mathrm{d}z_0}{z_0} \oplus  \mathbb{C}[z_0, z_2] \mathrm{d}z_2 \oplus \mathbb{C}[z_1, z_2] \dfrac{\mathrm{d}z_1}{z_1} \oplus \mathbb{C}[z_1, z_2] \mathrm{d}z_2$ \\
		$j_*\Omega^1_{D/S_0} \cong \mathbb{C}[z_2] \cdot \dfrac{\mathrm{d}z_0}{z_0} \oplus \mathbb{C}[z_2] \cdot \mathrm{d}z_2$.
	\end{center}
	We can write down the local expression of the map $\Omega^1_{X/S_0} \to \nu_*\Omega^1_{\overline{S}/S_0}$:
	\begin{center}
		$f_0 \frac{\mathrm{d}z_0}{z_0} + f_2 \mathrm{d}z_2 \mapsto f_0(z_0,0,z_2) \frac{\mathrm{d}z_0}{z_0} + f_2(z_0,0,z_2) \mathrm{d}z_2 - f_0(0,z_1,z_2+s) \frac{\mathrm{d}z_1}{z_1} + f_2(0,z_1,z_2+s) \mathrm{d}z_2$
	\end{center}
	and for the map $\nu_*\Omega^1_{\overline{S}/S_0} \to j_*\Omega^1_{D/S_0}$:
	\begin{center}
		$f_0 \frac{\mathrm{d}z_0}{z_0} + f_2 \mathrm{d}z_2 + f_1 \frac{\mathrm{d}z_1}{z_1} +g_2 \mathrm{d}z_2 \mapsto (f_0(0,z_2)+f_1(0,z_2-s)) \frac{\mathrm{d}z_0}{z_0} + (f_2(0,z_2)-g_2(0,z_2-s)) \mathrm{d}z_2$.
	\end{center}
	%The same proof as the one of the proposition \ref{ses_sections_gluing} shows the exactness of the sequence on the stalks
    Using the explicit description of the maps given above, one can prove the exactness at the level of stalks.
\end{preuve}

% \begin{remarque}
%     One can make the same remark as before here: the kernel of $\nu^*$ is supported on $\underline{D}$ and is therefore a torsion subsheaf of $\Omega^1_{X/S_0}$ but, in virtue of the proposition \ref{cotangent_sheaf_locally_free_log_smooth}, the latter is a locally free $\mathscr{O}_X$-module since $X \to S_0$ is log smooth, so the kernel of $\nu^*$ must be zero.
% \end{remarque}

\subsubsection{The case of 2-forms.}

Before introducing a similar sequence for logarithmic 2-forms, we relate the sheaf of logarithmic forms of $\overline{S}$ to the logarithmic de Rham complex of $S$ with logarithmic poles along $\tilde{D}$. We prove that in several steps using various intermediate lemmas.

We would like to construct a morphism of logarithmic schemes $h : \overline{S} \to \tilde{S}$ making the diagram
\begin{equation}\label{diagram_morph_connection_S}
    \begin{tikzcd}
	X & {\overline{S}} & {\tilde{S}} & S \\
	{S_0} & {S_0} & {S_0} & \star
	\arrow[from=1-1, to=2-1]
	\arrow[Rightarrow, no head, from=2-1, to=2-2]
	\arrow["\nu"', from=1-2, to=1-1]
	\arrow[from=1-2, to=2-2]
	\arrow["h", dashed, from=1-2, to=1-3]
	\arrow[Rightarrow, no head, from=2-2, to=2-3]
	\arrow[from=1-3, to=1-4]
	\arrow[from=1-4, to=2-4]
	\arrow[from=2-3, to=2-4]
	\arrow[from=1-3, to=2-3]
	\arrow["\lrcorner"{anchor=center, pos=0.125}, draw=none, from=1-3, to=2-4]
\end{tikzcd}
\end{equation}
commutative, where $S$ is equipped with the divisorial logarithmic structure associated to $\tilde{D}$.

We start by identifying the ghost sheaves $\overline{\mathscr{M}}_{\overline{S}}$ and $\overline{\mathscr{M}}_{\tilde{S}}$ of $\overline{S}$ and $\tilde{S}$.

\begin{lemme} \label{identification_ghost_sheaves}
    As \'{e}tale and analytic sheaves, we have $\overline{\mathscr{M}}_X \cong \underline{\mathbb{N}}_X \oplus \underline{\mathbb{N}}_D$. In particular, $\overline{\mathscr{M}}_{\overline{S}} \cong \underline{\mathbb{N}}_{\overline{S}} \oplus \underline{\mathbb{N}}_{D_0} \oplus \underline{\mathbb{N}}_{D_\infty}$.
\end{lemme}

\begin{preuve}
    We define a morphism $\varphi$ of analytic sheaves \footnote{The same formula also works for the \'{e}tale topology.} by
        \begin{align*}
\underline{\mathbb{N}}_X \oplus \underline{\mathbb{N}}_D&\rightarrow \nu_*\underline{\mathbb{N}}_{\overline{S}}\\
(f_1,f_2)&\mapsto \nu^*f_1 + \nu^*f_2
\end{align*}
We can check its injectivity and bijectivity on the stalks. Around all points of $X \setminus D$, the result is clear since $\nu$ is an isomorphism there. Let $x \in D$.
\begin{itemize}
    \item \textbf{Injectivity:} if $(f_1,f_2)$ and $(g_1, g_2)$ are such that $\nu^*f_1+\nu^*f_2=\nu^*g_1+\nu^*g_2$, then since
    \begin{center}
        $f_1|_{X \setminus D} = g_1|_{X \setminus D}$,
    \end{center}
    the fact that $f_1$ and $g_1$ are locally constant implies that $f_1 = g_1$. Now, the fact that $\mathbb{N}$ is an integral monoid in the sense of Definition \ref{definition_integral_ff_monoid}, it implies that $\nu^*f_2 = \nu^*g_2$ which implies that $f_2 = g_2$. Therefore, the morphism is injective.
    \item \textbf{Surjectivity:} if $f$ is a (germ of) section of $\nu_*\underline{\mathbb{N}}_{\overline{S}}$, we can suppose \footnote{For the \'{e}tale topology, the argument is slightly different: one can take an \'{e}tale chart of the form $\mathrm{Spec}(\mathbb{C}[z_0,z_1,z_2]/(z_0z_1)) \to X$ around $x$ and the normalization will then have two connected components that we will denote by $U_0$ and $U_\infty$ in that case.} that it is defined over a disconnected open subset $U:=U_0 \coprod U_\infty$ of $\overline{S}$ where $U_i \cap D_i \neq \emptyset$. Since $U_0 \cap U_\infty = \emptyset$, we have
    \begin{center}
        $f = f|_{U_0} + f|_{U_\infty}$.
    \end{center}
    and therefore, $f$ is the image of $(\overline{f|_{U_0}}, f|_{U_\infty})$ where $\overline{f|_{U_0}}$ is the extension of $f|_{U_0}$ to $X$ by zero --- note that it is still locally constant. 
\end{itemize}
Therefore, the morphism $\varphi$ is an isomorphism. By Lemma \ref{ses_ghost_sheaf_log_str_X}, we have $\overline{\mathscr{M}}_X = \nu_*\underline{\mathbb{N}}_{\overline{S}}$, which terminates the proof.
\end{preuve}

Now, we construct a sheaf that will encode the locally existing morphisms $h$.

\begin{paragraphe}
    We consider the sheaf of sets $\mathscr{P}$ on $S$ defined by
    \begin{center}
        $\Gamma(U, \mathscr{P}):= \lbrace \varphi \in \mathrm{Hom}_{\mathbf{Sh}(\mathrm{Mon})}(\mathscr{M}_{\tilde{S}}|_U, \mathscr{M}_{\overline{S}}|_U) \mid \varphi~\text{fits in the commutative square}~ (\ref{comm_square_sheaf_P})\rbrace$
    \end{center}
    for any $U \subseteq S$ open, where the commutative square in question is given by
    \begin{equation} \label{comm_square_sheaf_P}
        \begin{tikzcd}
	{\mathscr{M}_{\tilde{S}}|_U} & {\mathscr{M}_{\overline{S}}|_U} \\
	{\mathscr{O}_{\tilde{S}}|_U} & {\mathscr{O}_{\overline{S}}|_U}
	\arrow["\varphi", from=1-1, to=1-2]
	\arrow["{\alpha_{\tilde{S}}}"', from=1-1, to=2-1]
	\arrow[Rightarrow, no head, from=2-1, to=2-2]
	\arrow["{\alpha_{\overline{S}}}", from=1-2, to=2-2]
\end{tikzcd}
    \end{equation}
    We define an action of $\mathscr{O}_S^\times$ on that sheaf as follows: if $\varphi \in \Gamma(U, \mathscr{P})$ and $u \in \Gamma(U, \mathscr{O}_S^\times)$, we define
    \begin{center}
        $u \cdot \varphi : m \in \Gamma(V, \mathscr{M}_{\tilde{S}}) \mapsto u(\overline{m})+\varphi(m) \in \Gamma(V, \mathscr{M}_{\overline{S}})$
    \end{center}
    where $u(\overline{m})$ is defined as the image of $m$ by the composition
    \begin{center}
        $\mathscr{M}_{\tilde{S}} \to \overline{\mathscr{M}}_{\tilde{S}} \cong \underline{\mathbb{N}}_S \oplus \underline{\mathbb{N}}_{\tilde{D}} \to \underline{\mathbb{N}}_S \to \mathscr{O}_S^\times \to \mathscr{M}_{\overline{S}}$
    \end{center}
    where the morphism of sheaves of monoids $\underline{\mathbb{N}}_S \to \mathscr{O}_S^\times$ maps 1 to $u$ and where the identification comes from Lemma \ref{identification_ghost_sheaves}.
\end{paragraphe}

\begin{lemme}
    The sheaf of sets $\mathscr{P}$ is an $\mathscr{O}_S^\times$-pseudo-torsor.
\end{lemme}

\begin{preuve}
    Let us consider two sections $\varphi_1$ and $\varphi_2$ of $\mathscr{P}$ on an open set $U$. For every section $m$ of $\mathscr{M}_{\tilde{S}}$, we have
	\begin{center}
		$\overline{\varphi_1(m)}=\overline{\varphi_2(m)}$
	\end{center}
	where the bar stands for the image of a section of in the ghost sheaf. This holds because sections of $\mathscr{P}$ are uniquely determined at the level of ghost sheaves. Now, from this equality, we deduce that there exists an element $u(m)$ of $\Gamma(U, \mathscr{O}_S^\times)$ such that
	\begin{center}
		$\varphi_1(m)=u(m)+\varphi_2(m)$.
	\end{center}
	In fact, $u(m)$ is even unique. Indeed, if also $v(m)$ were such an element, we would have
	\begin{center}
		$u(m)+\varphi_2(m)=v(m)+\varphi_2(m)$
	\end{center}
	which implies that
	\begin{center}
		$u(m)=v(m)$
	\end{center}
	since $\mathscr{M}_{\overline{S}}$ is a sheaf of \emph{integral} monoids. By uniqueness, we also have that the morphism of sheaves of sets $u : m \in \mathscr{M}_{\overline{S}} \mapsto u(m) \in \mathscr{O}_S^\times$ is a morphism of sheaves of monoids. The fact that $\varphi_1$ and $\varphi_2$ coincide on $\mathscr{M}_{\tilde{S}}^\times:=\alpha^{-1}\mathscr{O}_{\tilde{S}}^\times$ implies that the morphism $u$ is constant equal to 1 on the subsheaf $\mathscr{M}_{\tilde{S}}^\times$ and therefore descends to a morphism of sheaves of monoids 
	\begin{center}
		$\overline{u} : \overline{\mathscr{M}}_{\tilde{S}} \to \mathscr{O}_S^\times$.
\end{center}
By Lemma \ref{identification_ghost_sheaves}, we have an isomorphism
\begin{center}
	$\overline{\mathscr{M}}_{\tilde{S}} \cong \underline{\mathbb{N}}_S \oplus \underline{\mathbb{N}}_{\tilde{D}}$
\end{center}
and we have $\overline{u}(p_1, p_2)=\overline{u}(p_1, 0)$ by construction of $\overline{u}$. In particular, $\overline{u}$ is a monoid homomorphism $\underline{\mathbb{N}}_S \to \mathscr{O}_S^\times$ and is completely determined by the image of 1 in $\mathscr{O}_S^\times$, hence the result.
\end{preuve}

\begin{remarque}
	If we consider the subsheaf $\mathscr{Q}$ of $\mathscr{P}$ consisting of morphisms $\varphi : \mathscr{M}_{\tilde{S}} \to \mathscr{M}_{\overline{S}}$ such that $\varphi(\rho_{\tilde{S}})=\rho_{\overline{S}}$ \footnote{Here, $\rho_{\tilde{S}} \in \Gamma(\tilde{S}, \mathscr{M}_{\tilde{S}})$ is the image of $1$ by the composition $\mathbb{N} \to \mathscr{M}_{S_0}=\mathbb{N} \oplus \mathbb{C}^\times \to \mathscr{M}_{\tilde{S}}$.}, it becomes a $\lbrace 1 \rbrace$-pseudo-torsor. Indeed, if $\varphi_1$ and $\varphi_2$ are two sections of that subsheaf, the fact that
	\begin{center}
		$\rho_{\overline{S}}=\varphi_1(\rho_{\tilde{S}})=\varphi_2(\rho_{\tilde{S}})$
	\end{center}
	implies that $u(\rho_{\tilde{S}})=1$ and thus that $\overline{u}(\overline{\rho_{\tilde{S}}})=1$. Since $\rho_{\tilde{S}}=(1,1)$ through the isomorphism
	\begin{center}
		$\overline{\mathscr{M}}_{\overline{S}} \cong \underline{\mathbb{N}}_S \oplus \underline{\mathbb{N}}_{\tilde{D}}$,
	\end{center}
	we deduce that $u=1$.
\end{remarque}

We are not always sure that $\mathscr{Q}$ is not locally empty but we have the following result.

\begin{corollaire} \label{globalization_morphism_h}
	If $h : \overline{S} \to \tilde{S}$ exists locally \footnote{By this abusive formulation, we mean that there are two open covers of $\tilde{S}$ and $\overline{S}$ such that on each open subset of the cover, a morphism making the diagram (\ref{diagram_morph_connection_S}) restricted to that open subset commutative exists.}, then it exists globally.
\end{corollaire}

\begin{preuve}
	If $h$ exists locally, then the sheaf $\mathscr{Q}$ is a torsor under the trivial sheaf of groups so the gluing relation follow from that. 
\end{preuve}

For a scheme $T$, we denote by $\eta : T_{\text{\'{e}t}} \to T_{\mathrm{Zar}}$ the canonical morphism of sites.

\begin{proposition} \label{morph_connection_S}
    There exists a morphism of logarithmic schemes $h_{\mathrm{Zar}} : (\underline{\overline{S}}, \eta_*^{\log}\mathscr{M}_{\overline{S}})  \to (\underline{\tilde{S}}, \eta_*^{\log}\mathscr{M}_{\tilde{S}})$ making the diagram (\ref{diagram_morph_connection_S}) with all the logarithmic structures restricted to the Zariski site commutative.
\end{proposition}

\begin{preuve}
Using the explicit equation of the elliptic curve, one can write down explicitly the expression of $h$ in each of the standard affine charts intersected with a trivializing open cover for the line bundle $\mathscr{L}$. Now, Corollary \ref{globalization_morphism_h} implies that the morphism of logarithmic schemes $h$ exists \emph{where the (\'{e}tale) logarithmic structures on the source and the target are restricted to the Zariski site}.
\end{preuve}

So far, we just obtained a morphism of sheaves of monoids between the Zariski logarithmic structures $\eta_*^{\log}\mathscr{M}_{\tilde{S}}$ and $\eta_*^{\log}\mathscr{M}_{\overline{S}}$. Let us see why it induces a morphism between $\mathscr{M}_{\tilde{S}}$ and $\mathscr{M}_{\overline{S}}$.

\begin{lemme}
    The morphism $h_{\mathrm{Zar}}^{\flat} : \eta_*^{\log}\mathscr{M}_{\tilde{S}} \to \eta_*^{\log}\mathscr{M}_{\overline{S}}$ of Zariski sheaves of monoids induces by functoriality of $\eta^*_{\log}$ a morphism $\mathscr{M}_{\tilde{S}} \to \mathscr{M}_{\overline{S}}$ of \'{e}tale sheaves of monoids. In particular, it implies the existence of a morphism of logarithmic schemes $h : \overline{S} \to \tilde{S}$ making the diagram (\ref{diagram_morph_connection_S}) commutative.
\end{lemme}

\begin{preuve}
    One just has to understand why
\begin{center}
    $\eta^*_{\log}\eta_*^{\log}\mathscr{M}_{\tilde{S}}\cong \mathscr{M}_{\tilde{S}}$ and $\eta^*_{\log}\eta_*^{\log}\mathscr{M}_{\overline{S}} \cong \mathscr{M}_{\overline{S}}$.
\end{center}
For $\tilde{S}$, this follows from the fact that $\mathscr{M}_{\tilde{S}}$ is the base change of the divisorial logarithmic structure associated to a smooth divisor --- in particular unibranch --- so we do not lose any information when we restrict to the Zariski topology; see \cite[III., Proposition 1.6.5.]{Og18}. For $\overline{S}$, this follows from Corollary \ref{criterion_zariski_log_str_constant_sheaf} with the Zariski sheaf $\overline{\mathscr{M}}_{\overline{S}}=\underline{\mathbb{N}}_S \oplus \underline{\mathbb{N}}_{\tilde{D}}$ and Lemma \ref{identification_ghost_sheaves}.
\end{preuve}

Let us now define a sheaf of ideals of monoids (see Definition \ref{definition_idealized_log_scheme}) on $\tilde{S}$ and $\overline{S}$ that will make $h$ behave almost like a log \'{e}tale morphism.

\begin{paragraphe}
     On $\overline{S}$ and $\tilde{S}$, we define the following sheaves of ideals of monoids, defined \'{e}tale-locally as follows:
     \begin{itemize}
         \item On $\overline{S}$, we consider the ideal generated by the image of $(1,1)$ in $\mathscr{M}_{\overline{S}}$. We denote it by $\overline{\mathscr{I}}$.
         \item On $\tilde{S}$, we consider the ideal generated by the image of $(0,1)$ in $\mathscr{M}_{\tilde{S}}$. We denote it by $\tilde{\mathscr{I}}$.
     \end{itemize}
     Endowed with these, $(\tilde{S}, \tilde{\mathscr{I}})$ and $(\overline{S}, \overline{\mathscr{I}})$ have the structure of idealized logarithmic schemes since in each case, the image of the sheaf of ideals in the structure sheaf is zero as the diagram (\ref{diagram_morph_connection_S}) shows. Moreover, $h$ is now a morphism is idealized logarithmic schemes since the morphism of monoids $\overline{\mathscr{M}}_{\tilde{S}} \to \overline{\mathscr{M}}_{\overline{S}}$ maps $(0,1)$ to $(1,1)$.
\end{paragraphe}

\begin{proposition} \label{h_ideally_log_etale}
    The morphism $h : (\overline{S}, \overline{\mathscr{I}}) \to (\tilde{S}, \tilde{\mathscr{I}})$ is ideally log \'{e}tale.
\end{proposition}

\begin{preuve}
    The property of being ideally log \'{e}tale is local for the \'{e}tale topology, so we can work locally for the \'{e}tale topology. By the remark above, $h$ is given locally by the commutative square of monoids
    \begin{center}
        \begin{tikzcd}
	{\mathbb{C}[x,z]} & {\mathbb{C}[x,z]} \\
	{\mathbb{N}^2} & {\mathbb{N}^2}
	\arrow[Rightarrow, no head, from=1-1, to=1-2]
	\arrow[from=2-1, to=1-1]
	\arrow[from=2-2, to=1-2]
	\arrow["a", from=2-2, to=2-1]
\end{tikzcd}
    \end{center}
    where $a$ is the (injective) morphism defined by $a(0,1):=(1,1)$ and $a(1,0):=(1,0)$. It induces a ring homomorphism $\mathbb{Z}[\tilde{x}, \tilde{y}] \to \mathbb{Z}[\overline{x}, \overline{y}]$ given by 
    \begin{center}
        $\tilde{x} \mapsto \overline{x}$ and $\tilde{y} \mapsto \overline{x} \cdot \overline{y}$.
    \end{center}
    The way the variable $\tilde{y}$ is mapped creates in some sense a "normal crossing singularity", so that $h$ is not log \'{e}tale. Nevertheless, we can make this singularity disappear using ideals of monoids (logarithmic blow-up), we denote by $\tilde{I}$ the ideal of the monoid $\mathbb{N}^2$ generated by $(0,1)$ and by $\overline{I}$ the ideal of the monoid $\mathbb{N}^2$ generated by $(1,1)$. Since $a(0,1)=(1,1)$, the morphism $a$ induces a ring homomorphism 
    \begin{center}
        $\mathbb{Z}[\mathbb{N}^2, \tilde{I}] \to \mathbb{Z}[\mathbb{N}^2, \overline{I}]$
    \end{center}
    which is just the identity. In particular, the morphism of schemes
    \begin{center}
        $\underline{\overline{S}} \to \underline{\tilde{S}} \times_{\mathrm{Spec}(\mathbb{Z}[\mathbb{N}^2, \tilde{I}])} \mathrm{Spec}(\mathbb{Z}[\mathbb{N}^2, \overline{I}])$
    \end{center}
    is \'{e}tale --- in fact, it is the identity. Theorem \ref{idealized_toroidal_criterion} implies that $h$ is ideally log \'{e}tale.
\end{preuve}

\begin{corollaire} \label{identification_forms_derham_pullback}
    We have an isomorphism of $\mathscr{O}_S$-modules $\Omega^1_{\overline{S}/S_0} \cong \Omega^1_{S/\star}=\Omega^1_S(\log \tilde{D})$.
\end{corollaire}

\begin{preuve}
    Since $h$ is ideally log \'{e}tale by Proposition \ref{h_ideally_log_etale}, the pullback morphism
    \begin{center}
        $h^* : \Omega^1_{\tilde{S}/S_0} \to h_*\Omega^1_{\overline{S}/S_0}$
    \end{center}
    is an isomorphism. In addition, since $\tilde{S}=S \times_\star S_0$, \cite[IV., Proposition 1.2.15.]{Og18} implies that
    \begin{center}
        $\Omega^1_{\tilde{S}/S_0} \cong \mathrm{pr}_1^*\Omega^1_{S/\star}$
    \end{center}
    but the scheme morphism underlying $\mathrm{pr}_1 : \tilde{S} \to S$ is the identity, which proves the result.
\end{preuve}

\begin{paragraphe} \label{definition_phi2}
    We have an injective morphism of $\mathscr{O}_X$-modules
    \begin{center}
        $\nu^* : \Omega^2_{X/S_0} \to \nu_*\Omega^2_{\overline{S}/S_0}$.
    \end{center}
    In addition, we consider the map $\varphi_2$ defined by
    \begin{center}
        $\varphi_2 : \alpha \in \Gamma(U, \nu_*\Omega^2_{\overline{S}/S_0}) \mapsto \sigma_0^*\mathrm{Res}_{D_0}(\alpha)+t_{-s}^*\sigma_\infty^*\mathrm{Res}_{D_\infty}(\alpha) \in \Gamma(U, j_*\Omega^1_{\underline{D}})$
    \end{center}
    through the identification $\Omega^2_{\overline{S}/S_0} \cong \Omega^2_S(\log \tilde{D})$ given by Corollary \ref{identification_forms_derham_pullback} and where $j : D \hookrightarrow X$ is the inclusion. We will give another construction that does not use this identification right after the proof of the proposition.
\end{paragraphe}

\begin{proposition} \label{ses_log_2_forms}
    The short sequence of $\mathscr{O}_X$-modules
    \begin{center}
        \begin{tikzcd}
	0 & {\Omega^2_{X/S_0}} & {\nu_*\Omega^2_{\overline{S}/S_0}} & {j_*\Omega^1_{\underline{D}}} & 0
	\arrow[from=1-1, to=1-2]
	\arrow["{\nu^*}", from=1-2, to=1-3]
	\arrow["{\varphi_2}", from=1-3, to=1-4]
	\arrow[from=1-4, to=1-5]
\end{tikzcd}   
    \end{center}
    is exact.
\end{proposition}

% \begin{preuve}
%     This can be checked locally in the \'{e}tale topology, therefore we can assume that
% 	\begin{center}
% 		$X=\mathrm{Spec}(\mathbb{C}[z_0,z_1,z_2]/(z_0z_1))$
% 	\end{center}
% 	as in the case of logarithmic 1-forms. We can also describe the sheaves of logarithmic forms:
% 	\begin{center}
% 		$\Omega^2_{X/S_0}= \mathbb{C}[z_0,z_1,z_2]/(z_0z_1) \cdot \dfrac{\mathrm{d}z_0}{z_0}\wedge \mathrm{d}z_2$ \\
% 		$\nu_*\Omega^2_{\bar{S}/S_0}=\mathbb{C}[z_0,z_2] \cdot \dfrac{\mathrm{d}z_0}{z_0} \wedge \mathrm{d}z_2 \oplus \mathbb{C}[z_1,z_2] \cdot \dfrac{\mathrm{d}z_1}{z_1} \wedge \mathrm{d}z_2$\\
% 		$j_*\Omega^1_{\underline{D}} \cong \mathbb{C}[z_2] \cdot \mathrm{d}z_2$
% 	\end{center}
% 	and we can describe the maps in the sequence (3) \'{e}tale-locally. The first one is given by
% 	\begin{center}
% 		$f(z_0,z_1,z_2) \frac{\mathrm{d}z_0}{z_0} \wedge \mathrm{d}z_2 \mapsto f(z_0,0,z_2) \frac{\mathrm{d}z_0}{z_0} \wedge \mathrm{d}z_2 - f(0,z_1,z_2+s) \frac{\mathrm{d}z_1}{z_1} \wedge \mathrm{d}z_2$
% 	\end{center}
% 	and the second one is given by
% 	\begin{center}
% 		$f_0(z_0,z_2) \frac{\mathrm{d}z_0}{z_0} \wedge \mathrm{d}z_2 + f_1(z_1,z_2) \frac{\mathrm{d}z_1}{z_1} \wedge \mathrm{d}z_2 \mapsto (f_0(0,z_2)+f_1(0,z_2-s))\mathrm{d}z_2$.
% 	\end{center}
% 	This proves that the sequence is exact.
% \end{preuve}

We now give the other construction, that does not use the identification $\Omega^2_{\overline{S}/S_0} \cong \Omega^2_S(\log \tilde{D})$ 
in the construction of the cokernel of the pullback of logarithmic 2-forms along $\nu$.

\begin{paragraphe} \label{constr_log_residus_2forms_gluing}
    We consider the logarithmic schemes $F_0$ and $F_\infty$ whose underlying schemes are respectively $\underline{D_0}$ and $\underline{D_\infty}$ such that the morphisms of log schemes $F_0 \to S_0$ and $F_\infty \to S_0$ are strict in the sense of Definition \ref{definition_strict_morphism}. The composition
    \begin{center}
        \begin{tikzcd}
	{D_0} & {\overline{S}} & {S_0}
	\arrow["{j_0}", from=1-1, to=1-2]
	\arrow[from=1-2, to=1-3]
        \end{tikzcd}
    \end{center}
yields a morphism of sheaves
\begin{center}
    \begin{tikzcd}
	{\Omega^2_{\overline{S}/S_0}} & {(j_0)_*\Omega^2_{D_0/S_0}}
	\arrow[from=1-1, to=1-2]
\end{tikzcd}
\end{center}
since $j_0$ is a strict morphism.
\par On the other hand, the morphism of log schemes $D_0 \to S_0$ can be factored as
\begin{center}
    \begin{tikzcd}
	{D_0} & {F_0} & {S_0}
	\arrow[from=1-1, to=1-2]
	\arrow[from=1-2, to=1-3]
\end{tikzcd}
\end{center}
where $D_0 \to F_0$ is the identity on the underlying schemes and $F_0 \to S_0$ is strict. This composition induces a short exact sequence of sheaves
\begin{equation} \label{ses_log_2_forms_intr}
    \begin{tikzcd}
	{\Omega^1_{F_0/S_0}} & {\Omega^1_{D_0/S_0}} & {\Omega^1_{D_0/F_0}} & 0
	\arrow[from=1-1, to=1-2]
	\arrow[from=1-2, to=1-3]
	\arrow[from=1-3, to=1-4]
\end{tikzcd}
\end{equation}
Moreover, we have a similar factorization and short exact sequence of sheaves for $D_\infty$ and $F_\infty$.
\end{paragraphe}

\begin{lemme} \label{lemma_two_forms_D_i_over_S_0}
    The sequence (\ref{ses_log_2_forms_intr}) is in fact exact on the left. Moreover, we have an isomorphism
    \begin{center}
        $\Omega^2_{D_0/S_0} \cong \Omega^1_{\underline{D_0}} \otimes_{\mathscr{O}_{D_0}} \Omega^1_{D_0/F_0}$
    \end{center}
    of $\mathscr{O}_{D_0}$-modules.
\end{lemme}

\begin{preuve}
    The left exactness can be checked \'{e}tale locally. Locally for the \'{e}tale topology, we have
    \begin{center}
        $\Omega^1_{F_0/S_0} \cong \mathbb{C}[z_2] \cdot \mathrm{d}z_2$ and $\Omega^1_{D_0/S_0} \cong \mathbb{C}[z_2] \cdot \mathrm{d}z_2 \oplus \mathbb{C}[z_2] \cdot \frac{\mathrm{d}z_0}{z_0}$
    \end{center}
    and the map on the left of the sequence is given by
    \begin{center}
        $f(z_2)\cdot \mathrm{d}z_2 \mapsto (f(z_2) \cdot \mathrm{d}z_2, 0)$,
    \end{center}
    which proves the injectivity. The local description of the first two sheaves of the sequence already shows that
    \begin{itemize}
        \item $\Omega^1_{F_0/S_0}$ is locally free of rank 1
        \item $\Omega^1_{D_0/S_0}$ is locally free of rank 2.
    \end{itemize}
    Now, we can also describe the cokernel locally, it is given by
    \begin{center}
        $\Omega^1_{D_0/F_0} \cong \mathbb{C}[z_2] \cdot \frac{\mathrm{d}z_0}{z_0} \oplus \mathbb{C}[z_2] \cdot \frac{\mathrm{d}z_1}{z_1}/(\frac{\mathrm{d}z_0}{z_0} + \frac{\mathrm{d}z_1}{z_1}) \cong \mathbb{C}[z_2] \cdot \frac{\mathrm{d}z_0}{z_0}$
    \end{center}
    and it is therefore locally free of rank one. Since the three sheaves are locally free, we have
    \begin{center}
        $\Omega^2_{D_0/S_0}=\mathrm{det}(\Omega^1_{D_0/S_0}) \cong \Omega^1_{F_0/S_0} \otimes_{\mathscr{O}_{D_0}} \Omega^1_{D_0/F_0}$
    \end{center}
    and since $F_0 \to S_0$ is strict, we have $\Omega^1_{F_0/S_0}=\Omega^1_{\underline{D_0}}$ which implies the result.
\end{preuve}

\begin{remarque}
    The same proof also leads to an isomorphism of $\mathscr{O}_{D_\infty}$-modules $\Omega^2_{D_\infty/S_0} \cong \Omega^1_{\underline{D_\infty}} \otimes_{\mathscr{O}_{D_\infty}} \Omega^1_{D_\infty/F_\infty}$.
\end{remarque}

\begin{lemme} \label{triviality_log_forms_D0_over_F0}
    For $i \in \lbrace 0, \infty \rbrace$, we have $\Omega^1_{D_i/F_i} \cong \mathscr{O}_{D_i}$. In particular, we have $\Omega^2_{D_i/S_0} \cong \Omega^1_{\underline{D_i}}$.
\end{lemme}

\begin{preuve}
    Let $i \in \lbrace 0, \infty \rbrace$. Since the morphism $F_i \to S_i$ is strict, we have $\overline{\mathscr{M}}_{F_i} \cong \underline{\mathbb{N}}_{F_i}$. Moreover, we have an isomorphism $\overline{\mathscr{M}}_{D_i} \cong \underline{\mathbb{N}}^2_{D_i}$ by Lemma \ref{identification_ghost_sheaves} so it implies that the $\mathscr{O}_{D_i}$-module $\Omega^1_{D_i/F_i}$ is \emph{globally} generated by $\frac{\mathrm{d}z_0}{z_0}$ if $i=0$ and $\frac{\mathrm{d}z_1}{z_1}$ if $i=\infty$.  
\end{preuve}

Let us now prove that these a priori two different constructions lead to the same result through the two 
identifications $\Omega^2_{\overline{S}/S_0} \cong \Omega^2_S(\log \tilde{D})$ and $\Omega^2_{D_0/S_0} \cong \Omega^1_{\underline{D_0}}$. We start with a lemma about the sheaf of 2-forms of $\overline{S}$ over $S_0$.

\pagebreak

\begin{lemme} \label{lemme_commutative_square_2forms}
    The square of $\mathscr{O}_X$-modules
    \begin{center}
        \begin{tikzcd}
	{\nu_*\Omega^2_S(\log \tilde{D})} && {j_*\Omega^1_{\underline{D}}} \\
	{\nu_*\Omega^2_{\overline{S}/S_0}} & {(\nu \circ j_0)_* \Omega^1_{\underline{D_0}} \oplus (\nu \circ j_\infty)_* \Omega^1_{\underline{D_\infty}}} & {j_*\Omega^1_{\underline{D}}}
	\arrow["{\sigma_0^*\mathrm{Res}_{D_0} + t_{-s}^*\sigma_\infty^*\mathrm{Res}_{D_\infty}}", from=1-1, to=1-3]
	\arrow["\sim", from=1-1, to=2-1]
	\arrow[equals, from=1-3, to=2-3]
	\arrow[from=2-1, to=2-2]
	\arrow["\sim"', from=2-2, to=2-3]
\end{tikzcd}
    \end{center}
    is commutative.
\end{lemme}

\begin{preuve}
    Let us make the bottom horizontal arrow explicit: restriction of logarithmic forms give two morphisms of $\mathscr{O}_S$-modules
    \begin{center}
        $\Omega^2_{\overline{S}/S_0} \to (j_0)_* \Omega^2_{D_0/S_0}$ and $\Omega^2_{\overline{S}/S_0} \to (j_\infty)_* \Omega^2_{D_\infty/S_0}$
    \end{center}
    where $j_i : D_i \hookrightarrow \overline{S}$ is the (strict) restriction; we therefore get a morphism of $\mathscr{O}_S$-modules
    \begin{center}
        $\varphi : \Omega^2_{\overline{S}/S_0} \to (j_0)_* \Omega^2_{D_0/S_0} \oplus (j_\infty)_* \Omega^2_{D_\infty/S_0}$.
    \end{center}
    %For all $i \in \lbrace 0, \infty \rbrace$, we have an isomorphism of $\mathscr{O}_{D_i}$-modules
    %\begin{equation}\label{residu_partpol_partreg}
        %\begin{tikzcd}
	%{\Omega^1_{F_i/S_0} \otimes_{\mathscr{O}_{D_i}} \Omega^1_{D_i/F_i}} & {\Omega^2_{D_i/S_0}}
	%\arrow["\sim", from=1-1, to=1-2]
        %\end{tikzcd}
    %\end{equation}
    %given on stalks by wedging a section of $\Omega^1_{F_i/S_0}$ and any preimage in $\Omega^1_{D_i/S_0}$ of a section of $\Omega^1_{D_i/F_i}$ together. We have seen that 
    %\begin{center}
       % $\Omega^1_{D_0/F_0} \cong \mathscr{O}_{D_0} \langle \mathrm{dlog}(z_0) \rangle$
    %\end{center}
    %and that
    %\begin{center}
        %$\Omega^1_{D_\infty/F_\infty} \cong \mathscr{O}_{D_\infty} \langle \mathrm{dlog}(z_1) \rangle$
   %\end{center}
    Combining the identifications given by Lemma \ref{lemma_two_forms_D_i_over_S_0} and Lemma \ref{triviality_log_forms_D0_over_F0}, we get a morphism of $\mathscr{O}_X$-modules
    \begin{center}
        \begin{tikzcd}
	{\nu_*\Omega^2_{\overline{S}/S_0}} & {(\nu \circ j_0)_*\Omega^2_{D_0/S_0} \oplus (\nu \circ j_\infty)_*\Omega^2_{D_\infty/S_0}} & {(\nu \circ j_0)_*\Omega^1_{\underline{D_0}} \oplus (\nu \circ j_\infty)_*\Omega^1_{\underline{D_\infty}}}
	\arrow["\varphi", from=1-1, to=1-2]
	\arrow["\sim"', from=1-2, to=1-3]
\end{tikzcd}
    \end{center}
    which is the first morphism of the horizontal composite in the square of the lemma. Note that, via the identification $\Omega^2_{\overline{S}/S_0} \cong \Omega^2_S(\log \tilde{D})$, this composition is given by taking the residue along $D_0$ and $D_\infty$ respectively. The second map of the horizontal composite is given by
    \begin{center}
        $(\alpha_0, \alpha_\infty) \in \Gamma(U, (\nu \circ j_0)_*\Omega^1_{\underline{D_0}} \oplus (\nu \circ j_\infty)_*\Omega^1_{\underline{D_\infty}}) \mapsto \sigma_0^*\alpha_0+t_{-s}^*\sigma_\infty^*\alpha_\infty \in \Gamma(U, j_*\Omega^1_{\underline{D}})$.
    \end{center}

    This implies that the square is commutative. %Now, we have to understand the effect of the left vertical arrow. This is given by the pullback along $h : \overline{S} \to \tilde{S}$ and since the morphism of sheaves of monoids
    %\begin{center}
        %$\overline{h}^\flat : \overline{\mathscr{M}}_{\tilde{S}} \cong \underline{\mathbb{N}}_S \oplus \underline{\mathbb{N}}_{\tilde{D}} \to \overline{\mathscr{M}}_{\overline{S}} \cong \underline{\mathbb{N}}_S \oplus \underline{\mathbb{N}}_{\tilde{D}}$
    %\end{center}
    %maps $(0,1)$ to $(0,1)$, the divisorial part of $h$ is preserved which implies that the square of the lemma is indeed commutative.
\end{preuve}

\begin{lemme} \label{trivial_log_dr_sheaf_2forms}
    We have $\Omega^2_{\overline{S}/S_0} \cong \mathscr{O}_S$.
\end{lemme}

\begin{preuve}
    By Corollary \ref{identification_forms_derham_pullback}, we have an isomorphism $\Omega^2_{\overline{S}/S_0} \cong \Omega^2_S(\log \tilde{D})$. Moreover, since $\Omega^2_S(\log \tilde{D})$ is a line bundle on $S$, it is of the form
    \begin{center}
        $\Omega^2_S(\log \tilde{D})=:p^*\mathscr{M} \otimes_{\mathscr{O}_S} \mathscr{O}_S(n)$
    \end{center}
    where $\mathscr{M}$ is a line bundle on $E$ and $n \in \mathbb{Z}$. Now, if $F \cong \mathbb{P}^1_\mathbb{C}$ is a fiber of $p$, we have
    \begin{center}
        $\Omega^2_S(\log \tilde{D})|_F \cong \Omega^1_{\mathbb{P}^1}(2) \cong \mathscr{O}_F$
    \end{center}
    which implies that $n=0$. Now, since the tensor product is right-exact, the Poincar\'{e} residue restricted to $D_0$ gives a surjection of $\mathscr{O}_{D_0}$-modules
    \begin{center}
        $\Omega^2_S(\log \tilde{D}) \otimes_{\mathscr{O}_S} \mathscr{O}_{D_0} \twoheadrightarrow \mathscr{O}_{D_0}$
    \end{center}
    but the source and the target of this surjection are line bundles on $D_0$; therefore, this surjection is an isomorphism. Now, this implies that
    \begin{center}
        $\mathscr{M} \cong \sigma_0^*(\Omega^2_S(\log \tilde{D})|_{D_0}) \cong \mathscr{O}_E$
    \end{center}
    and therefore, $\Omega_S^2(\log \tilde{D}) \cong \mathscr{O}_S$.
\end{preuve}

\subsection{Consequences for the dualizing sheaf of $X$.}

The short exact sequence for logarithmic 2-forms constructed above has a very nice consequence on $X$. Note that by Proposition \ref{dualizing_sheaf_d_semistable variety}, we have $\omega_X=\Omega^2_{X/S_0}$.

\begin{proposition} \label{triv_dualizing_sheaf_gluing}
	The variety $X$ is log Calabi-Yau, i.e. $\omega_X \cong \mathscr{O}_X$.
\end{proposition}
	
\begin{preuve}
	Combining Proposition \ref{ses_sections_gluing}, Proposition \ref{ses_log_2_forms} and Lemma \ref{lemme_commutative_square_2forms}, we have the following (a priori not commutative) diagram:
	\begin{center}
		\begin{tikzcd}
0 \arrow[r] & \mathscr{O}_X \arrow[r, "\nu^*"] \arrow[d] & \nu_*\mathscr{O}_S \arrow[r, "\delta"] \arrow[d, "\sim"] & \mathscr{O}_{D} \arrow[r] \arrow[d, "\sim"] & 0 \\
0 \arrow[r] & \omega_X \arrow[r, "\nu^*"]                      & \nu_*\Omega^2_{\overline{S}/S_0} \arrow[r, "\varphi_2"]          & \Omega^1_{\underline{D}} \arrow[r]                               & 0
\end{tikzcd}
	\end{center}
	where the vertical middle arrow is the isomorphism of Lemma \ref{trivial_log_dr_sheaf_2forms}. We now want to prove that the right square of the diagram above is commutative so that we can conclude that
 \begin{center}
    $\omega_X \cong \mathscr{O}_X$.    
 \end{center}
 We have the short exact sequence of $\mathscr{O}_S$-modules
 \begin{center}
     $0 \to p^* \Omega^1_E \to \Omega^1_{\underline{S}} \to \Omega^1_{\underline{S}/E} \to 0$
 \end{center}
 and taking the determinant yields an isomorphism
 \begin{center}
     $\Omega^2_{\underline{S}} \cong p^*\Omega^1_E \otimes_{\mathscr{O}_S} \Omega^1_{\underline{S}/E}$.
 \end{center}
Now, the map
\begin{center}
    \begin{tabular}{ccc}
      $\Omega^2_{\underline{S}} \otimes_{\mathscr{O}_S} \mathscr{O}_S(\tilde{D})$ & $\to$ & $\Omega^2_S(\log \tilde{D})$  \\
         $\alpha \otimes f$ & $\mapsto$ & $f\alpha$  
    \end{tabular}
\end{center}
is well-defined because $f$ has at most simple poles along $\tilde{D}$ so $f\alpha$ has at most simple poles along $\tilde{D}$ and $\mathrm{d}(f\alpha)$ as well since it is zero for dimension reasons; so $f\alpha$ is indeed a section of $\Omega^2_S(\log \tilde{D})$. Moreover, one can check locally that this morphism is an isomorphism. This yields isomorphisms:
 \begin{center}
     $\Omega^2_S(\log \tilde{D}) \cong \Omega^2_{\underline{S}} \otimes_{\mathscr{O}_S} \mathscr{O}_S(\tilde{D}) \cong p^*\Omega^1_E \otimes_{\mathscr{O}_S} \Omega^1_{\underline{S}/E} \otimes_{\mathscr{O}_S} \mathscr{O}_S(\tilde{D})$.
 \end{center}
 so by Lemma \ref{trivial_log_dr_sheaf_2forms} and Corollary \ref{identification_forms_derham_pullback}, we have
 \begin{center}
     $\Omega^1_{\underline{S}/E} \otimes_{\mathscr{O}_S} \mathscr{O}_S(\tilde{D}) \cong \mathscr{O}_S$.    
 \end{center}
since $\Omega^1_E \cong \mathscr{O}_E$. Now, if we call $\alpha$ a global generator of $\Omega^1_{\underline{S}/E} \otimes_{\mathscr{O}_S} \mathscr{O}_S(\tilde{D})$ and $\beta$ a global generator of $\Omega^1_E$, the images of elements in the right square of the diagram that starts the proof are given as follows:
\begin{center}
    \begin{tikzcd}
	f & {\sigma_0^*f|_{D_0}-t_{-s}^*\sigma_\infty^*f|_{D_\infty}} \\
	{f(p^*\beta \wedge \alpha)} & \gamma
	\arrow[maps to, from=1-1, to=1-2]
	\arrow[maps to, from=1-1, to=2-1]
	\arrow[maps to, from=2-1, to=2-2]
\end{tikzcd}
\end{center}
where $\gamma$ is given by
\begin{center}
    $\gamma=\sigma_0^*(f|_{D_0}p^*\beta|_{D_0} \cdot \mathrm{Res}_{D_0}(\alpha)) + t_{-s}^*\sigma_\infty^*(f|_{D_\infty}p^*\beta|_{D_\infty} \cdot \mathrm{Res}_{D_\infty}(\alpha))$
\end{center}
 and where $s$ is the gluing shift defining $X$. Since on a fiber $F \cong \mathbb{P}^1$, $\alpha$ is the global generator of $\Omega^1_{\mathbb{P}^1}(2)$ given by $\frac{\mathrm{d}x}{x}$ whose sum of residues at 0 and at $\infty$ is zero, we have
\begin{center}
    $\sigma_0^*\mathrm{Res}_{D_0}(\alpha)=-t_{-s}^*\sigma^*_\infty \mathrm{Res}_{D_\infty}(\alpha)$.
\end{center}
Note that, since $\alpha$ exists globally as a section of $\Omega^1_{\underline{S}/E} \otimes_{\mathscr{O}_S} \mathscr{O}_S(\tilde{D})$, its residue also exists globally so it must be a constant since it is a section of $\mathscr{O}_{D_\infty}$ and in particular, it is invariant under translation. Since  $\beta$ is also translation invariant by \cite[Proposition 5.1.]{Sil09}, we can rewrite $\gamma$ as 
\begin{center}
    $\gamma=(\sigma_0^*f|_{D_0}-t_{-s}^*\sigma_\infty^*f|_{D_\infty}) \beta$
\end{center}
as sections of $\Omega^1_E$. This is the image of the section 
    \begin{center}
        $\sigma_0^*f|_{D_0}-t_{-s}^*\sigma_\infty^*f|_{D_\infty}$
    \end{center}
    of $\mathscr{O}_E$ in $\Omega^1_E$; in other words, the square on the right commutes and $\omega_X \cong \mathscr{O}_X$. 
\end{preuve}

From the fact that $X$ has a trivial dualizing sheaf, we are able to compute the cohomology of its structure sheaf.

\begin{corollaire}\label{computation_coho_str_sheaf}
	We have $h^i(X, \mathscr{O}_X)= \binom{2}{i}$ for all $i \geq 0$.
\end{corollaire}

\begin{preuve}
	As $X$ is a compact connected complex variety, we have that $h^0(X, \mathscr{O}_X)=1$ and since $\omega_X \cong \mathscr{O}_X$ by Proposition \ref{triv_dualizing_sheaf_gluing}, Serre duality for coherent sheaves implies that $h^2(X, \mathscr{O}_X)=1$. Now, the short exact sequence of Proposition \ref{ses_sections_gluing} implies that 
 \begin{center}
     $\chi(X, \mathscr{O}_X)=\chi(S, \mathscr{O}_S)-\chi(D, \mathscr{O}_D)=0$,
 \end{center}
 since $\chi(S, \mathscr{O}_S)=0$ by \cite[V., Corollary 2.5.]{Har77}. This equality implies that $h^1(X,\mathscr{O}_X)=2$ as we wanted to prove.

\end{preuve}

\section{Hodge theory of $X$.}

\subsection{Computation of Deligne's mixed Hodge structure of $X$.}

From \cite{Del74}, we know that the cohomology group $\mathrm{H}^k(X, \mathbb{Q})$ with rational coefficients of any complex algebraic variety $X$ carries a mixed Hodge structure. This is done using a semi-simplicial resolution of cohomological descent of $X$. In this section, we will try to compute this mixed Hodge structure explicitly, using the following semi-simplicial resolution of $X$:
\begin{center}
	\begin{tikzcd}
	X & {S=:X_0} & {D=:X_1}
	\arrow["\nu"', from=1-2, to=1-1]
	\arrow["{\sigma_0}"', shift right=1, from=1-3, to=1-2]
	\arrow["{\sigma_\infty \circ t_s}", shift left=1, from=1-3, to=1-2]
\end{tikzcd}
\end{center}

We start by computing the second page of the weight spectral sequence for $X$.

\begin{paragraphe}
	By \cite[Scholie 8.1.9.]{Del74}, we have a spectral sequence (called the \emph{weight spectral sequence}) whose first page is given by:
	\begin{center}
		$\mathrm{E}^{p,q}_1:= \mathrm{H}^q(X_p, \mathbb{Q}) \Rightarrow \mathrm{H}^{p+q}(X, \mathbb{Q})$
	\end{center}
\end{paragraphe}

Now, let us compute the differentials in the first page of the weight spectral sequence.

\begin{proposition}
	The differentials $d_i : \mathrm{H}^i(X_0, \mathbb{Q}) \to \mathrm{H}^i(X_1, \mathbb{Q})$ are zero for all $i \geq 0$. In particular the weight spectral sequence for $X$ degenerates on the first page.
\end{proposition}

\begin{preuve}
	Let us first start with the cases $i \in \llbracket 0,2 \rrbracket$. We can work with complex coefficients by the universal coefficient theorem. For all $i$, we have the diagram
            \begin{center}
                 \begin{tikzcd}
	               {\mathrm{H}^i(E, \mathbb{Q})} & {\mathrm{H}^i(S, \mathbb{Q})} & {\mathrm{H}^i(E, \mathbb{Q})}
	               \arrow["{p^*}", from=1-1, to=1-2]
	               \arrow["{\sigma_0^*}", shift left=1, from=1-2, to=1-3]
	               \arrow["{\sigma_\infty^*}"', shift right=1, from=1-2, to=1-3]
                \end{tikzcd}
		\end{center}
    Since $\sigma_0^*$ is a section of $p^*$, we have $\sigma_0^\star \circ p^*=(p \circ \sigma_0)^*=\mathrm{id}$, so $p^*$ is injective.
    \begin{itemize}
        \item \textbf{If $i=0$}, the linear map $t_s^*$ acts as the identity on $\mathrm{H}^0(E, \mathbb{C})$ since these cohomology classes are represented by constant functions by properness of $E$. Since $i=0$, the map $p^*$ is an injective linear map between two vector spaces of the same dimension so it is an isomorphism; in particular, it has at most one left inverse. Since $\sigma_0^*$ and $\sigma_\infty^*$ are two left inverses of $p^*$, we deduce that $\sigma_0^*=\sigma_\infty^*=(\sigma_\infty \circ t_s)^*$: in other words: $d_0=0$.
        \item \textbf{If $i=1$}, we also want to prove that $t_s^*$ acts as the identity on $\mathrm{H}^1(E, \mathbb{C})$. Since $E$ is a projective curve, the Hodge decomposition yields an isomorphism
  \begin{center}
      $\mathrm{H}^1(E, \mathbb{C}) \cong \mathrm{H}^0(E, \Omega^1_E) \oplus \overline{\mathrm{H}^0(E, \Omega^1_E)}$
  \end{center}
  and any global section of $\Omega^1_E$ is of the form $\alpha \cdot \beta$ where $\beta$ is a global generator and $\alpha \in \mathbb{C}$. Since $t_s^*\beta=\beta$ by \cite[proposition 5.1.]{Sil09}, it implies that for any class $\omega \in \mathrm{H}^1(E, \mathbb{C})$, we have $t_s^*\omega=\omega$. By the same argument as above, $p^*$ is an isomorphism and we deduce that $d_1=0$.
        \item \textbf{If $i=2$}, the linear map $t_s^*$ acts as the identity on $\mathrm{H}^2(E, \mathbb{C})$ since we have
        \begin{center}
            $\mathrm{H}^2(E, \mathbb{C}) \cong \mathrm{H}^1(E, \Omega^1_E) \cong \mathrm{H}^1(E, \mathscr{O}_E) \cong \overline{\mathrm{H}^0(E, \Omega^1_E)}$.
        \end{center}
        The same arguments as in the case where $i=1$ imply that $t_s^*$ acts as the identity on $\mathrm{H}^2(E, \mathbb{C})$. To prove that $\sigma_0^*=\sigma_\infty^*$, we cannot use literally the same argument as in the two cases before as $p^*$ is not surjective but since the image of $p^*$ is exactly the source of the $\sigma_i^*$'s, the argument still works, which proves that $d_2=0$.
    \end{itemize}
    The differentials are obviously zero for $i \geq 3$ because their target is zero.  
\end{preuve}

We can now compute the Betti numbers of $X$.

\begin{corollaire} \label{computation_weight_filtration}
	% We have
	% \begin{itemize}
	% 	\item $\mathrm{Gr}_0^W\mathrm{H}^0(X, \mathbb{Q})=\mathrm{H}^0(S, \mathbb{Q})$.
	% 	\item $\mathrm{Gr}_0^W\mathrm{H}^1(X, \mathbb{Q})=\mathrm{H}^0(E, \mathbb{Q})$ and $\mathrm{Gr}_1^W\mathrm{H}^1(X, \mathbb{Q})=\mathrm{H}^1(S, \mathbb{Q})$.
	% 	\item $\mathrm{Gr}_1^W\mathrm{H}^2(X, \mathbb{Q})=\mathrm{H}^1(E, \mathbb{Q})$ and $\mathrm{Gr}_2^W\mathrm{H}^2(X, \mathbb{Q})=\mathrm{H}^2(S, \mathbb{Q})$.
	% 	\item $\mathrm{Gr}_2^W\mathrm{H}^3(X, \mathbb{Q})=\mathrm{H}^2(E, \mathbb{Q})$ and $\mathrm{Gr}_3^W\mathrm{H}^3(X, \mathbb{Q})=\mathrm{H}^3(S, \mathbb{Q})$.
	% 	\item $\mathrm{Gr}_4^W\mathrm{H}^4(X, \mathbb{Q})=\mathrm{H}^4(S, \mathbb{Q})$.
	% \end{itemize}
	The Betti numbers of $X$ are given by
	\begin{center}
		$b_0(X)=1$, $b_1(X)=3$, $b_2(X)=4$, $b_3(X)=3$ and $b_4(X)=1$. 
	\end{center}
\end{corollaire}

\begin{remarque}
	Note that the weights are consistent: first, $h^{p,q}(\mathrm{H}^k(X, \mathbb{C}))=0$ when $p+q>k$ which was predictable since $X$ is proper (see, for example, \cite[Th\'{e}or\`{e}me 8.2.4., (iii)]{Del74}); second, $\mathrm{Gr}^W_p\mathrm{H}^k(X, \mathbb{C})$ has a \emph{pure} Hodge structure of weight $p$.	
\end{remarque}

We can now easily compute the dimension of the non-zero pieces of the Hodge filtration.

\begin{corollaire} \label{comp_hodge_fil_gluing}
We have
\begin{itemize}
	\item $F^0\mathrm{H}^1(X, \mathbb{C}) \cong \mathbb{C}^2$, $F^1 \mathrm{H}^1(X, \mathbb{C}) \cong \mathbb{C}$.
	\item $F^0\mathrm{H}^2(X, \mathbb{C}) \cong \mathbb{C}^4$, $F^1 \mathrm{H}^2(X, \mathbb{C}) \cong \mathbb{C}^3$ and $F^2\mathrm{H}^2(X, \mathbb{C})=0$.
	\item $F^0 \mathrm{H}^3(X, \mathbb{C}) \cong \mathbb{C}^3$, $F^1 \mathrm{H}^3(X, \mathbb{C}) \cong \mathbb{C}^3$ and $F^2 \mathrm{H}^3(X, \mathbb{C}) \cong \mathbb{C}$.
\end{itemize}
\end{corollaire}

\subsection{The torsion-free Hodge--de Rham spectral sequence.}

\subsubsection{The sheaf of torsion-free differentials.}

Let us start this subsection by briefly discussing the sheaf of torsion-free differentials.

\begin{definition}
    Let $Y$ be a complex algebraic variety and let $p \geq 0$.
    \begin{itemize}
        \item A \emph{torsion $p$-form} on $Y$ is a local section of $\Omega^p_Y$ that is supported on the singular locus of $Y$. We denote by $\tau^p_Y$ the subsheaf of $\Omega^p_Y$ formed of torsion $p$-forms.
        \item The sheaf $\tilde{\Omega}^p_Y$ of \emph{torsion-free p-forms} on $Y$ is defined by $\tilde{\Omega}^p_Y:=\Omega^p_Y/\tau^p_Y$.
    \end{itemize}
\end{definition}

% \begin{remarque}
%     In fact, the subsheaf $\tau^p_Y$ of $\Omega^p_Y$ is its torsion subsheaf in the following sense: if $\mathscr{F}$ is a coherent sheaf on $Y$, its \emph{torsion subsheaf} $\mathscr{F}_{\mathrm{tors}}$ is the subsheaf of $\mathscr{F}$ defined as
%     \begin{center}
%         $\mathscr{F}_{\mathrm{tors}}:=\mathrm{Ker}(\mathscr{F} \to \mathscr{F}^{\vee \vee})$.
%     \end{center}
%     If $\mathscr{F}_{\mathrm{tors}}=0$, we say that $\mathscr{F}$ is a \emph{torsion-free} $\mathscr{O}_Y$-module. The sheaf $\mathscr{F}/\mathscr{F}_{\mathrm{tors}}$ is torsion-free by construction. In particular, $\tilde{\Omega}^p_Y$ is a torsion-free $\mathscr{O}_Y$-module.
% \end{remarque}

We finish this subsubsection by computing the sheaf of torsion 1-forms on $X$ explicitly.

\begin{lemme} \label{triviality_torsion_1forms_gluing}
	We have $\tau_X^1=\mathscr{O}_D$.
\end{lemme}

\begin{preuve}
	Since the normalization $\nu : S \to X$ is an isomorphism outside the singular locus, it implies that the support of $\tau_X^1$ is contained in $D$. Now by \cite[p. 88]{Fri83}, d-semistability and Serre duality, we have
	\begin{center}
		$\mathrm{H}^0(X, \tau_X^1)^*=\mathrm{Ext}^2_{\mathscr{O}_X}(\tau^1_X, \mathscr{O}_X)=\mathrm{H}^1(D, \mathscr{O}_D)$ ;
	\end{center}
	therefore, $\tau^1_X$ has a non-zero global section. On the other hand, $\tau^1_X$ is a line bundle on $D$; to prove this claim, we can work \'{e}tale locally and therefore suppose that
    \begin{center}
        $\underline{X}=\mathrm{Spec}(\mathbb{C}[x,y,z]/(xy))$ and $\underline{S}=\mathrm{Spec}(\mathbb{C}[x,z] \times \mathbb{C}[y,z])$.
    \end{center}
    We can now compute the modules of K\"{a}hler differentials
    \begin{center}
        $\Omega^1_{\underline{X}}= \big( \mathbb{C}[x,y,z]/(xy) \cdot \mathrm{d}x \oplus \mathbb{C}[x,y,z]/(xy) \cdot \mathrm{d}y \oplus \mathbb{C}[x,y,z]/(xy) \cdot \mathrm{d}z \big)/(x\mathrm{d}y+y \mathrm{d}x)$
    \end{center}
    and the pullback map of these to the normalization
    \begin{center}
        $\nu^* : [\overline{f_0} \mathrm{d}x + \overline{f_1} \mathrm{d}y + \overline{f_2} \mathrm{d}z] \mapsto (f_0(x,0,z) \mathrm{d}x + f_2(x,0,z) \mathrm{d}z, f_1(0,y,z+s) \mathrm{d}y+f_2(0,y,z+s) \mathrm{d}z)$.
    \end{center}
    Therefore, the image of $[\overline{f_0} \mathrm{d}x + \overline{f_1} \mathrm{d}y + \overline{f_2} \mathrm{d}z]$ is zero if and only if
    \begin{center}
        $y \mid f_0$, $x \mid f_1$, $x \mid f_2$ and $y \mid f_2$
    \end{center}
    which can be rewritten 
    \begin{center}
        $y \mid f_0$, $x \mid f_1$ and $\overline{f_2}=0$.
    \end{center}
    In particular, we have
    \begin{center}
        $\tau^1_X=\mathrm{Ker}(\nu^*)=\lbrace \overline{f_0} y\mathrm{d}x + \overline{f_1}x \mathrm{d}y \in \Omega^1_{\underline{X}} \rbrace= \lbrace \overline{f_1-f_0} \cdot x \mathrm{d}y \in \Omega^1_{\underline{X}} \rbrace$,
    \end{center}
    which proves that $\tau^1_X$ is a line bundle on $D$, locally generated by the 1-form $x\mathrm{d}y$.
 
    Now since $\nu^*\tau_X^1=0$, it has degree zero as a line bundle on $D$. Since it has also a non-zero global section, it must be trivial.
\end{preuve}

\subsubsection{The torsion-free Hodge--de Rham spectral sequence.}

By \cite[Proposition 1.5.]{Fri83}, we know that if $Y$ is a \emph{global} normal crossing complex variety and if each of its component are K\"{a}hler manifolds then, there is the \emph{torsion-free Hodge--de Rham} spectral sequence
\begin{center}
	$\mathrm{E}^{p,q}_1=\mathrm{H}^q(Y, \tilde{\Omega}^p_Y) \Longrightarrow \mathbb{H}^{p+q}(X, \tilde{\Omega}^\bullet_Y) \cong \mathrm{H}^{p+q}(Y, \mathbb{C})$
\end{center}
which degenerates on the first page and whose induced filtration (on the complex cohomology of $Y$) is the Hodge filtration of the Deligne mixed Hodge structure on $Y$.

In this subsubsection, we prove this result for $X$, which is only \emph{local} normal crossing. Before stating the theorem, let us set up some notations.

\begin{paragraphe}
    We define $Y_0:=X$ and $Y_1:=D$ and write $Y^{[i]}$ for the normalization of $Y_i$ and $\nu_i : Y^{[i]} \to Y_i$ for the normalization morphism. We therefore have the double complex $\mathscr{F}^{\bullet, \bullet}$ of $\mathscr{O}_X$-modules defined by
    \begin{center}
        $\mathscr{F}^{p,q}:=(\nu_p)_*\Omega^q_{Y^{[p]}}$
    \end{center}
    for $p \in \lbrace 0,1 \rbrace$ and $q \geq 0$ with differentials $d : \mathscr{F}^{p,q} \to \mathscr{F}^{p,q+1}$ being the exterior derivative and $\delta : \mathscr{F}^{0,q} \to \mathscr{F}^{1,q}$ defined by $\delta:=\sigma_0^*-(\sigma_\infty \circ t_{-s})^*$. We now consider $\mathscr{F}^\bullet$ to be the total complex associated to this double complex, it is equipped with the differential $D:=d+\delta$. As in \cite[Lemma 4.6.]{GS73}, one can show, using Poincar\'{e} lemma and the spectral sequence of a double complex, that $(\mathscr{F}^\bullet, D)$ is a resolution of the constant sheaf $\underline{\mathbb{C}}_X$.
\end{paragraphe}

We can now phrase the theorem of this subsection.

\begin{theoreme} \label{hodge_derham_mod_tors_degen}
    The torsion-free Hodge--de Rham spectral sequence given by
    \begin{center}
        $\mathrm{E}_1^{p,q}=\mathrm{H}^q(X, \tilde{\Omega}^p_X) \Longrightarrow \mathrm{H}^{p+q}(X, \mathbb{C})$
    \end{center}
    degenerates on the first page.
\end{theoreme}

\begin{preuve}
	For all $q \geq 0$, we have a short exact \footnote{This can be checked locally and in that case, we can suppose that $X$ has \emph{global} normal crossings so it follows from \cite{Fri83}.} sequence of $\mathscr{O}_X$-modules
	\begin{center}
		\begin{tikzcd}
	0 & {\tilde{\Omega}^q_X} & {\mathscr{F}^{0,q}=\nu_*\Omega^q_{\underline{S}}} & {\mathscr{F}^{1,q}=j_*\Omega^q_{\underline{D}}} & 0
	\arrow[from=1-1, to=1-2]
	\arrow["{\nu^*}", from=1-2, to=1-3]
	\arrow["{\delta}", from=1-3, to=1-4]
	\arrow[from=1-4, to=1-5]
\end{tikzcd}
	\end{center}
	Taking hypercohomology yields a spectral sequence 
	\begin{center}
		$\mathrm{F}_1^{r,s}=\mathrm{H}^s(X^{[r]}, \Omega^q_{X^{[r]}}) \Longrightarrow \mathbb{H}^{r+s}(X, \mathscr{F}^{\bullet,q}) \cong \mathrm{H}^{r+s}(X, \tilde{\Omega}^q_X)$
	\end{center}
	since $(\mathscr{F}^{\bullet,q}, \delta)$ is a resolution of the sheaf $\tilde{\Omega}^q_X$. By the same argument as in the degeneration of the weight spectral sequence on the first page, we get that the two potentially non-zero differentials on the first page of the spectral sequence $(\mathrm{F}_d^{r,s})_{d,r,s}$ are zero, which exactly means that this spectral sequence degenerates on the first page \footnote{More generally, since the differentials on the first page of the weight spectral sequence for Deligne's mixed Hodge structure on the singular cohomology of $X$ are morphisms of Hodge structures, they will restrict to the differentials on the first page of the hypercohomology spectral sequence $(\mathrm{F}_d^{r,s})_{d,r,s}$ for all $q$. In particular, the degeneration of the weight spectral sequence on the second page implies that the direct sum spectral sequence $\oplus_q \mathrm{F}_d^{r,s}$ degenerates on the second page.}. Therefore, for all $q \geq 0$, we have the equality
	\begin{center}
		$\displaystyle \dim_\mathbb{C} \Big( \bigoplus_{r,s \geq 0} \mathrm{H}^s(X^{[r]}, \Omega^q_{X^{[r]}}) \Big) = \dim_\mathbb{C}(\mathrm{H}^\bullet(X, \tilde{\Omega}^q_X))$
	\end{center}
	and taking the direct sum over $q \geq 0$ gives the equality
	\begin{center}
		$\displaystyle \dim_\mathbb{C} \Big( \bigoplus_{q,r,s \geq 0} \mathrm{H}^s(X^{[r]}, \Omega^q_{X^{[r]}}) \Big) = \dim_\mathbb{C}\Big( \bigoplus_{q=0}^{+ \infty}\mathrm{H}^\bullet(X, \tilde{\Omega}^q_X)\Big)$
	\end{center}
	Now, since the $X^{[r]}$ are compact K\"{a}hler, we have that for all $r \geq 0$,
	\begin{center}
		$\displaystyle \bigoplus_{q,s \geq 0} \mathrm{H}^s(X^{[r]}, \Omega^q_{X^{[r]}})=\mathrm{H}^\bullet(X^{[r]}, \mathbb{C})$.
	\end{center}
	and since, in our case, the weight spectral sequence does not only degenerate on the second page but on the \emph{first} page, we also have:
	\begin{equation}
		\displaystyle \bigoplus_{q,r \geq 0} \mathrm{H}^q(X^{[r]}, \mathbb{C}) = \mathrm{H}^\bullet(X, \mathbb{C}).
	\end{equation}
	Now, if the spectral sequence
	\begin{center}
		$\mathrm{E}_1^{p,q}=\mathrm{H}^q(X, \tilde{\Omega}^p_X) \Longrightarrow \mathrm{H}^{p+q}(X, \mathbb{C})$
	\end{center}
	did not degenerate on the first page, there would exist a non-zero differential on the second page and therefore
	\begin{center}
		$\displaystyle \dim_\mathbb{C} \Big( \bigoplus_{q=0}^{+ \infty} \mathrm{H}^\bullet(X, \tilde{\Omega}^q_X) \Big) > \dim_{\mathbb{C}}(\mathrm{H}^\bullet(X, \mathbb{C}))$.
	\end{center}
	Putting everything together, we would obtain
	\begin{center}
		$\displaystyle \dim_\mathbb{C}(\mathrm{H}^\bullet(X, \mathbb{C}))=\dim_\mathbb{C} \Big( \bigoplus_{r=0}^{+ \infty} \mathrm{H}^\bullet(X^{[r]}, \mathbb{C}) \Big)=\dim_\mathbb{C} \Big( \bigoplus_{q,r,s \geq 0} \mathrm{H}^s(X^{[r]}, \Omega^q_{X^{[r]}}) \Big) > \dim_\mathbb{C}(\mathrm{H}^\bullet(X, \mathbb{C}))$,
	\end{center}
	which is absurd. Thus, the torsion-free Hodge--de Rham spectral sequence must degenerate on the first page. 	 
\end{preuve}

% \begin{remarque}
% 	Here, the only important points in the degeneration of the spectral sequence $(\mathrm{F}_k^{r,s})_{k,r,s}$ (or some related spectral sequence like direct sum) and the degeneration of the weight spectral sequence on the \emph{first} page; in our case, the first fact relied on the fact that the semi-simplicial resolution of $X$ was particularly short (of length strictly lower than 2). The fact that the (ordinary) Hodge-de Rham spectral sequence for the $X^{[r]}$'s degenerated on the first page was also used (here, it happened because they were compact K\"{a}hler but we can relax the hypothesis by just asking that they belong to the class $\mathscr{C}$ of Fujiki, for example).
% \end{remarque}

\section{Deformation theories of $X$.}

In this section, we want to understand the tangent spaces of the functor of log-smooth deformations of $X$ and that of the (classical) flat deformations of $\underline{X}$; we would like to understand the link between these two theories as well. We start by introducing the deformation functors we will use. We refer to \cite{Ser06} for more background on deformation theory.

\begin{paragraphe}
    If $k$ is a field, we denote by $\mathbf{Art}_k$ the category of artinian local $k$-algebras with residue field $k$. If $\Lambda \in \mathbf{Art}_k$, we denote by $\mathbf{Art}_{\Lambda}$ the category of artinian local $\Lambda$-algebras with residue field equal to $k$. We consider the functor of Artin rings $\mathrm{Def}^{\mathrm{lt}}_{X, \Lambda}$ of locally trivial deformations of $\underline{X}$ defined by
    \begin{center}
        $\mathrm{Def}^{\mathrm{lt}}_{X, \Lambda} : A \in \mathbf{Art}_{\Lambda} \mapsto \lbrace \mathscr{X} \to \mathrm{Spec}(A)~\text{locally trivial deformation of}~X \rbrace \in \mathbf{Set}$.
    \end{center}
    Similarly, we define the functor of Artin rings $\mathrm{Def}_{X, \Lambda}$ of deformation of $\underline{X}$. We now define the functor $\mathrm{LD}_X$ of log smooth deformations of $X$ as
    \begin{center}
        $\mathrm{LD}_X : A \in \mathbf{Art}_{\mathbb{C}\llbracket t_1, \dots, t_m \rrbracket}  \mapsto \lbrace \mathscr{X} \to S_A~\text{log smooth whose special fiber is}~X \to S_0 \rbrace \in \mathbf{Set}$
    \end{center}
    where $m$ is the number of connected components of the singular locus of $X$ and $S_A$ is the logarithmic scheme whose underlying scheme is $\mathrm{Spec}(A)$ and whose logarithmic structure is given by
        \begin{align*}
\mathbb{N}^m \oplus A^\times&\rightarrow A\\
(d_1, \dots, d_m, u) &\mapsto t_1^{d_1}\dots t_m^{d_m} u
\end{align*}
We have a natural transformation $F : \mathrm{LD}_X \to \mathrm{Def}_{X, \mathbb{C}\llbracket t_1, \dots, t_m \rrbracket}$ by forgetting the logarithmic structure. Moreover, if we identify the subcategory of $\mathbf{Art}_{\mathbb{C}\llbracket t_1, \dots, t_m \rrbracket}$ of algebras such that all the $t_i$'s act as zero with the category $\mathbf{Art}_{\mathbb{C}}$, this natural transformation restricts to a natural transformation $F|_{\mathbf{Art}_{\mathbb{C}}} : \mathrm{LD}_X|_{\mathbf{Art}_{\mathbb{C}}} \to \mathrm{Def}_{X, \mathbb{C}}^{\mathrm{lt}}$.
\end{paragraphe}

\begin{remarque}
    For the sake of brevity and when no confusion can arise, we will usually omit the mention of $\Lambda$ in the functor of flat deformations when $\Lambda=\mathbb{C}$.
\end{remarque}

\subsection{Locally trivial flat deformations of $X$.}

To get more information about the classical flat deformations of $\underline{X}$, we compute the cohomology of its tangent sheaf.

\begin{proposition} \label{comp_coh_tangent_sheaf_gluing}
	We have $h^0(X, \mathscr{T}_{\underline{X}})=2$, $h^1(X, \mathscr{T}_{\underline{X}})=3$ and $h^2(X, \mathscr{T}_{\underline{X}})=1$.
\end{proposition}

\begin{preuve}
	For all $i \geq 0$, by \cite[Proof of Lemma 2.10.]{Fri83}  we have
	\begin{center}
		$\mathrm{H}^i(X, \mathscr{T}_{\underline{X}})=\mathrm{Ext}^i_{\mathscr{O}_X}(\tilde{\Omega}^1_X, \mathscr{O}_X)$,
	\end{center}
	which implies, via Serre duality, that
	\begin{center}
		$\mathrm{H}^i(X, \mathscr{T}_{\underline{X}})=\mathrm{H}^{2-i}(X, \tilde{\Omega}^1_X)$
	\end{center}
	since $\omega_X=\mathscr{O}_X$ by Proposition \ref{triv_dualizing_sheaf_gluing}. The degeneration of the torsion-free Hodge--de Rham spectral sequence on the first page implies that, for all $i \in \llbracket 0, 2 \rrbracket$, we have
	\begin{center}
		$\mathrm{H}^i(X, \mathscr{T}_{\underline{X}})=\mathrm{Gr}_1^F \mathrm{H}^{2-i}(X, \mathbb{C})$ 
	\end{center}
	Therefore, using Corollary \ref{comp_hodge_fil_gluing}, we get
	\begin{center}
		$h^0(X,\mathscr{T}_{\underline{X}})=2$, $h^1(X, \mathscr{T}_{\underline{X}})=3$ and $h^2(X, \mathscr{T}_{\underline{X}})=1$.
	\end{center}
\end{preuve}

% \begin{remarque}
% 	One can also check that the numbers $h^i(X, \mathscr{O}_X)$ coincide with the one computed respectively in Lemma \ref{computation_coho_str_sheaf}.
% \end{remarque}

\subsection{Flat deformations of $X$.}

Now that we have understood the tangent space of the locally trivial flat deformations of $\underline{X}$, we want to understand it for the flat deformations. In virtue of Serre duality for coherent sheaves, it is enough to compute the cohomology with coefficients in $\Omega^1_{\underline{X}}$. 

\subsubsection{Another gluing to construct $X$.}

In order to compute the cohomology groups of $\Omega^1_{\underline{X}}$, we introduce another construction of $X$ by performing some identifications on a fiber bundle that will be called $X_0$ --- this variety has nothing to do with the variety $X_0$ introduced in the third section. This will have the major advantage that the resulting quotient morphism $X_0 \twoheadrightarrow X$ will be \emph{\'{e}tale} and we will thus be able to descend differentials from $X_0$ to $X$ using \'{e}tale descent --- not that this was not possible with the initial gluing construction of $X$ (starting from $S$) because the normalization morphism $S \to X$ is \emph{not} \'{e}tale. Thanks to this new construction of $X$, we will be able to construct a non-torsion 1-form on $X$ and therefore compute $h^0(X, \Omega^1_{\underline{X}})$, which will give us the dimension of the space of flat deformations of $\underline{X}$.

\begin{lemme}
	As line bundles over $E$, $S \setminus D_0$ and $S \setminus D_\infty$ are respectively isomorphic to $\mathrm{Tot}(\mathscr{L}^\vee)$ and $\mathrm{Tot}(\mathscr{L})$. 
\end{lemme}

\begin{preuve}
	Let $i \in \lbrace 0, \infty \rbrace$. Since $S$ is a fiber bundle with fiber $\mathbb{P}^1$ over $E$ and $D_i$ is a section of $S$, it's clear that $S \setminus D_i$ is a fiber bundle with fiber $\mathbb{P}^1 \setminus \lbrace i \rbrace \cong \mathbb{C}$. Let us treat the case $i=\infty$ first. The transition functions for $\mathscr{O}_E \oplus \mathscr{L}$ are of the form
	\begin{align*}
(U \cap V) \times \mathbb{C}^2&\rightarrow (U \cap V) \times \mathbb{C}^2\\
(p,(x,y))&\mapsto (p, (x, \lambda_{U,V}(p)y))
\end{align*}	 
for some $\lambda_{U,V}(p) \in \mathbb{C}^\times$.  Therefore, the transition functions of the fiber bundle $S \to E$ are given by
\begin{align*}
(U \cap V) \times \mathbb{P}^1&\rightarrow (U \cap V) \times \mathbb{P}^1\\
(p,[x:y])&\mapsto (p, [x:\lambda_{U,V}(p)y])
\end{align*}
and those of the fiber bundle $S \setminus D_\infty \to E$ are given by
\begin{center}
\begin{tikzcd}
	{(p,[x:y])} & {(p,[x:\lambda_{U,V}(p)y])} \\
	{(p, \frac{y}{x})} & {(p, \frac{\lambda_{U,V}(p)y}{x})}
	\arrow[maps to, from=1-1, to=1-2]
	\arrow["\sim", maps to, from=1-2, to=2-2]
	\arrow["\sim"', maps to, from=1-1, to=2-1]
	\arrow[maps to, from=2-1, to=2-2]
\end{tikzcd}
\end{center}
where the vertical arrows are the identifications $(U \cap V) \times \mathbb{P}^1 \setminus \lbrace \infty \rbrace \cong (U \cap V) \times \mathbb{C}$. We see that these transition functions are exactly the transition functions of the line bundle $\mathscr{L} \otimes \mathscr{O}_E^\vee \cong \mathscr{L}$. 

For the case $i=0$, we get the following commutative square:
\begin{center}
	\begin{tikzcd}
	{(p,[x:y])} & {(p,[x:\lambda_{U,V}(p)y])} \\
	{(p, \frac{x}{y})} & {(p, \frac{x}{\lambda_{U,V}(p)y})}
	\arrow[maps to, from=1-1, to=1-2]
	\arrow["\sim", maps to, from=1-2, to=2-2]
	\arrow["\sim"', maps to, from=1-1, to=2-1]
	\arrow[maps to, from=2-1, to=2-2]
\end{tikzcd}
\end{center}
and we see that the transition functions of $S \setminus D_0 \to E$ are exactly those of the line bundle $\mathscr{O}_E \otimes \mathscr{L}^\vee \cong \mathscr{L}^\vee$. 
\end{preuve}

This lemma implies that we have two open immersions
\begin{center}
	$j_0 : L:=\mathrm{Tot}(\mathscr{L}) \hookrightarrow S$ and $j_\infty : L^\vee:=\mathrm{Tot}(\mathscr{L}^\vee) \hookrightarrow S$.
\end{center}

We now construct the fiber bundle $X_0 \to E$.

\begin{paragraphe}
Let us now consider the fiber bundle $X_0$ over $E$ defined as
\begin{center}
	$X_0:=\mathrm{Spec}_E \Big( \mathrm{Sym}_{\mathscr{O}_E}^\bullet(\mathscr{L}^\vee \oplus \mathscr{L})/(\mathscr{L} \cdot \mathscr{L}^\vee) \Big)$
\end{center}
and let us denote by $\pi : X_0 \to E$ the structural morphism of this fiber bundle. For $x \in E$, we have
\begin{center}
	$\pi^{-1}(x) \cong Z:=V(xy) \subseteq \mathbb{A}^2_{\kappa(x)}$.
\end{center}
The space $X_0$ is not normal and its normalization is equal to $L^\vee \coprod L$.
%We define $X$ as $X:=S/\sim$ where $\sim$ is the equivalence relation generated by
%\begin{center}
	%$\sigma_0(x) \sim \sigma_\infty(x)+s$
%\end{center}
%for all $x \in E$. It's a non-normal scheme with normalization $\nu : S \to X$. 
We have a diagram
\begin{center}
	\begin{tikzcd}
	{L^\vee \coprod L} & S \\
	{X_0} & X
	\arrow["{\tilde{j_0} \coprod j_\infty}"{pos=0.5}, from=1-1, to=1-2]
	\arrow["\nu", from=1-2, to=2-2]
	\arrow["{\nu_0}"', from=1-1, to=2-1]
\end{tikzcd}
\end{center}
where $\nu_0$ is the normalization morphism for $X_0$ given by the natural morphisms
\begin{center}
	$L^\vee \to X_0$ and $L \to X_0$
\end{center}
and where $\tilde{j_0}$ is the morphism $L^\vee \hookrightarrow S$ which is defined on a trivializing open $U \subseteq E$ for $\mathscr{L}$ by
\begin{center}
    $\tilde{j_0}$ : \begin{tabular}{ccc}
        $U \times \mathbb{C}$ & $\to$ & $U \times \mathbb{P}^1$  \\
        $(e,x)$ & $\mapsto$ & $(e+s, [x:1])$.
    \end{tabular}
\end{center}
One can check that it is \emph{globally} well-defined since it is compatible with the transition functions of the fiber bundles $X_0$ on the source and $S$ on the target since $\mathscr{L}$ has degree zero.

% Now, the two injections
% \begin{center}
% 	$\mathscr{L} \hookrightarrow \mathscr{L} \oplus \mathscr{L}^\vee$ and $\mathscr{L}^\vee \hookrightarrow \mathscr{L} \oplus \mathscr{L}^\vee$
% \end{center}
% define two injective morphisms
% \begin{center}
% 	$\mathrm{Tot}(\mathscr{L}) \hookrightarrow X_0$ and $\mathrm{Tot}(\mathscr{L}^\vee) \hookrightarrow X_0$.
% \end{center}
% We denote by $Z_0$ resp. $Z_\infty$ their image, they are the two irreducible components of $X_0$. Moreover, for both $\mathscr{L}$ and $\mathscr{L}^\vee$, we have the zero sections that furnishes sections
% \begin{center}
% 	$s_0 : E \to Z_0$ and $s_\infty : E \to Z_\infty$
% \end{center}
% to the respective structural morphisms of $\mathscr{L}$ and $\mathscr{L}^\vee$ and hence to $\pi : X_0 \to E$. By construction of $X_0$, the sections $s_0(E)$ and $s_\infty(E)$ get identified together --- they are also identified with the zero section of $\mathrm{Tot}(\mathscr{L} \oplus \mathscr{L}^\vee)$.
\end{paragraphe}

We now construct $X$ as a gluing of $X_0$. The idea is the following: a projective nodal curve is obtained by identifying the points $0$ and $\infty$ on $\mathbb{P}^1$, we would like to do this relatively over $E$ and taking the line bundle $\mathscr{L}$ in account.

\begin{paragraphe}\label{equivalence-relation-gluing-x0}
    Let $U \subseteq E$ be a trivializing open for $\mathscr{L}$. On $U \times Z$, we define a binary relation $\mathscr{R}$ as follows:
    \begin{center}
        $(e,(x,0)) \mathscr{R} (e+s, (0, x^{-1}))$ and $(e, (0,y)) \mathscr{R} (e-s, (y^{-1}, 0))$
    \end{center}
    for all $(e, (x,y)) \in U \times Z$ such that $e+s \in U$ and $x$ and $y$ are not both zero. 
\end{paragraphe}

% \begin{paragraphe}\label{equivalence-relation-gluing-x0}
%     We construct a binary relation $\mathscr{R}$ on $X_0$ as follows: if $U \subseteq E$ is a trivializing open subset for $\mathscr{L}$, we define, for $x \neq 0$ and $p \in U$,
%     \begin{center}
%         $(p,(x,0)) \mathscr{R} (p,(0, x^{-1}))$.
%     \end{center}
%     For $\mathscr{R}$ to be symmetric, we add
%     \begin{center}
%         $(p, (0, y)) \mathscr{R} (p, (y^{-1},0))$.
%     \end{center}
%     Out of the total spaces of $\mathscr{L}$ and $\mathscr{L}^\vee$, which are affine lines relatively over $E$, this identification constructs the $\mathbb{P}^1$ relatively over $E$ and the fact that we glued together the zero sections of $\mathscr{L}$ and $\mathscr{L}^\vee$ implies that the points $0$ and the points $\infty$, that are respectively relatively incarnated by the zero section of $\mathscr{L}$ and the zero section of $\mathscr{L}^\vee$, are identified together. Now, we have to take the shift $s$ in account. To do that, we add the identification
%     \begin{center}
%         $(p, (0,0)) \mathscr{R} (p+ns, (0,0))$
%     \end{center}
%     for all $n \in \mathbb{Z}$ such that $p+ns \in U$. By construction, $\mathscr{R}$ is a symmetric and transitive relation on $U \times V(xy)$, we take the reflexive closure of it, that we still denote with the same letter.
% \end{paragraphe}

\begin{lemme} \label{globalization_equiv_relation}
    The binary relation $\mathscr{R}$ defined in the paragraph \ref{equivalence-relation-gluing-x0} globalizes to binary relation on $X_0$, that we still denote by $\mathscr{R}$.
\end{lemme}

\begin{preuve}
     Let us take two trivializing open subsets $U$ and $V$ of $E$ for $\mathscr{L}$, we have the transition function
    \begin{center}
        $\lambda_{U,V} : U \cap V \to \mathbb{C}^\times$
    \end{center}
    for $L$. The transition function for $L^\vee$ is given by
    \begin{center}
        $\lambda^{-1}_{U,V} : p \in U \cap V \mapsto \lambda_{U,V}(p)^{-1} \in \mathbb{C}^\times$.
    \end{center}
    This implies that the transition functions $\varphi_{U,V} : \pi^{-1}(U \cap V) \times Z \to \pi^{-1}(U \cap V) \times Z$ for the fiber bundle $X_0$ are given by
    \begin{center}
        \begin{tabular}{ccc}
            $(e,(x,0))$ & $\mapsto$ & $(e, (\lambda_{U,V}(e)^{-1}x, 0))$  \\
            $(e,(0,y))$ & $\mapsto$ & $(e, (0,\lambda_{U,V}(e)y))$
        \end{tabular}
    \end{center}
    We must now check that the identifications given by $\mathscr{\mathscr{R}}$ are compatible with the transition functions of $X_0$ in the following sense: if $(e,z) \mathscr{R} (e',z')$, then $\varphi_{U,V}(e,z) \mathscr{R} \varphi_{U,V}(e',z')$ for all elements $(e,z)$ and $(e',z')$ of $(U \cap V) \times Z$.

    Let us start with the first type of identifications: if $(e, (x,0)) \in (U \cap V) \times Z$, then, on one side, $(e,(x,0))$ gets identified with $(e+s, (0, x^{-1}))$ which is mapped to the element $(e+s,(0, \lambda_{U,V}(e+s)x^{-1}))$ by $\varphi_{U,V}$. On the other side, $(e, (x,0))$ is mapped to $(e, (\lambda_{U,V}(e)^{-1}x, 0))$ by $\varphi_{U,V}$ which gets identified with $(e+s, (0, \lambda_{U,V}(e)x^{-1})$.

    Now, since $t_s^*\mathscr{L} \cong \mathscr{L}$ (because $\mathscr{L}$ is of degree zero), we have $\lambda_{U,V}(e)=\lambda_{U,V}(e+s)$ for all $e \in U \cap V$, which implies that
    \begin{center}
        $(e+s,(0, \lambda_{U,V}(e+s)x^{-1}))=(e+s, (0, \lambda_{U,V}(e)x^{-1})$
    \end{center}
    and similarly for points of the form $(e, (0,y))$, where we use that $t^*_{-s}\mathscr{L} \cong \mathscr{L}$. This implies that the identifications given by $\mathscr{R}$ are compatible with the transition functions for $X_0$. Thus, the binary relation $\mathscr{R}$ globalizes to a binary relation on $X_0$.
\end{preuve}

We consider the equivalence relation on $X_0$ generated by the binary relation $\mathscr{R}$, we still denote by $\mathscr{R}$ the resulting relation on $X_0$.

We now have candidates to write globally the gluing relations for $S$ starting from $X_0$.

\begin{proposition} \label{quotient_map_new_gluing}
	The quotient \footnote{Since the equivalence relation is \'{e}tale, it surely does exist at least as an algebraic space.} $X_0/\mathscr{R}$ is isomorphic to $X$. Moreover, the quotient map $q : X_0 \twoheadrightarrow X_0/\mathscr{R}$ fits in the following commutative square:
	\begin{center}
		\begin{tikzcd}
	{\mathrm{Tot}(\mathscr{L}^\vee) \coprod \mathrm{Tot}(\mathscr{L})} & S \\
	{X_0} & X
	\arrow["{\tilde{j_0} \coprod j_\infty}"{pos=0.7}, from=1-1, to=1-2]
	\arrow["\nu", from=1-2, to=2-2]
	\arrow["{\nu_0}"', from=1-1, to=2-1]
	\arrow["q"', from=2-1, to=2-2]
\end{tikzcd}
	\end{center}
\end{proposition}

\begin{preuve}
    If $U \subseteq E$ is a trivializing open for $\mathscr{L}$, we define the map
    \begin{center}
        $f$ : \begin{tabular}{ccc}
           $U \times Z$  & $\to$     & $(U \times \mathbb{P}^1)/\sim$ \\
           $(e,(x,0))$   & $\mapsto$ & $(e+s, [x:1])$ \\
           $(e,(0,y))$   & $\mapsto$ & $(e, [1:y])$
        \end{tabular}
    \end{center}
    factors through (the equivalence relation generated by) $\mathscr{R}$ since
    \begin{center}
        $f(e+s, (0, x^{-1}))=(e, [x:1])=f(e,(x,0))$ and $f(e-s, (y^{-1},0))=(e, [1:y])=f(e,(0,y))$.
    \end{center}
    We still denote by $f$ the resulting map. Let us now construct a map in the other direction: we consider
    \begin{center}
        $g$ : \begin{tabular}{ccc}
           $U \times \mathbb{P}^1$  & $\to$     & $(U \times Z)/\mathscr{R}$ \\
           $(e,[x:y])$   & $\mapsto$ & $\begin{cases} (e-s,(xy^{-1},0)) & \textbf{if}~ y \neq 0 \\
                   (e, (0,yx^{-1}))& \textbf{if}~x \neq 0\end{cases}$
        \end{tabular},
    \end{center}
    which factors through (the equivalence relation generated by) $\sim$ since
    \begin{center}
        $g(e+s, [0:1])=(e,(0,0))=g(e,[1:0])$.
    \end{center}
    We still denote by $g$ the induced map; one can check that $f$ and $g$ are mutual inverses. 

    Let us now prove that $f$ and $g$ are compatible with the transition functions of the fiber bundle $X_0$. 

    Let $U$ and $V$ be two trivializing open subsets of $E$ for $\mathscr{L}$ such that $U \cap V \neq \emptyset$. Let us first check the compatibility for $f$: on one hand, $(e, (x,0))$ is mapped to $(e, (\lambda_{U,V}(e)^{-1}x,0))$ via the transition function $\varphi_{U,V}$ and this is mapped to $(e+s,[\lambda_{U,V}(e)^{-1}x:1])$ by $f$; on the other hand, $(e,(x,0))$ is mapped to $(e+s,[x:1])$ which is mapped to $(e+s,[x:\lambda_{U,V}(e+s)])$ by the transition function of the fiber bundle $S \to E$; these two points coincide since $t_s^*\mathscr{L} \cong \mathscr{L}$. One can check that it is also the case for the points of the form $(e,(0,y))$ so $f$ induces a morphism $X_0/\mathscr{R} \to X$.

    Let us now check this for $g$: let $(e,[x:y]) \in U \times \mathbb{P}^1$ and suppose for example that $x \neq 0$ (the case $y \neq 0$ is treated similarly). On one hand, $(e,[x:y])$ is mapped to $(e,(0, yx^{-1}))$ via $g$ and to $(e, (0, \lambda_{U,V}(x)yx^{-1}))$ via the transition function $\varphi_{U,V}$; on the other hand, $(e,[x:y])$ is mapped to $(e,[x:\lambda_{U,V}(e)y])$ via the transition function of $S \to E$ and this point is mapped to $(e, (0,\lambda_{U,V}(e)yx^{-1}))$ via $g$. the both point coincide so $g$ induces a morphism $X \to X_0/\mathscr{R}$ and since $f$ and $g$ were mutual inverses in the local construction, it is still true globally so $X_0/\mathscr{R}$ is isomorphic to $X$.

    With the explicit description of the isomorphism above, one can check that the square of the proposition is commutative.
\end{preuve}

\begin{lemme}
    The morphism $q : X_0 \to X$ is \'{e}tale.
\end{lemme}

\begin{preuve}
    Let $x_0 \in X_0$. We have to check that the morphism
    \begin{equation} \label{morphisme_q_completion}
        \hat{\mathscr{O}}_{X, q(x_0)} \to \hat{\mathscr{O}}_{X_0,x_0}
    \end{equation}
    is an isomorphism; this can be done locally around $x_0$ so we can suppose that $X_0=U \times Z$ and one can use the description of the isomorphism $X_0/\mathscr{R} \cong X$ given above.
\end{preuve}

\subsubsection{Existence of a global non-torsion 1-form on $X$.}

We now construct a global non-torsion 1-forms on $X$. Let $\beta$ be a global generator of the $\mathscr{O}_E$-module $\Omega^1_E$ and let us consider $\alpha:=\pi^*\beta \in \Gamma(X_0, \Omega^1_{X_0})$.

\begin{proposition} \label{existence_non_torsion_form}
	There exists $\omega \in \Gamma(X, \Omega^1_{\underline{X}})$ such that $q^*\omega=\alpha$. Moreover, we have $\nu^*\omega=p^*\beta$; in particular $\omega \notin \Gamma(X, \tau^1_X)$.
\end{proposition}

\begin{preuve}
If $U \subseteq E$ is a trivializing open for $\mathscr{L}$, we have
\begin{center}
    $\alpha|_{\pi^{-1}(U)}=(\beta|_U, 0) \in \Omega^1_{\pi^{-1}(U)} \cong \Omega^1_U \boxplus \Omega^1_{V(xy)}$.
\end{center}
Since $\beta$ is translation invariant, $\beta|_U$ is also translation invariant and therefore, by \'{e}tale descent, $\alpha|_U$ descends to a form $\omega_U \in \Gamma(q(U), \Omega^1_{\underline{X}})$. Now, if we take $V \subseteq E$ to be another trivializing open for $\mathscr{L}$, the transition function $\varphi_{U,V} : (U \cap V) \times Z \to (U \cap V) \times Z$  for the fiber bundle $X_0 \to E$ are given by
\begin{center}
    $(p, (x,0)) \mapsto (p, (\lambda_{U,V}(p)^{-1}x, 0))$ and $(p, (0,y)) \mapsto (p, (0, \lambda_{U,V}(p)y))$.
\end{center}
As in the proof of Lemma \ref{globalization_equiv_relation}, the map $\varphi_{U,V}$ descends to a map
\begin{center}
    $\tilde{\varphi}_{U,V} : q(U) \cap q(V) \to q(U) \cap q(V)$
\end{center}
and if we denote by $\omega_V$ the 1-form obtained from $\alpha|_{\pi^{-1}(V)}$ by \'{e}tale descent, we get
\begin{center}
    $\omega_U|_{q(U) \cap q(V)}=\tilde{\varphi}_{U,V}^*\omega_V|_{q(U) \cap q(V)}$.
\end{center}
Moreover, since the transition functions for $\mathscr{L}$ form a 1-cocycle, the $\varphi_{U,V}$'s do and thus, the $\tilde{\varphi}_{U,V}$ do as well. This implies that the 1-forms $\omega_U$ glue together to $\omega \in \Gamma(X, \Omega^1_{\underline{X}})$ such that $q^*\omega=\alpha$.
 
Now, by commutativity of the square of Proposition \ref{quotient_map_new_gluing}, we have that
	\begin{center}
		$\nu_0^*\pi^*\beta=\nu_0^*q^*\omega=(\tilde{j_0} \coprod j_\infty)^*\nu^*\omega$.
	\end{center}
	Thus, the relation $\nu^*\omega=p^*\beta$ is true on $\mathrm{Tot}(\mathscr{L}) \cong S \setminus D_\infty$ and on $\mathrm{Tot}(\mathscr{L}^\vee) \cong S \setminus D_0$; it is then true on $S$ since these two open subsets form a covering of $S$. The fact that $\omega$ is not a torsion 1-form follows from the fact that $\nu^*\omega=p^*\beta \neq 0$ and the observation that the pullback of a torsion 1-form on $X$ along $\nu$ is zero. 
\end{preuve}

% \begin{remarque}
% 	Since the gluing relations describing $X$ as a gluing of $X_0$ and as a gluing of $S$ are similar, one could wonder why this argument does not work directly to show that the form $p^*\beta$ on $S$ descends to a 1-form on $X$; in fact, it does not work since $\nu$ is not even flat since it's the normalization of a non-normal scheme, whereas $q$ which is even \'{e}tale so descent theory applies. 
% \end{remarque}

\begin{corollaire} \label{comp_coh_tgt_gluing}
	We have $h^1(X, \Omega^1_X)=4$; in particular, this number does neither depend on $\mathscr{L}$ nor on $s$. Therefore, the $\mathbb{C}$-vector space $\mathrm{Ext}^1_{\mathscr{O}_X}(\Omega^1_X, \mathscr{O}_X)$ has dimension 4 and in particular, $X$ has a non-locally trivial first order deformation. 
\end{corollaire} 

\begin{preuve}
	The short exact sequence
	\begin{center}
		$0 \longrightarrow \tau_X^1 \longrightarrow \Omega^1_{\underline{X}} \longrightarrow \tilde{\Omega}^1_X \longrightarrow 0$
	\end{center}
	yields a long exact sequence in sheaf cohomology
	\begin{center}
		$0 \to \mathrm{H}^0(X, \tau^1_X) \to \mathrm{H}^0(X, \Omega^1_{\underline{X}}) \to \mathrm{H}^0(X, \tilde{\Omega}^1_X) \to \mathrm{H}^1(X, \tau^1_X) \to \mathrm{H}^1(X, \Omega^1_{\underline{X}}) \to \mathrm{H}^1(X, \tilde{\Omega}^1_X) \to 0$.
	\end{center}
	Replacing all the known data, one proves that two cases would a priori be possible:
	\begin{enumerate}
		\item $h^0(X, \Omega^1_{\underline{X}})=1$ and $h^1(X, \Omega^1_{\underline{X}})=3$
		\item $h^0(X, \Omega^1_{\underline{X}})=2$ and $h^1(X, \Omega^1_{\underline{X}})=4$.
	\end{enumerate}
	The first case cannot occur: indeed, since $\tau^1_X \cong \mathscr{O}_D$, we have $h^0(X, \tau^1_X)=1$ and therefore $h^0(X, \Omega^1_{\underline{X}}) \geq 1$; by Proposition \ref{existence_non_torsion_form}, there exists a non-torsion 1-form on $X$ so $h^0(X, \Omega^1_{\underline{X}}) \geq 2$, we thus deduce that
    \begin{center}
        $h^0(X, \Omega^1_{\underline{X}})=2$.
    \end{center}
    Now, by Serre's duality for coherent sheaves, we have
    \begin{center}
        $\dim_\mathbb{C}(\mathrm{Ext}^1_{\mathscr{O}_X}(\Omega^1_{\underline{X}}, \mathscr{O}_X))=h^1(X, \Omega^1_{\underline{X}})=4$.
    \end{center}
\end{preuve}

\subsection{Link between log-smooth and classical flat deformations of $X$.}

In this subsection, we discuss the link between the logarithmic and the classical deformations of $X$.

\begin{paragraphe} \label{construction_morphism_torsionfree_forms_to_log_forms}
    By construction of the sheaf of logarithmic forms of a logarithmic scheme (see Definition \ref{definition_log_cotangent_sheaf} in the appendix), we have a morphism of $\mathscr{O}_X$-modules
    \begin{center}
        $\Omega^1_{\underline{X}} \to \Omega^1_{X/S_0}$.
    \end{center}
    \'{E}tale locally, it is given by
    \begin{center}
        $[\overline{f_0} \mathrm{d}x + \overline{f_1} \mathrm{d}y + \overline{f_2} \mathrm{d}z] \mapsto (x\overline{f_0}-y\overline{f_1}) \frac{\mathrm{d}x}{x} + \overline{f_2}\mathrm{d}z$.
    \end{center} 
    In particular, the kernel of $\Omega^1_{\underline{X}} \to \Omega^1_{X/S_0}$ is equal to $\tau^1_X$ so we have an injective morphism of $\mathscr{O}_X$-modules
    \begin{center}
        $\tilde{\Omega}^1_X \hookrightarrow \Omega^1_{X/S_0}$.
    \end{center}
   % One can also see that noting that the kernel of this morphism is supported on $D$ and is then a torsion subsheaf of $\tilde{\Omega}^1_X$; the latter being by construction a torsion-free sheaf, the kernel has to be zero.
\end{paragraphe}

\begin{proposition} \label{comparison_torsionfree_differentials_with_log_differentials}
    The morphism constructed in \ref{construction_morphism_torsionfree_forms_to_log_forms} takes place in the short exact sequence of $\mathscr{O}_X$-modules
    \begin{center}
        $0 \to \tilde{\Omega}^1_X \to \Omega^1_{X/S_0} \to j_*\mathscr{O}_D \to 0$
    \end{center}
    where the morphism on the right is given by the composition
    \begin{center}
        $\Omega^1_{X/S_0} \xrightarrow{\nu^*} \nu_*\Omega^1_{\overline{S}/S_0} \xrightarrow{\cdot|_{D_0}} (\nu \circ j_0)_* \Omega^1_{D_0/S_0} \to (\nu \circ j_0)_*\Omega^1_{D_0/F_0} \cong j_*\mathscr{O}_D$.
    \end{center}
\end{proposition}

\begin{preuve}
    Recall that the morphism $(\nu \circ j_0)_* \Omega^1_{D_0/S_0} \to (\nu \circ j_0)_*\Omega^1_{D_0/F_0}$ comes from pushing forward the morphism on the right of the exact sequence of $\mathscr{O}_{D_0}$-modules
    \begin{center}
        $\Omega^1_{F_0/S_0} \to \Omega^1_{D_0/S_0} \to \Omega^1_{D_0/F_0} \to 0$
    \end{center}
    and that the isomorphism in the end comes from Corollary \ref{triviality_log_forms_D0_over_F0}. One can check that \'{e}tale locally, the composition is given by the morphism of $R$-modules
    \begin{center}
        $f \frac{\mathrm{d}x}{x} + g \mathrm{d}z \in R \cdot \frac{\mathrm{d}x}{x} \oplus R \cdot \mathrm{d}z \mapsto f(0,0,z) \in \mathbb{C}[z]$
    \end{center}
    where $R:=\mathbb{C}[x,y,z]/(xy)$. Let us now prove exactness at each spot separately.
    \begin{itemize}
        \item \textbf{Exactness on the left}: it was already proven in the paragraph above.
        \item \textbf{Exactness on the right}: any element $f(z) \in \mathbb{C}[z]$ is the image of the element $f \frac{\mathrm{d}x}{x}$.
        \item \textbf{Exactness in the middle}: if we have $\overline{f} \frac{\mathrm{d}x}{x} + \overline{g} \mathrm{d}z \in R \cdot \frac{\mathrm{d}x}{x} \oplus R \cdot \mathrm{d}z$ such that $f(0,0,z)=0$, it implies that
        \begin{center}
            $x \mid f(x,0,z)$ and $y \mid f(0,y,z)$
        \end{center}
        where the divisibility relations respectively take place in $\mathbb{C}[x,z]$ and $\mathbb{C}[y,z]$. In particular, we can write
        \begin{center}
            $f(x,0,z)=x \cdot \tilde{f}_0(x,z)$ and $f(0,y,z)=y \cdot \tilde{f}_1(y,z)$
        \end{center}
        for some $\tilde{f}_0 \in \mathbb{C}[x,z]$ and $\tilde{f}_1 \in \mathbb{C}[y,z]$. Now, the image of the element $[\overline{\tilde{f}_0} \mathrm{d}x-\overline{\tilde{f}_1} \mathrm{d}y+\overline{g} \mathrm{d}z] \in \tilde{\Omega}^1_{\underline{X}}$ in $\Omega^1_{X/S_0}$ is equal to
        \begin{center}
            $(x \overline{\tilde{f}_0}+y \overline{\tilde{f}_1}) \frac{\mathrm{d}x}{x} + \overline{g} \mathrm{d}z$
        \end{center}
        but because of the congruence
        \begin{center}
            $x \tilde{f}_0 + y \tilde{f_1}= f(x,0,z)+f(0,y,z) \equiv f(x,y,z) \pmod{xy}$,
        \end{center}
        its image is equal to the element we started with. Conversly, if we start with a germ $\alpha$ of section of $\tilde{\Omega}^1_{\underline{X}}$, its image in $j_*\mathscr{O}_D$ will be equal, using the same notations as above for the coefficients of $\alpha$, to $xf_0(x,0,z)|_{x=0}$ which is indeed zero.
    \end{itemize}
    This proves the exactness of the sequence of $\mathscr{O}_X$-modules of the proposition.
    \end{preuve}

We now construct a short exact sequence of $\mathscr{O}_X$-modules involving $\mathscr{T}_{X/S_0}$ and $\mathscr{T}_{\underline{X}}$ in order to derive a map between the tangent spaces of the functor of log smooth deformations of $X$ and locally trivial (classical) deformations of $\underline{X}$.

\begin{proposition} \label{comparison_classical_tangent_sheaf_with_log_tangent_sheaf}
   Dualizing the short exact sequence of $\mathscr{O}_X$-modules in Proposition \ref{comparison_torsionfree_differentials_with_log_differentials} yields a short exact sequence of $\mathscr{O}_X$-modules
   \begin{center}
       $0 \to \mathscr{T}_{X/S_0} \to \mathscr{T}_{\underline{X}} \to \mathscr{E}xt^1_{\mathscr{O}_X}(j_*\mathscr{O}_D, \mathscr{O}_X) \to 0$.
   \end{center}
\end{proposition}

\begin{theoreme} \label{log_defo_first_order_gluing}
    The non-locally trivial first-order deformation of $X$ can be equipped with a logarithmic structure of semi-stable type.
\end{theoreme}

\begin{preuve}
Let us denote by $\mathscr{X}$ this family of deformation, it sits in the following cartesian square:
    \begin{center}
        \begin{tikzcd}
	X & {\mathscr{X}} \\
	{\mathrm{Spec}(\mathbb{C})} & {\mathrm{Spec}(\mathbb{C}[\varepsilon])}
	\arrow[from=1-1, to=2-1]
	\arrow["\iota", from=1-1, to=1-2]
	\arrow[from=2-1, to=2-2]
	\arrow["\psi", from=1-2, to=2-2]
	\arrow["\lrcorner"{anchor=center, pos=0.125}, draw=none, from=1-1, to=2-2]
\end{tikzcd}
    \end{center}
    Let us construct a logarithmic structure on $\mathscr{X}$ by hand relying on the following observation: if $X$ is smoothable, we should obtain eventually an algebraic deformation $\mathscr{X} \to \mathrm{Spec}(\mathbb{C}\llbracket t \rrbracket)$ of $X$; in that setting, we can equip $\mathscr{X}$ with the divisorial logarithmic structure associated to the divisor $X$ in $\mathscr{X}$. Therefore, if $U \cong \mathrm{Spec}(\mathbb{C}[\varepsilon][z_0,z_1,z_2]/(z_0z_1-\alpha \varepsilon)) \to \mathscr{X}$ is an \'{e}tale chart of $\mathscr{X}$, we can define a sheaf of monoids $\mathscr{M}_U$ by
    \begin{center}
        $\Gamma(U, \mathscr{M}_U):= \lbrace z_0^pz_1^qu \in \Gamma(U, \mathscr{O}_U) \mid (p,q) \in \mathbb{N}^2, u \in \Gamma(U, \mathscr{O}_U^\times) \rbrace$.
    \end{center}
    We have a natural morphism of sheaves $\mathscr{M}_U \to \mathscr{O}_U$ given by the inclusion. 

    Let us prove that we can cover $\mathscr{X}$ with charts whose coordinates restrict to coordinates of a logarithmic altas (in the sense of \cite[1.]{KN94}) on $X$. In order to do that, we pick a logarithmic atlas $\lbrace (U_i, x_0^{(i)}, x_1^{(i)}, x_2^{(i)}) \rbrace_{i \in I}$ of $X$ and let us fix $ i\in I$. 
    Since the sequence
    \begin{center}
        $0 \to \varepsilon \mathscr{O}_{\mathscr{X}} \to \mathscr{O}_{\mathscr{X}} \to \iota_*\mathscr{O}_X \to 0$
    \end{center}
    is exact, it is also exact on the stalks so, if we shrink $U_i$ enough and restrict the $x_j^{(i)}$'s, we can lift them to functions $z_j^{(i)}$'s on a small open of $\mathscr{X}$. Since the $z_j^{(i)}$'s restrict to the $x_j^{(i)}$'s, we have that
    \begin{center}
        $z_0^{(i)}z_1^{(i)}=\varepsilon f^{(i)}$
    \end{center}
    where $f^{(i)}$ is a section of $\mathscr{O}_{\mathscr{X}}$. Now, since $\mathscr{X}$ is an infinitesimal smoothing of $X$, it implies that $f^{(i)}$ is not identically zero so, if we shrink again the open subset of $\mathscr{X}$ on which it is defined, we can suppose that $f^{(i)}$ is invertible. We can now replace one of the $z_j^{(i)}$'s by $z_j^{(i)} \alpha (f^{(i)})^{-1}$.

    This implies that for all $\ell \in \lbrace 0, 1 \rbrace$
    \begin{center}
        $z_\ell^{(i)}=\tilde{u}_{\ell}^{(i,j)} \cdot z_\ell^{(j)} + \varepsilon f_\ell^{(i,j)}$
    \end{center}
    where $f_\ell^{(i,j)} \in \Gamma(U_i \cap U_j, \mathscr{O}_{\mathscr{X}})$ and $\tilde{u_\ell}$ is of the form $u_\ell+\varepsilon g_\ell$ with $g_\ell \in \Gamma(U_i \cap U_j, \mathscr{O}_{\mathscr{X}})$; it is an invertible function since $u_\ell$ is and $\varepsilon g_\ell$ is nilpotent so its value at any point in zero. Since the $u_\ell$'s and the $\varepsilon g_\ell$'s are unique, the $\tilde{u_\ell}$'s are also unique. In fact, for all $\ell \in \lbrace 0, 1 \rbrace$, $f_\ell^{(i,j)}$ is divisible by $z_\ell^{(j)}$. 
    
    % Indeed, by construction, we have $z_0^{(i)}z_1^{(i)}=\alpha \varepsilon$, which means that 
    % \begin{center}
    %     $\tilde{u}_0^{(i,j)}\tilde{u}_1^{(i,j)} z_0^{(j)}z_1^{(j)} + \varepsilon (f_0^{(i,j)}z_0^{(j)}\tilde{u}_1^{(i,j)} + f_1^{(i,j)}z_0^{(j)} \tilde{u}_0^{(i,j)}) = \alpha \varepsilon$
    % \end{center}
    % which rewrites as
    % \begin{center}
    %     $\alpha \varepsilon (1- \tilde{u}_0^{(i,j)}\tilde{u}_1^{(i,j)}) = \varepsilon (f_0^{(i,j)}z_1^{(j)}\tilde{u}_1^{(i,j)} + f_1^{(i,j)}z_0^{(j)} \tilde{u}_0^{(i,j)})$
    % \end{center}
    % since $z_0^{(j)}z_1^{(j)}=\alpha \varepsilon$. Since $u_0^{(i,j)}u_1^{(i,j)}=1$ and $\varepsilon^2=0$, we have that $\alpha \varepsilon(1-\tilde{u}_0^{(i,j)}\tilde{u}_1^{(i,j)})=0$, which implies that
    % \begin{center}
    %     $f_0^{(i,j)}z_1^{(j)}\tilde{u}_1^{(i,j)}+f_1^{(i,j)}z_0^{(j)}\tilde{u}_0^{(i,j)}=0$
    % \end{center}
    % which in turn implies that
    % \begin{center}
    %     $f_0^{(i,j)}z_1^{(j)}+f_1^{(i,j)}z_0^{(j)}\tilde{u}_0^{(i,j)}(\tilde{u}_1^{(i,j)})^{-1}=0$
    % \end{center}
    % so $z_\ell^{(j)}$ divides $f_\ell^{(i,j)}$ for all $\ell \in \lbrace 0, 1 \rbrace$. 
    In particular, we can write
    \begin{center}
        $z_\ell^{(i)}= \tilde{u}_i^{(i,j)}(1+\frac{\varepsilon f_\ell^{(i,j)}}{\tilde{u}_\ell^{(i,j)} z_\ell^{(j)}})z_\ell^{(j)}$
    \end{center}
    for all $\ell \in \lbrace 0,1 \rbrace$ and $\tilde{u}_\ell^{(i,j)}(1+\frac{\varepsilon f_\ell^{(i,j)}}{\tilde{u}_\ell^{(i,j)} z_\ell^{(j)}})$ is invertible. Moreover, it is unique since $\varepsilon f_\ell^{(i,j)}$ and $\tilde{u}_{\ell}^{(i,j)}$ are unique. In particular, this unicity property implies that the gluing maps for the $\mathscr{M}_U$'s satisfy a 1-cocyle relation. As a consequence, the \'{e}tale-locally defined sheaves $\mathscr{M}_U$ glue to a sheaf of monoids $\mathscr{M}_{\mathscr{X}}$ and this gluing is compatible with the inclusion to $\mathscr{O}_U$ so we obtain a prelogarithmic structure $\alpha : \mathscr{M}_{\mathscr{X}} \to \mathscr{O}_{\mathscr{X}}$. In fact, it is a logarithmic structure since the image of a local section $z_0^pz_1^qu$ in $\mathscr{O}_{\mathscr{X}}$ is invertible if and only if $p=q=0$, hence $\alpha$ induces an isomorphism
    \begin{center}
        $\alpha^{-1}(\mathscr{O}_{\mathscr{X}}^\times) \to \mathscr{O}_{\mathscr{X}}^\times$
    \end{center}
    and therefore $\alpha$ is a logarithmic structure on $\mathscr{X}$. By construction, we have \'{e}tale locally a chart for $\mathscr{M}_{\mathscr{X}}$ given by
    \begin{center}
        $(p,q) \in \mathbb{N}^2 \mapsto \lambda^p z_0^pz_1^q \in \Gamma(U, \mathscr{M}_U)$.
    \end{center}

    Now, we want to see whether the morphisms
    \begin{center}
        $(k,v) \in \mathbb{N} \oplus \mathbb{C}[\varepsilon] \mapsto z_0^kz_1^kv \in \Gamma(U, \mathscr{M}_U)$
    \end{center}
    defined for an \'{e}tale chart $U=\mathrm{Spec}(\mathbb{C}]z_1,z_1,z_2]/(z_0z_1)) \to X$ globalize to a morphism of sheaves of monoids \footnote{Here, $S_{1, \alpha}$ is the logarithmic scheme whose underlying scheme is $\mathrm{Spec}(\mathbb{C}[\varepsilon])$ and whose log structure is defined by $1 \in \mathbb{N} \mapsto \alpha \varepsilon \in \mathbb{C}[\varepsilon]$.} $\psi^\flat : \mathscr{M}_{S_{1, \alpha}} \to \psi_*\mathscr{M}_{\mathscr{X}}$.
    
    Let us take $U \to \mathscr{X}$ and $U' \to \mathscr{X}$ two \'{e}tale charts. We have to check that the square
    \begin{center}
        \begin{tikzcd}
	{\psi_*\mathscr{M}_{U}|_{U \times_\mathscr{X} U'}} & {\mathscr{M}_{S_{1, \alpha}}} \\
	{\psi_*\mathscr{M}_{U'}|_{U \times_\mathscr{X} U'}} & {\mathscr{M}_{S_{1, \alpha}}}
	\arrow[Rightarrow, no head, from=1-2, to=2-2]
	\arrow["{\psi^\flat_U}", from=1-2, to=1-1]
	\arrow["{\psi^\flat_{U'}}"', from=2-2, to=2-1]
	\arrow["\varphi_{U,U'}", from=1-1, to=2-1]
\end{tikzcd}
    \end{center}
    is commutative and that if $U'' \to \mathscr{X}$ is a third \'{e}tale chart such that $U \cap U' \cap U'' \neq \emptyset$, then we have the 1-cocyle relation
    \begin{center}
        $\varphi_{U,U'} \circ \varphi_{U',U''}=\varphi_{U,U''}$    
    \end{center} 
    on the triple intersection. Again, if we change the chart, the coordinates $z_0$ and $z_1$ will be multiplied by a section of $\mathscr{O}_{\mathscr{X}}^\times$ so the map $\varphi_{U,U'}$ 
    will just be given by
    \begin{center}
        $\varphi_{U,U'} : z_0^pz_1^qu \in \Gamma(U \times_\mathscr{X} U', \mathscr{M}_{U}) \mapsto \tilde{u_0}^pz_0^p\tilde{u_1}^qz_1^qu \in \Gamma(U \times_\mathscr{X} U', \mathscr{M}_{U'})$.    
    \end{center}
    and, again, by unicity of the $\tilde{u_i}$'s, it satisfies the 1-cocycle relation above and the square is commutative. Moreover, because of the chart described above, the morphism of log schemes $(\mathscr{X}, \mathscr{M}_{\mathscr{X}}) \to S_{1, \alpha}$ is log smooth by Kato's toroidal criterion (Theorem \ref{kato_toroidal_criterion}) for log smoothness --- \'{e}tale locally, $\psi$ is just a base change of $X \to S_0$, which is log smooth. Therefore, it defines a logarithmic deformation of $X$ at first order.
\end{preuve}

\section{Smoothability result and Kuranishi germs.}

By theorems \cite[Th\'{e}or\`{e}me principal, p. 598]{Dou74} of Douady and \cite[Hauptsatz, p. 140]{Gra74} of Grauert, we know that for any compact complex analytic space $Z$, there exists a versal deformation family $f : \mathscr{Z} \to \mathrm{Def}(Z)$ of $Z$ whose base is called the \emph{Kuranishi space} of $Z$. This deformation of $Z$ is called the \emph{Kuranishi family} of $Z$. It is usual to write $0 \in \mathrm{Def}(Z)$ for the point such that $f^{-1}(0) \cong Z$.

By \cite[(0.3) Corollary]{FK87}, there exists a versal locally trivial deformation family $\mathscr{Z}^{\mathrm{lt}} \to \mathrm{Def}^{\mathrm{lt}}(Z)$ of $Z$ which is called the \emph{locally trivial Kuranishi family} of $Z$ and whose base is called the \emph{locally trivial Kuranishi space} of $Z$, it is a closed complex subspace of $\mathrm{Def}(Z)$.

\subsection{Smoothability result.}

We now come back to our initial smoothability question.

\begin{theoreme} \label{functor_log_def_smooth_gluing}
    The functor $\mathrm{LD}_X$ is smooth.
\end{theoreme} 

\begin{preuve}
    We want to apply the $T^1$-lifting principle (\cite[Theorem A]{FM99}), i.e., we have to prove that, for all $n \geq 0$, for any logarithmic deformation $\mathscr{X}_n$ of $X$ over $S_n:=\mathrm{Spec}(\mathbb{N} \xrightarrow{1 \mapsto \overline{t}} \mathbb{C}[t]/(t^{n+1}))$ and for any lifting $\mathscr{X}_{n+1}$ to order $n+1$, the map
	\begin{center}
		$\mathrm{H}^1(\mathscr{X}_{n+1}, \mathscr{T}_{\mathscr{X}_{n+1}/S_n}) \to \mathrm{H}^1(\mathscr{X}_n, \mathscr{T}_{\mathscr{X}_{n}/S_{n-1}})$
	\end{center}
	is surjective. This follows from \cite[Theorem 1.10.]{FRF21} that says that $\mathrm{H}^q(X, \Omega^p_{X/S_n})$ is a free $\mathbb{C}[t]/(t^{n+1})$-module whose formation commutes with base change. This ensures that the map is surjective for all $n > 0$ but for $n=0$, the right hand side is just the singleton $\lbrace X \rbrace$ so the surjectivity of the map means in that case that for any $\alpha \in \mathbb{C}$, $X$ has a first-order logarithmic deformation over $S_{1, \alpha}:=\mathrm{Spec}(\mathbb{N} \xrightarrow{1 \mapsto \alpha \varepsilon} \mathbb{C}[\varepsilon])$. But this is exactly guarateed by Proposition \ref{log_defo_first_order_gluing}.
\end{preuve}

\begin{remarque}
    This proves that the set of first order logarithmic deformations of $X$ is non-empty and therefore a torsor under the action of $\mathrm{H}^1(X, \mathscr{T}_X)$.
\end{remarque}

\begin{corollaire} \label{smoothability_result}
	The variety $X$ is smoothable as a complex space.
\end{corollaire}

\begin{preuve}
	By Theorem \ref{functor_log_def_smooth_gluing}, the functor $\mathrm{LD}_X$ is smooth, so the map
	\begin{center}
		$\mathrm{LD}_X(\mathbb{C}[\varepsilon]) \twoheadrightarrow \mathrm{LD}_X(\mathbb{C})=\lbrace [X] \rbrace$
	\end{center}
	where the logarithmic structure on $\mathrm{Spec}(\mathbb{C}[\varepsilon])$ is given by $1 \mapsto \varepsilon$, is surjective. In particular, we have a first order log deformation of $X$ which is not locally trivial. By smoothness again, we can lift it to any order and at the limit, we get a formal logarithmic deformation $\hat{f} : \mathscr{Z} \to \mathrm{Spf}(\mathbb{C} \llbracket t \rrbracket)$. On the other hand, we have the Kuranishi family $\mathscr{X} \to W:=\mathrm{Def}(X)$ which induces a germ of formal family 
	\begin{center}
		$\widehat{(\mathscr{X}, [X])} \longrightarrow \widehat{(W, w)}$.
	\end{center}
	By versality of the Kuranishi family, we get a morphism \footnote{This morphism corresponds to a (unique) morphism of germs of formal complex spaces $\widehat{(\Delta, 0)} \to \widehat{(W, w)}$.} of $\mathbb{C}\llbracket t \rrbracket$-algebras
	\begin{center}
		$\hat{\mathscr{O}}_{W,w} \longrightarrow \mathbb{C} \llbracket t \rrbracket=\hat{\mathscr{O}}_{\Delta, 0}$
	\end{center}
	%inducing the following commutative triangle:
	%\begin{equation} \label{commutative_triangle-kuranishi}
		%\begin{tikzcd}
	%{\widehat{(\mathscr{X}, [X])} \times_{\widehat{(W,w)}} \widehat{(\Delta, 0)}} & (\mathscr{Z}, [X]) \\
	%{\widehat{(\Delta, 0)}}
	%\arrow["{\mathrm{pr}_2}"', from=1-1, to=2-1]
	%\arrow["{\hat{f}}", from=1-2, to=2-1]
	%\arrow["\cong", from=1-1, to=1-2]
%\end{tikzcd}
	%\end{equation}
% 	Now, by Artin's approximation theorem (\cite[Theorem A. 10. 9.]{Alp23}), there exists a residually trivial \'{e}tale morphism $(W',w') \to (W,w)$ and a morphism of $\mathbb{C}\llbracket t \rrbracket$-algebras $\mathscr{O}_{W',w'} \to \hat{\mathscr{O}}_{\Delta, 0}$ coinciding with $\overline{\alpha}$ up to first order --- by this, we mean that the composition
% 	\begin{center}
% 		$\overline{\alpha}' : \mathscr{O}_{W,w} \to \mathscr{O}_{W',w'} \to \hat{\mathscr{O}}_{\Delta, 0}$
% \end{center}
% also coincides with $\overline{\alpha}$ up to first order. This is not so much of an abuse of language because of the residual triviality of the \'{e}tale morphism $(W',w') \to (W,w)$ which implies that 
% 	\begin{center}
% 		$\mathscr{O}_{W',w'}/\mathfrak{m}_{w'}^2 \cong \mathscr{O}_{W,w}/\mathfrak{m}_w^2$.
% \end{center}
and precomposing with the canonical morphism $\mathscr{O}_{W,w} \to \hat{\mathscr{O}}_{W,w}$, we get a morphism
\begin{center}
    $\overline{\alpha} : \mathscr{O}_{W,w} \to \hat{\mathscr{O}}_{\Delta,0}$.
\end{center}
Now, by Artin approximation theorem \cite[Theorem 1.5a, (i)]{Art68}, we obtain a morphism of $\mathbb{C}$-algebras
\begin{center}
	$\alpha : \mathscr{O}_{W,w} \to \mathscr{O}_{\Delta, 0}$
\end{center}
which coincides with $\overline{\alpha}$ up to first order. This morphism $\alpha$ induces a morphism of germs of complex spaces 
\begin{center}
	$(\Delta, 0) \to (W,w)$
\end{center}			 
and we can consider the (proper) morphism of germs of analytic spaces
\begin{center}
	$f : \Big( \mathscr{X} \times_W \Delta, ([X], 0) \Big) \longrightarrow (\Delta, 0)$.
\end{center}
Since locally around the singular locus of the central fiber, $\mathscr{Z}$ is defined by the equation $xy=at$ for $a \in \mathbb{C}^\times$, that $t$ is not killed at first order so the germ $(\mathscr{X} \times_W \Delta, ([X], 0))$ is smooth and therefore that $f$ is a germ of smoothing of $X$ --- which can be defined on a sufficiently small disk around the origin. Hence, $X$ is smoothable as an analytic space.
\end{preuve}

%\begin{remarque}
  %  We do not know if the logarithmic structure also algebraizes.
%\end{remarque}

\subsection{Possible semistable smoothings of $X$.}

In that section, we prove that the nearby fiber of any semistable smoothing of $X$ must be an abelian surface. This is consistent with the fact that $X$ appears as a degeneration of an abelian surface.

\begin{paragraphe}
	If $\pi : \mathscr{X} \to \Delta$ is a semistable degeneration of $X$, then by construction --- \cite[3.3.]{Del70} for the logarithmic de Rham sheaf and \cite[2.]{KN94} for the sheaf of logarithmic differentials --- we have that, for all $p \geq 0$:
	\begin{center}
		$\Omega^p_{X/S_0}=\Omega^p_{\mathscr{X}/\Delta}(\log X) \otimes_{\mathscr{O}_{\mathscr{X}}} \mathscr{O}_X$,
	\end{center}
 where $X$ is equipped with a logarithmic structure of semistable type.
\end{paragraphe}

From \cite[Proposition 2.16]{St76} and Ehresmann lemma, we deduce the following result: 

\begin{proposition} \label{cohomology_nearby_fiber}
	For any semistable degeneration $\pi : \mathscr{X} \to \Delta$ of $X$ and for any $t \in \Delta$, we have an isomorphism $\mathrm{H}^k(X_t, \mathbb{C}) \cong \mathbb{H}^k(X, \Omega^\bullet_{X/S_0})$ for all $k \geq 0$, where $X_t:=\pi^{-1}(t)$.
\end{proposition}

% \begin{remarque}
% 	This isomorphism is not canonical since it depends on the choice of $t \in \Delta - \lbrace 0 \rbrace$.
% \end{remarque}

From the triviality of the dualizing sheaf of $X$, we can compute some logarithmic Hodge numbers of $X \to S_0$.

\begin{lemme} \label{calcul h^1,0 log}
    The linear map
    \begin{center}
        $\varphi : \mathrm{H}^0(\overline{S}, \Omega^1_{\overline{S}/S_0}) \to \mathrm{H}^0(D, \Omega^1_{D/S_0})$
    \end{center}
    induced by the short exact sequence of Proposition \ref{ses_log1-forms} is zero. In particular, we have
    \begin{center}
        $h^0(X, \Omega^1_{X/S_0})=2$.
    \end{center}
\end{lemme}

\begin{preuve}
    By Lemma \ref{identification_ghost_sheaves}, we have $\overline{\mathscr{M}}_D \cong \underline{\mathbb{N}}_D^{\oplus 2}$, which implies that
    \begin{center}
        $\Omega^1_{D/S_0} \cong \Omega^1_{\underline{D}} \oplus \mathscr{O}_D$.
    \end{center}
    In particular, $\varphi$ splits in two linear maps
    \begin{center}
        $\mathrm{H}^0(\overline{S}, \Omega^1_{\overline{S}/S_0}) \to \mathrm{H}^0(D, \Omega^1_{\underline{D}})$ and $\mathrm{H}^0(\overline{S}, \Omega^1_{\overline{S}/S_0}) \to \mathrm{H}^0(D, \mathscr{O}_D)$.
    \end{center}
    The second map is given by the sum of the residues along $D_0$ and $D_\infty$. It is zero because each of the residue will be a global section of $\mathscr{O}_{D_i}$, i.e., a constant; now, the sum of these constant is zero because of the residue theorem applied to a fiber of $p : S \to E$ (which is isomorphic to $\mathbb{P}^1$). Let us now identify the first map and show that it is zero. If $\alpha \in \mathrm{H}^0(\overline{S}, \Omega^1_{\overline{S}/S_0})$, then for all $i \in \lbrace 0, \infty \rbrace$, we have
    \begin{center}
        $\alpha|_{D_i} \in \mathrm{H}^0(D_i, \Omega^1_{D_i/S_0}) \cong \mathrm{H}^0(D_i, \Omega^1_{\underline{D_i}}) \oplus \mathrm{H}^0(D_i, \mathscr{O}_{D_i})$.
    \end{center}
    By construction of the pullback of logarithmic forms, the first component of $\alpha|_{D_i}$ in the decomposition above comes from $\mathrm{H}^0(\overline{S}, \Omega^1_{\underline{\overline{S}}})$ which is one-dimensional generated by the pullback of a global generator of the $\mathscr{O}_E$-module $\Omega^1_E$. Since the first components of $\alpha|_{D_0}$ and $\alpha|_{D_\infty}$ both come from the same element of $\mathrm{H}^0(\overline{S}, \Omega^1_{\underline{\overline{S}}})$, we have
    \begin{center}
        $\alpha|_{D_0}-\alpha|_{D_\infty}=0$,        
    \end{center}
    which implies that the map is zero. 
\end{preuve}

We continue with a degeneration result of the logarithmic Hodge--de Rham spectral sequence for a semistable smoothing.

\begin{proposition} \label{degeneration_hodge_de_rham_lissage}
    Let $f : \mathscr{X} \to \Delta$ be a semistable smoothing of $X$. If we equip $\mathscr{X}$ with the divisorial logarithmic structure associated to $X \subseteq \mathscr{X}$ and $\Delta$ with the divisorial logarithmic structure associated to the origin, then
     \begin{enumerate}
        \item The logarithmic Hodge--de Rham spectral sequence  given by
        \begin{center}
            $\mathrm{E}_1^{p,q}=\mathrm{R}^qf_*\Omega^p_{\mathscr{X}/\Delta} \Rightarrow \mathbb{R}^{p+q}f_*\Omega^\bullet_{\mathscr{X}/\Delta}$
        \end{center}
        degenerates on the first page.
        \item For all natural numbers $p$ and $q$, the $\mathscr{O}_\Delta$-modules $\mathrm{R}^qf_*\Omega^p_{\mathscr{X}/\Delta}$ are locally free and their formation commutes with base change.
    \end{enumerate}
\end{proposition}

\begin{preuve}
    The first statement and the local freeness are local for the analytic topology on $\Delta$; by \cite[Proposition 7.8.5.]{Gro63} and \cite[Proposition 7.8.4.]{Gro63}, the commutation with base change follow from local freeness and is hence local for the analytic topology. We can therefore work on a germ of $\Delta$ around the origin; in this case, $\mathrm{E}_1^{p,q}$ corresponds to a finite-type $\mathbb{C}$-algebra. By faithful exactness of the completion, we can suppose that we work over the completion of the germ $(\Delta, 0)$. By the analytic version of the theorem on formal functions \cite[Satz 3.1.]{Bin78}, one can suppose that we work over $\mathrm{Spec}(\mathbb{N} \xrightarrow{1 \mapsto t} \mathbb{C}[t]/(t^{n+1}))$ and in that case, the result is known, see \cite[Theorem 1.10.]{FRF21}.
\end{preuve}

% \begin{corollaire}
%     Let $f : \mathscr{X} \to S$ be a proper log smooth morphism between fine and saturated logarithmic complex analytic spaces such that there exists a morphism $p : S \to C$ of log analytic spaces (here, $C$ is a Riemann surface equipped with the divisorial logarithmic structure associated to a point $x$) admitting locally around $x$ a chart of the form
%     \begin{center}
%         \begin{tikzcd}
% 	{\mathbb{N}} & R \\
% 	{\mathbb{N}} & {\Gamma(\Delta, \mathscr{O}_{\Delta})}
% 	\arrow[from=1-1, to=1-2]
% 	\arrow[Rightarrow, no head, from=1-1, to=2-1]
% 	\arrow[from=2-1, to=2-2]
% 	\arrow["{p^\sharp}", from=2-2, to=1-2]
% \end{tikzcd}
%     \end{center}
%     where the arrow at the bottom of the square maps $1$ to the coordinate function $t$ on the disk $\Delta$. Then
%      \begin{enumerate}
%         \item The logarithmic Hodge--de Rham spectral sequence  given by
%         \begin{center}
%             $\mathrm{E}_1^{p,q}=\mathrm{R}^qf_*\Omega^p_{\mathscr{X}/S} \Rightarrow \mathbb{R}^{p+q}f_*\Omega^\bullet_{\mathscr{X}/S}$
%         \end{center}
%         degenerates on the first page.
%         \item For all natural numbers $p$ and $q$, the $\mathscr{O}_S$-modules $\mathrm{R}^qf_*\Omega^p_{\mathscr{X}/S}$ are locally free and their formation commute with base change.
%     \end{enumerate}
% \end{corollaire}

% \begin{preuve}
%     Since all the statements are local for the analytic topology on $S$, we .  
% \end{preuve}

Now, we study the topological constraints on the nearby fiber $X_t$ of a semistable degeneration of $X$.

\begin{theoreme} \label{determination_possible_nearby_fiber}
	The nearby fiber of any semistable degeneration $f : \mathscr{X} \to \Delta$ of $X$ must be an abelian surface.
\end{theoreme}

\begin{preuve}
    Let us first prove that the dualizing sheaves of $X_t$ are trivial for $t$ in the small disk contained in $\Delta$. By Corollary \ref{degeneration_hodge_de_rham_lissage}, the $f_*\omega_{\mathscr{X}/\Delta}$ is a locally free $\mathscr{O}_{\Delta}$-module, which has rank one because it is true on the central fiber $X$. We can consider a global section $\sigma$ of $f_*\omega_{\mathscr{X}/\Delta}$, its non-vanishing locus is an open subset of $\Delta$ which contains the origin since $X$ is $K$-trivial; by openness, this locus contains a small disk centered at the origin which means that for $t$ close to the origin, $\sigma$ does not vanish, which implies that $X_t$ is $K$-trivial for $t \in \Delta$ close enough to the origin.

    Since $X$ is a proper complex algebraic variety which is normal crossing, \cite[Theorem 1.9.]{FRF21} yields a logarithmic analogue of the Hodge-de Rham spectral sequence, namely:
	\begin{center}
		$\mathrm{E}_1^{p,q}=\mathrm{H}^q(X, \Omega^p_{X/S_0}) \Longrightarrow \mathbb{H}^{p+q}(X, \Omega_{X/S_0}^\bullet)$
	\end{center}
	which degenerates on the first page. 
 
 % Therefore, we have (non-canonical) isomorphisms:
	% \begin{itemize}
	% 	\item $\mathbb{H}^0(X, \Omega_{X/S_0}^\bullet) \cong \mathrm{H}^0(X, \mathscr{O}_X)\cong\mathbb{C}$.
	% 	\item $\mathbb{H}^1(X, \Omega_{X/S_0}^\bullet) \cong \mathrm{H}^0(X, \Omega^1_{X/S_0}) \oplus \mathrm{H}^1(X, \mathscr{O}_X) \cong \mathrm{H}^0(X, \Omega^1_{X/S_0}) \oplus \mathbb{C}^2$.
	% 	\item $\mathbb{H}^2(X, \Omega_{X/S_0}^\bullet) \cong \mathrm{H}^2(X, \mathscr{O}_X) \oplus \mathrm{H}^1(X, \Omega^1_{X/S_0}) \oplus \mathrm{H}^0(X, \omega_X) \cong \mathrm{H}^1(X, \Omega^1_{X/S_0}) \oplus \mathbb{C}^2$.
	% 	\item $\mathbb{H}^3(X, \Omega_{X/S_0}^\bullet) \cong \mathrm{H}^2(X, \Omega^1_{X/S_0}) \oplus \mathrm{H}^1(X, \omega_X) \cong  \mathrm{H}^2(X, \Omega^1_{X/S_0}) \oplus \mathbb{C}^2$.
	% 	\item $\mathbb{H}^4(X, \Omega_{X/S_0}^\bullet) \cong \mathrm{H}^2(X, \omega_X) \cong \mathbb{C}$.
	% \end{itemize}
	Using Lemma \ref{calcul h^1,0 log} and Proposition \ref{cohomology_nearby_fiber}, we get that
    \begin{center}
        $b_1(X_t)=4$
    \end{center}
	and the $K$-triviality of the nearby fiber implies that the latter has Kodaira dimension zero, the result therefore follows from Enriques-Kodaira classification \cite[Chapter VI., Theorem 1.1.]{BHPV04} of compact complex surfaces.
\end{preuve}

From the smoothability of $X$, we can compute all its logarithmic Hodge numbers.

\begin{corollaire} \label{computation_tgt_space_log_smooth_deformations}
    We have $h^1(X, \Omega^1_{X/S_0})=4$ and $h^2(X, \Omega^1_{X/S_0})=2$. In particular $t_{\mathrm{LD}_X}$ is four-dimensional.
\end{corollaire}

\begin{preuve}
    By Lemma \ref{calcul h^1,0 log}, we know that $h^0(X, \Omega^1_{X/S_0})=2$. Since $X$ admits a smoothing by Theorem \ref{smoothability_result}, Proposition \ref{cohomology_nearby_fiber} implies that $b_1(X_t)=4$ and since the nearby fiber $X_t$ is smooth, Poincar\'{e} duality implies that
    \begin{center}
        $4=b_3(X_t)=\underbrace{h^1(X, \omega_X)}_{=2}+h^2(X, \Omega^1_{X/S_0})$,
    \end{center}
    so $h^2(X, \Omega^1_{X/S_0})=2$. Now, the short exact sequence of Proposition \ref{ses_log1-forms} and the additivity of the Euler characteristic implies that $\chi(X, \Omega^1_{X/S_0})=0$, which in turn implies that $h^1(X, \Omega^1_{X/S_0})=4$.
    
    The second claim follows using Serre duality, the fact that $X$ is log Calabi-Yau and that the tangent space of the functor $\mathrm{LD}_X$ is canonically isomorphic to $\mathrm{H}^1(X, \mathscr{T}_{X/S_0})$.
\end{preuve}

\subsection{Kuranishi germs.}

In this subsection, we use the information we have obtained on the functor $\mathrm{LD}_X$ to get information on the functor $\mathrm{Def}_{X, \mathbb{C}}^{\mathrm{lt}}$.

\begin{proposition}
    The differential $\mathrm{d}(F|_{\mathbf{Art}_{\mathbb{C}}}) : \mathrm{H}^1(X, \mathscr{T}_{X/S_0}) \to \mathrm{H}^1(X, \mathscr{T}_{\underline{X}})$ of the natural transformation $F|_{\mathrm{Art}_{\mathbb{C}}}$ is surjective.
\end{proposition}

\begin{preuve}
    The long exact sequence in cohomology associated to the short exact sequence of the proposition \ref{comparison_classical_tangent_sheaf_with_log_tangent_sheaf} is given by
    \begin{center}
        $0 \to \mathrm{H}^0(X, \mathscr{T}_{X/S_0}) \to \mathrm{H}^0(X, \mathscr{T}_{\underline{X}}) \to \mathrm{H}^0(X, \mathscr{E}xt^1_{\mathscr{O}_X}(j_*\mathscr{O}_D, \mathscr{O}_X)) \to \mathrm{H}^1(X, \mathscr{T}_{X/S_0}) \xrightarrow{\mathrm{d}(F|_{\mathbf{Art}_{\mathbb{C}}})} \mathrm{H}^1(X, \mathscr{T}_{\underline{X}}) \to \dots$
    \end{center}
    By Proposition \ref{comp_coh_tangent_sheaf_gluing} and Corollary \ref{computation_tgt_space_log_smooth_deformations}, we have that $h^0(X, \mathscr{T}_{X/S_0})=h^0(X, \mathscr{T}_{\underline{X}})$ so the first map of the long exact sequence is an isomorphism. By exactness, the next arrow has to be zero and therefore, the linear map
    \begin{center}
        $\mathrm{H}^0(X, \mathscr{E}xt^1_{\mathscr{O}_X}(j_*\mathscr{O}_D, \mathscr{O}_X)) \to \mathrm{H}^1(X, \mathscr{T}_{X/S_0})$
    \end{center}
    has to be injective. Let us now compute the rank of this linear map, which is, in this case, the dimension of the source space. We have the local to global spectral sequence
    \begin{center}
        $\mathrm{E}_2^{p,q}=\mathrm{H}^p(X, \mathscr{E}xt^q_{\mathscr{O}_X}(j_*\mathscr{O}_D, \mathscr{O}_X)) \Rightarrow \mathrm{Ext}^{p+q}_{\mathscr{O}_X}(j_*\mathscr{O}_D, \mathscr{O}_X)$.
    \end{center}
    For $p+q=1$, the term that potentially contribute are $\mathrm{H}^0(X, \mathscr{E}xt^1_{\mathscr{O}_X}(j_*\mathscr{O}_D, \mathscr{O}_X))$ and $\mathrm{H}^1(X, \mathscr{H}om_{\mathscr{O}_X}(j_*\mathscr{O}_D, \mathscr{O}_X))$ but the second is zero since $j_*\mathscr{O}_D$ is a torsion $\mathscr{O}_X$-module. This implies that
    \begin{center}
        $\mathrm{H}^0(X, \mathscr{E}xt^1_{\mathscr{O}_X}(j_*\mathscr{O}_D, \mathscr{O}_X)) \cong \mathrm{Ext}^1_{\mathscr{O}_X}(j_*\mathscr{O}_D, \mathscr{O}_X) \cong \mathrm{H}^1(D, \mathscr{O}_D)^*$
    \end{center}
    where the last isomorphism is given by Serre duality for coherent sheaves. In particular
    \begin{center}
        $h^0(X, \mathscr{E}xt^1_{\mathscr{O}_X}(j_*\mathscr{O}_D, \mathscr{O}_X))=1$
    \end{center}
    and since $h^1(X, \mathscr{T}_{X/S_0})=4$ by Corollary \ref{computation_tgt_space_log_smooth_deformations}, we get that
    \begin{center}
        $\mathrm{rk}(\mathrm{d}(F|_{\mathbf{Art}_{\mathbb{C}}}))=h^1(X, \mathscr{T}_{X/S_0})-h^0(X, \mathscr{E}xt^1_{\mathscr{O}_X}(j_*\mathscr{O}_D, \mathscr{O}_X))=3$
    \end{center}
    and since $h^1(X, \mathscr{T}_{\underline{X}})=3$ by Proposition \ref{comp_coh_tangent_sheaf_gluing}, it implies that the map $\mathrm{d}(F|_{\mathbf{Art}_{\mathbb{C}}})$ has to be surjective.
\end{preuve}

\begin{corollaire} \label{smooth_locally_trivial_kuranishi_space}
    The functor $\mathrm{Def}_{X, \mathbb{C}}^{\mathrm{lt}}$ is smooth. In particular, the germ $(\mathrm{Def}^{\mathrm{lt}}(X), 0)$ of the locally trivial Kuranishi space of $X$ is smooth of dimension 3.
\end{corollaire}

\begin{preuve}
    By Theorem \ref{functor_log_def_smooth_gluing}, the functor $\mathrm{LD}_X$ is smooth, this implies that its restriction to the subcategory $\mathbf{Art}_{\mathbb{C}}$ is smooth as well. Since $X$ is a proper, fine and saturated (in particular integral) logarithmic scheme, the functor $\mathrm{LD}_X$ has a prorepresentable hull by \cite[Theorem 8.7.]{FKat96}. This implies that its restriction to the subcategory $\mathbf{Art}_{\mathbb{C}}$ has a prorepresentable hull as well. Indeed, for all extension $A' \twoheadrightarrow A$ and small extension $A'' \twoheadrightarrow A$ in $\mathbf{Art}_\mathbb{C}$, the canonical map
    \begin{center}
        $\alpha : \mathrm{LD}_X|_{\mathbf{Art}_\mathbb{C}}(A' \times_A A'') \to \mathrm{LD}_X|_{\mathbf{Art}_\mathbb{C}}(A') \times_{\mathrm{LD}_X|_{\mathbf{Art}_\mathbb{C}}(A)} \mathrm{LD}_X|_{\mathbf{Art}_\mathbb{C}}(A'')$
    \end{center}
    is surjective by Schlessinger criterion (\cite[Theorem 2.3.2.]{Ser06}) applied to the smooth functor $\mathrm{LD}_X$ --- note that the fiber product on the left is (set-theoretically) the same whether it is computed in the category $\mathbf{Art}_\mathbb{C}$ or $\mathbf{Art}_{\mathbb{C}\llbracket t \rrbracket}$; for the same reason, the map $\alpha$ is bijective when $A=\mathbb{C}$ and $A''=\mathbb{C}[\varepsilon]$; moreover, the tangent space of $\mathrm{LD}_X|_{\mathbf{Art}_\mathbb{C}}$ is the same to that of $\mathrm{LD}_X$ so it is finite dimensional; Schlessinger criterion therefore implies that $\mathrm{LD}_X|_{\mathbf{Art}_\mathbb{C}}$ has a prorepresentable hull. By \cite[Proposition 2.4.6.]{Ser06}, $\mathrm{Def}_{X, \mathbb{C}}^{\mathrm{lt}}$ has an obstruction space so \cite[Corollary 2.3.7.]{Ser06} implies that the functor $\mathrm{Def}_{X, \mathbb{C}}^{\mathrm{lt}}$ is smooth. In particular, the germ $(\mathrm{Def}^{\mathrm{lt}}(X), 0)$ is smooth and its dimension is equal to $h^1(X, \mathscr{T}_{\underline{X}})=3$.
\end{preuve}

We now move to the Kuranishi germ of $X$.

\begin{proposition}
    The differential $\mathrm{d}F : \mathrm{H}^1(X, \mathscr{T}_{X/S_0}) \to \mathrm{Ext}^1_{\mathscr{O}_X}(\Omega^1_{\underline{X}}, \mathscr{O}_X)$ is an isomorphism.
\end{proposition}

\begin{preuve}
    The image of the logarithmic non-locally trivial first-order deformation $[\mathscr{X} \to S_0]$ constructed in the theorem \ref{log_defo_first_order_gluing} by $\mathrm{d}(F|_{\mathbf{Art}_{\mathbb{C}}})$ is zero. A way of rephrase Theorem \ref{log_defo_first_order_gluing} would be to say that $[\mathscr{X} \to S_{1,1}] \in \mathrm{H}^1(X, \mathscr{T}_{X/S_0})$ is mapped to $[\underline{\mathscr{X}} \to \mathrm{Spec}(\mathbb{C}[\varepsilon])] \in \mathrm{Ext}^1_{\mathscr{O}_X}(\Omega^1_{\underline{X}}, \mathscr{O}_X)$ by $\mathrm{d}F$. In particular, the linear map $\mathrm{d}F$ is surjective but since its source and its target have the same dimension, $\mathrm{d}F$ is an isomorphism. 
 \end{preuve}

We now can deduce a non-locally trivial version of Corollary \ref{smooth_locally_trivial_kuranishi_space}.

\begin{corollaire} \label{smooth_kuranishi_space}
    The functor $\mathrm{Def}_{X, \mathbb{C}}$ is smooth. In particular, the germ $(\mathrm{Def}(X), 0)$ of the Kuranishi space of $X$ is smooth of dimension 4.
\end{corollaire}

\begin{preuve}
    Using the same arguments as in the proof of Corollary \ref{smooth_locally_trivial_kuranishi_space}, one obtain that the functor $\mathrm{Def}_{X, \mathbb{C}\llbracket t \rrbracket \rbrace}$ is smooth. Since the functor $\mathrm{Def}_{X, \mathbb{C}}$ is the restriction of the functor $\mathrm{Def}_{X, \mathbb{C}\llbracket t \rrbracket}$ to the subcategory $\mathbf{Art}_{\mathbb{C}}$, it is also smooth as we have seen at the beginning of the proof of Corollary \ref{smooth_locally_trivial_kuranishi_space}. The computation of the dimension of the germ comes from the isomorphism
    \begin{center}
        $\mathrm{T}_0\mathrm{Def}(X) \cong \mathrm{Ext}^1_{\mathscr{O}_X}(\Omega^1_{\underline{X}}, \mathscr{O}_X)$
    \end{center}
    and the fact that the dimension of the right hand side is equal to 4 in virtue of Corollary \ref{comp_coh_tgt_gluing}.
\end{preuve}

\section{Appendix.}

In this appendix, whenever the letter $X$ appears, unless stated otherwise, it has nothing to do with the variety $X$ defined in Section 1.

\subsection{Logarithmic geometry.}

The interested reader can read \cite{Og18} for more background on logarithmic geometry. We prove some results that are mentioned with no proof in this book and also prove a (partial) toroidal characterization of ideal log smoothness \`{a} la Kato.

\subsubsection{The category of logarithmic schemes.}

Let us start with the main definition of logarithmic geometry.

\begin{definition}
    Let $(X, \mathscr{O}_X)$ be a ringed space. 
    A \emph{prelogarithmic structure} on $(X, \mathscr{O}_X)$ is a pair $(\mathscr{M}_X, \alpha)$ where $\mathscr{M}_X$ is a sheaf of monoids on $X$ and $\alpha : \mathscr{M}_X \to (\mathscr{O}_X, \cdot)$ is a morphism of sheaves of monoids. If in addition $\alpha$ induces an isomorphism between $\mathscr{O}_X^\times$ and $\mathscr{M}_X^\times:=\alpha_X^{-1}\mathscr{O}_X^\times$, we say that $(\mathscr{M}_X, \alpha)$ is a \emph{logarithmic structure} on $(X, \mathscr{O}_X)$.
\end{definition}

\begin{remarque}
    This definition also makes sense when $(X, \mathscr{O}_X)$ is a ringed site, we will use it. It is usual to forget the map $\alpha$ in the notation and to say that $\mathscr{M}_X$ is a logarithmic structure on $X$.
\end{remarque}

Since $\mathscr{O}_X^\times$ is part of any logarithmic structure, we might want to get rid of it.

\begin{definition}
    Let $(X, \mathscr{O}_X)$ be a ringed space equipped with a logarithmic structure $\mathscr{M}$. \\
    The \emph{ghost sheaf} $\overline{\mathscr{M}}$ of $\mathscr{M}$ is the sheaf of monoids on $X$ defined as $\overline{\mathscr{M}}:=\mathscr{M}/\mathscr{M}^\times$.
\end{definition}

We can define morphisms between (pre)logarithmic structures.

\begin{definition}
    Let $\alpha_i : \mathscr{M}_i \to \mathscr{O}_X$ be prelogarithmic structures on $X$ for $i \in \lbrace 1,2 \rbrace$. \\
    A morphism $\varphi : \alpha_1 \to \alpha_2$ is a morphism of sheaves $\varphi : \mathscr{M}_1 \to \mathscr{M}_2$ fitting in the following commutative triangle:
    \begin{center}
        \begin{tikzcd}
	{\mathscr{M}_1} & {\mathscr{M}_2} \\
	{\mathscr{O}_X}
	\arrow["\varphi", from=1-1, to=1-2]
	\arrow["{\alpha_1}"', from=1-1, to=2-1]
	\arrow["{\alpha_2}", from=1-2, to=2-1]
\end{tikzcd}.
    \end{center}
    We define the category of prelogarithmic structures on $X$ as the category whose objects are prelogarithmic structures on $X$ and whose morphisms are defined above. The category of logarithmic structures on $X$ is the full subcategory of the category of prelogarithmic structures on $X$ whose objects are logarithmic structures on $X$.
\end{definition}

\begin{lemme}
    The forgetful functor from the category of logarithmic structures on $X$ to that of prelogarithmic structures on $X$ has a left adjoint.
\end{lemme}

\begin{remarque}
    This left adjoint can be described explicitely: if $(\mathscr{M}, \alpha)$ is a prelogarithmic structure on $X$, then we define a logarithmic structure $(\mathscr{M}^{\mathrm{log}}, \alpha^{\mathrm{log}})$ called the \emph{logarithmic structure associated to the prelogarithmic structure $(\mathscr{M}, \alpha)$} by: 
    \begin{center}
        $\mathscr{M}^{\mathrm{log}}:=(\mathscr{M} \oplus \mathscr{O}_X^\times)/\mathscr{P}$ where $\mathscr{P}:=\lbrace(x,\alpha(x)^{-1}) \in \mathscr{M} \oplus \mathscr{O}_X^\times \mid x \in \alpha^{-1}(\mathscr{O}_X^\times)\rbrace$
    \end{center}
    and $\alpha^{\mathrm{log}}(x, u):=u \cdot \alpha(x)$.    
\end{remarque}

\begin{definition}
    A \emph{logarithmic scheme} is a scheme equipped with a logarithmic structure. 
\end{definition}

\begin{remarque}
    Usually, we will consider logarithmic structures for the \'{e}tale topology but it will happen that for computational reasons, we consider logarithmic structures for the Zariski topology. Following \cite{Og18}, we will denote by $\underline{X}$ the underlying scheme to a logarithmic scheme $X$.
\end{remarque}

We also have functorial constructions for logarithmic schemes.

\begin{definition}
    Let $X$ and $Y$ be logarithmic schemes and $f : \underline{X} \to \underline{Y}$ be a scheme morphism.
    \begin{itemize}
        \item The \emph{pushforward} $f_*^{\log}\mathscr{M}_X$ is the logarithmic structure on $\underline{Y}$ defined by the cartesian square
        \begin{center}
            \begin{tikzcd}
	{f_*^{\log}\mathscr{M}_X} & {f_*\mathscr{M}_X} \\
	{\mathscr{O}_Y} & {f_*\mathscr{O}_X}
	\arrow[from=1-1, to=1-2]
	\arrow[from=1-1, to=2-1]
	\arrow["\lrcorner"{anchor=center, pos=0.125}, draw=none, from=1-1, to=2-2]
	\arrow["{f_*\alpha_X}", from=1-2, to=2-2]
	\arrow["{f^\sharp}", from=2-1, to=2-2]
\end{tikzcd}
        \end{center}
         \item The \emph{pullback} $f^*_{\log}\mathscr{M}_Y$ is the logarithmic structure on $\underline{X}$ associated to the prelogarithmic structure given by the composition
         \begin{center}
             $f^{-1}\mathscr{M}_Y \to f^{-1}\mathscr{O}_Y \to \mathscr{O}_X$.
         \end{center}
    \end{itemize}
\end{definition}

One can prove the following expected lemma.

\begin{lemme}
    Let $X$ and $Y$ be logarithmic schemes and $f : \underline{X} \to \underline{Y}$ be a scheme morphism. \\
    The pair $(f^*_{\log}, f_*^{\log})$ is a pair of adjoint functors.
\end{lemme}

We recall a lemma that is often useful, we refer to \cite[III., Remark 1.1.6.]{Og18} for a proof.

\begin{lemme} \label{ghost_sheaf_pullback}
    If $f : X \to Y$ is a scheme morphism and if $\alpha : \mathscr{M}_Y \to \mathscr{O}_Y$ is a logarithmic structure on $Y$, then the ghost sheaf of the logarithmic structure $f^*_{\log}\mathscr{M}_Y$ is canonically isomorphic to $f^{-1}\overline{\mathscr{M}}_Y$.
\end{lemme}

We now define the notion of morphism between logarithmic schemes.

\begin{definition} \label{definition_strict_morphism}
    A morphism between two logarithmic schemes $X$ and $Y$ is a pair $(f, f^\flat)$ where $f : \underline{X} \to \underline{Y}$ is a scheme morphism and $f^\flat : \mathscr{M}_Y \to f_* \mathscr{M}_X$ is a morphism of sheaves of monoids making the square
	\begin{center}
		\begin{tikzcd}
\mathscr{M}_Y \arrow[d, "\alpha_Y"'] \arrow[r, "f^\flat"] & f_* \mathscr{M}_X \arrow[d, "\alpha_X \circ f"] \\
\mathscr{O}_Y \arrow[r, "f^\sharp"']                      & f_*\mathscr{O}_X                 
\end{tikzcd}
	\end{center}
 commutative. A morphism $f : X \to Y$ of logarithmic schemes is said to be \emph{strict} if $f^\flat$ is an isomorphism of sheaves of monoids.
\end{definition}

\begin{remarque}
    Note that in that case, the logarithmic structure on $X$ is isomorphic to the pullback of that on $Y$.
\end{remarque}

Let us give some examples that will be useful.

\begin{paragraphe} \label{example_log_point}
    Let $k$ be a field. A \emph{logarithmic point} is a logarithmic scheme $S$ which is isomorphic to $\mathrm{Spec}(k)$ equipped with the logarithmic structure associated to the prelogarithmic structure defined by
    \begin{align*}
    Q&\rightarrow k\\
    m&\mapsto \mathbf{1}_{Q^\times}(m)
\end{align*}
where $Q$ is a monoid that has a single invertible element --- we say that $Q$ is a \emph{sharp} monoid --- and where $\mathbf{1}_{Q^\times}$ stands for the characteristic function of the set $Q^\times$. When $Q=\mathbb{N}$, the logarithmic scheme that we obtain is called the \emph{standard} logarithmic point, we will denote it by $S_0$. If $k$ is algebraically closed, this is the only logarithmic point by \cite[Example 2.5. (2)]{KK88}.
\end{paragraphe}

\begin{paragraphe}
    Let $D \subseteq X$ be a reduced divisor in a normal scheme and let $j : X \setminus D \hookrightarrow X$ be the open inclusion of the complement of $D$. The \emph{divisorial logarithmic structure on $X$ associated to $D$} is that for which the sheaf of monoids is given by $j_*\mathscr{O}_{X \setminus D}^\times \cap \mathscr{O}_X$ and $\alpha$ by the inclusion. It is a logarithmic structure on $X$.
\end{paragraphe}

The following example will be of specific importence because toric geometry is somehow the local model for logarithmic geometry.

\begin{paragraphe}
    Let $P$ be a monoid and $R$ be a ring. We can form the $R$-algebra $R[P]$ of the monoid $P$. \\
    We define the logarithmic scheme whose underlying scheme is the affine toric variety $\mathrm{Spec}(R[P])$ and whose logarithmic structure is the logarithmic structure associated to the prelogarithmic structure induced by the canonical map $P \to R[P]$. This is sometimes called the \emph{monoidal logarithmic structure associated to $P$} and the resulting logarithmic scheme is usually written $A_P$ when $R=\mathbb{Z}$.
\end{paragraphe}

We now introduce some objects that will allows us to locally control logarithmic structures. Right before that, let us introduce some notions related to monoids.

\begin{definition} \label{definition_integral_ff_monoid}
    We say that a monoid $M$ is
    \begin{itemize}
        \item \emph{finitely generated} if there exists a surjective morphism of monoids $\mathbb{N}^r \twoheadrightarrow M$.
        \item \emph{u-integral} if for all $m,m' \in M$ and all $u \in M^\star$, $m+u=m'+u$ implies that $m=m'$.
        \item \emph{integral} if for all $m,m' \in M$ and all $n \in M$, $m+n=m'+n$ implies that $m=m'$.
    \end{itemize}
\end{definition}

\begin{definition}
    Let $f : X \to Y$ be a morphism of logarithmic schemes.
    \begin{itemize}
        \item A \emph{chart} for $\mathscr{M}_X$ is a morphism of sheaves of monoids $\underline{P}_X \to \mathscr{M}$ that induces a isomorphism of logarithmic structures between the monoidal structure defined by $P$ and and $\mathscr{M}_X$.
        \item A \emph{chart} for $f$ is a triple $(\underline{P}_X \xrightarrow{q_X} \mathscr{M}_X, \underline{Q}_Y \xrightarrow{q_Y} \mathscr{M}_Y, Q \xrightarrow{a} P)$ where $q_X$ is a chart for $X$, $q_Y$ is a chart for $Y$ and $a$ is a morphism of monoids making the square
        \begin{center}
            \begin{tikzcd}
	{\underline{Q}_X} & {\underline{P}_X} \\
	{f^{-1}\mathscr{M}_Y} & {\mathscr{M}_X}
	\arrow["a", from=1-1, to=1-2]
	\arrow["{q_X}"', from=1-1, to=2-1]
	\arrow["{q_Y}", from=1-2, to=2-2]
	\arrow["{f^\flat}"', from=2-1, to=2-2]
\end{tikzcd}
        \end{center}
        commutative.
    \end{itemize}
\end{definition}

\begin{remarque}
    A chart for a logarithmic structure controls the behaviour of the ghost sheaf but the converse might not be true, namely when the logarithmic structure is not $u$-integral --- see \cite[II., Proposition 2.1.4.]{Og18}.
\end{remarque}

As in classical scheme theory, one can define some finiteness conditions for log schemes.

\begin{definition}
   Let $X$ be a logarithmic scheme. \\
   We say that $X$ is \emph{fine} (resp. \emph{saturated}) if it locally admits a chart whose source is a finitely generated and integral (resp. saturated) monoid.
\end{definition}

Let us finish this subsubsection by introducing a type of logarithmic structures that appears in semistable degenerations of smooth varieties.

\begin{definition}
    Let $X$ be a normal crossing complex variety of dimension $n$. \\
    A logarithmic structure $\mathscr{M}$ on $X$ is said to be of \emph{embedding type} if there exists an \'{e}tale cover $(h_i : U_i \to X)_{i \in I}$ of $X$ such that, for each $i \in I$, there exists a closed immersion $\iota_i : U_i \hookrightarrow V_i$ fitting in the commutative square
    \begin{center}
    \begin{tikzcd}
	{U_i} & {V_i} \\
	{V(x_0\dots x_r)} & {\mathbb{A}^{n+1}_{\mathbb{C}}}
	\arrow["{\iota_i}", hook, from=1-1, to=1-2]
	\arrow[from=1-1, to=2-1]
	\arrow[from=1-2, to=2-2]
	\arrow[hook, from=2-1, to=2-2]
    \end{tikzcd}
    \end{center}
     for some $0 \leq r \leq n$, where the two vertical arrows are \'{e}tale.
\end{definition}

\begin{remarque}
    By the remarks following \cite[Definition 11.4.]{FKat96}, we can describe \'{e}tale locally a chart in the case of a logarithmic structure of embedding type: with the notations above, it is given by $e_j \in \mathbb{N}^r \mapsto z_j \in \mathscr{O}_{U_i}$ where $z_j$ is the pullback to $U_i$ of the $j$-th coordinate function on $\mathbb{A}^{n+1}$.
\end{remarque}
For a logarithmic structure of embedding type, one can compute the ghost sheaf.

\begin{paragraphe} \label{construction_order_morphism}
    Let $X$ be a normal crossing complex algebraic variety equipped with a logarithmic structure of embedding type with (reduced) singular locus $D$. We define the morphism $\varphi$ of \'{e}tale sheaves of monoids on $X$ as follows: if $U \to X$ is \'{e}tale, we consider the morphism of monoids defined by
    \begin{center}
	$\varphi_U$ : \begin{tabular}{ccc}
	$\mathscr{M}(U \xrightarrow{h} X)$ & $\rightarrow$ & $(\nu_*\underline{\mathbb{N}}_{\tilde{X}})(U \xrightarrow{h} X)$ \\
	$m$ & $\mapsto$ & $\mathrm{ord}_{U \times_X D \times_X \tilde{X}}(\alpha_X(m) \circ \nu \circ \tilde{h})$
	
	\end{tabular}
	\end{center}
where $\tilde{h} : U \times_X \tilde{X} \to \tilde{X}$ is the base change of $h$ and $\nu : \tilde{X} \to X$ is the normalization. This is well defined because $U \times_X D \times_X \tilde{X}$ is a union of finitely many smooth divisors in $\tilde{X}$.
\end{paragraphe}

\begin{lemme} \label{ses_ghost_sheaf_log_str_X}
    If $X$ is a normal crossing complex variety endowed with a logarithmic structure of embedding type, we have a short exact sequence of \'{e}tale sheaves of abelian groups
    \begin{center}
        $0 \to \mathscr{O}_X^\times \to \mathscr{M}_X^{\mathrm{gp}} \xrightarrow{\varphi^{\mathrm{gp}}} \nu_*\underline{\mathbb{Z}}_{\tilde{X}} \to 0$
    \end{center}
    where $\nu : \tilde{X} \to X$ is the normalization.
\end{lemme}

\begin{preuve}
    Let us prove the exactness at each step separately. We can work \'{e}tale locally and therefore suppose that $\mathscr{M}_X$ is the divisorial logarithmic structure on $\mathbb{A}^{n+1}$ associated to the normal crossing divisor $X=V(x_1 \dots x_r)$ restricted to this divisor. In that case, we have $\tilde{X}=\coprod_{i=1}^r V(x_0 \dots \hat{x_i} \dots x_r)$ and therefore the stalk of $\nu_*\underline{\mathbb{N}}_{\tilde{X}}$ at any (geometric) point is equal to $\mathbb{N}^r$ and $\mathscr{M}_X  \cong (j_*\mathscr{O}_{\mathbb{A}^{n+1} \setminus X}^\times \cap \mathscr{O}_{\mathbb{A}^{n+1}})|_X$.
    \begin{itemize}
        \item \textbf{Exactness on the left:} the map is the composition of the two injections $\mathscr{O}_X^\times \hookrightarrow \mathscr{M}_X \hookrightarrow \mathscr{M}_X^{\mathrm{gp}}$ where the second map is injective because $\mathscr{M}_X$ is an integral logarithmic structure.
        \item \textbf{Exactness in the middle:} a section of $\mathscr{M}_X^{\mathrm{gp}}$ is nowhere vanishing if and only if it does not vanish on the divisor $X$ which is equivalent to the non-vanishing on each of the irreducible components of $X$, hence the exactness in the middle.
        \item \textbf{Exactness on the right:} if $(p_1, \dots, p_r) \in \mathbb{N}^r$, the section $x_1^{p_1} \dots x_r^{p_r}$ vanishes at order $p_i$ along $V(x_0 \dots \Hat{x_i} \dots x_r)$ so the map on the right is surjective.
    \end{itemize}
\end{preuve}

This allows us to define the main type of logarithmic structure we are interested in.

\begin{definition} \label{definition_log_structure_semistable_type}
    Let $X$ be a normal crossing complex variety. \\
    A logarithmic structure \emph{of semistable type} on $X$ is a logarithmic structure of embedding type $\mathscr{M}$ for which there exists a morphism of \'{e}tale sheaves of monoids $\psi : \underline{\mathbb{Z}}_X \to \mathscr{M}^{\mathrm{gp}}$ such that the triangle
    \begin{center}
        \begin{tikzcd}
	{\mathscr{M}^{\mathrm{gp}}} & {\nu_*\underline{\mathbb{Z}}_{\tilde{X}}} \\
	& {\underline{\mathbb{Z}}_X}
	\arrow["{\varphi^{\mathrm{gp}}}", from=1-1, to=1-2]
	\arrow["\psi", from=2-2, to=1-1]
	\arrow["\delta", from=2-2, to=1-2]
\end{tikzcd}
    \end{center}
    is commutative, where $\delta$ is the diagonal morphism and $\varphi^{\mathrm{gp}}$ was constructed in \ref{construction_order_morphism}.
\end{definition}

\subsubsection{\'{E}tale \emph{vs.} Zariski logarithmic structures.}

In this subsubsection, we discuss the difference between logarithmic structure for the \'{e}tale topology and those for the Zariski topology. Let $X$ be a scheme and let $\eta : X_{\text{\'{e}t}} \to X_{\mathrm{Zar}}$ be the canonical morphism of sites.

\begin{paragraphe}
    If $\mathscr{M}$ is a logarithmic structure on $X_{\text{\'{e}t}}$, we define the logarithmic strucure $\eta_*^{\log}\mathscr{M}$ on $X_{\mathrm{Zar}}$ as the restriction (any open immersion is an \'{e}tale morphism) of $\mathscr{M}$ to the (small) Zariski site of $X$. Conversly, if $\mathscr{M}$ is a logarithmic structure on $X_{\mathrm{Zar}}$, we define the logarithmic structure $\eta^*_{\mathrm{log}}\mathscr{M}$ on $X_{\text{\'{e}t}}$ as follows: if $U \xrightarrow{h} X$ is an \'{e}tale morphism, we put
    \begin{center}
        $(\eta^*_{\log}\mathscr{M})(U \xrightarrow{h} X):=\Gamma(U, h^*_{\log}\mathscr{M})$.
    \end{center}
    The pair $(\eta_{\log}^*, \eta_*^{\log})$ is a pair of adjoint functors.
\end{paragraphe}

We now prove the compatibility of these two operations with the formation of the ghost sheaf.

\begin{lemme}
    If $\alpha : \mathscr{M}_X \to \mathscr{O}_X$ is a logarithmic structure on $X_{\mathrm{Zar}}$, then the canonical morphism $\overline{\eta_{\log}^*\mathscr{M}} \to \eta^{-1}\overline{\mathscr{M}}$ is an isomorphism.
\end{lemme}

\begin{preuve}
    Let us recall a few things before starting the proof. First, if $\mathscr{F}$ is a sheaf on $X_{\mathrm{Zar}}$, $\eta^{-1}\mathscr{F}$ is the (\'{e}tale) sheafification of the presheaf $(\eta^{-1}\mathscr{F})^\mathrm{p}$ defined by $(U \xrightarrow{h} X) \mapsto \Gamma(U, h^{-1}\mathscr{F})$; second, the canonical morphism of the lemma is constructed as follows: if $U \xrightarrow{h} X$ is an \'{e}tale map, $h^{-1}\overline{\mathscr{M}}$ is the ghost sheaf of $h^*_{\log}\mathscr{M}$ by Lemma \ref{ghost_sheaf_pullback}, therefore we get a morphism of \'{e}tale sheaves of monoids
    \begin{center}
        $\eta^*_{\log}\mathscr{M} \to \eta^{-1}\overline{\mathscr{M}}$.
    \end{center}
    Since this morphism maps all the sections of $\mathscr{O}_{X_{\text{\'{e}t}}}^\times$ to zero in $\eta^{-1}\overline{\mathscr{M}}$, we get the sought-for canonical morphism 
    \begin{center}
        $\overline{\eta_{\log}^*\mathscr{M}} \to \eta^{-1}\overline{\mathscr{M}}$.
    \end{center}
    Let us now prove the lemma. We will prove that this morphism is an isomorphism at each stalk and since the sheafification operation does not change the stalks, we can forget about it. If we take an \'{e}tale neighborhood $h : U \to X$ of a geometric point $\overline{x}$ of $X$, we have, Lemma \ref{ghost_sheaf_pullback} implies that
    \begin{center}
        $\overline{\eta^*_{\log}\mathscr{M}}(U \xrightarrow{h} X)=\Gamma(U, \overline{h^*_{\log}\mathscr{M}}) \cong \Gamma(U, h^{-1}\overline{\mathscr{M}})=\Gamma(U, \overline{\eta_{\log}^*\mathscr{M}})$.
    \end{center}
    Taking the limit over all \'{e}tale neighborhoods of $\overline{x}$ in $X$, we obtain the result.
\end{preuve}

\begin{lemme}
    If $\alpha : \mathscr{M}_X \to \mathscr{O}_X$ is a $u$-integral logarithmic structure on $X_{\text{\'{e}t}}$, then the canonical morphism of Zariski sheaves $\overline{\eta^*_{\log}\mathscr{M}} \to \eta_*\overline{\mathscr{M}}$ is an isomorphism. 
\end{lemme}

\begin{preuve}
    The canonical morphism is obtained as in the case of the previous lemma. It is also enough to check that this morphism is an isomorphism on the stalks. The proof goes then as in \cite[III. Proposition 1.4.1., (2)]{Og18} whose proof can be found at the bottom of the page 284, noting that since we want a statement on the stalk, we can just work with local charts and not a \emph{global} chart as it is done there.
\end{preuve}

We give a criterion to for an \'{e}tale logarithmic structure to come from a Zariski one.

\begin{corollaire}
    Let $\mathscr{M}$ be a $u$-integral logarithmic structure on $X_{\text{\'{e}t}}$. The following are equivalent:
    \begin{enumerate}
        \item The logarithmic structure $\mathscr{M}$ comes from a Zariski logarithmic structure.
        \item The canonical morphism of \'{e}tale sheaves $\eta^{-1}\eta_* \overline{\mathscr{M}} \to \overline{\mathscr{M}}$ is an isomorphism.
    \end{enumerate}
\end{corollaire}

\begin{preuve}
    If $\mathscr{M}:=\eta_*\mathscr{N}$ for some Zariski logarithmic structure $\mathscr{N}$, then we have a commutative diagram of \'{e}tale sheaves of monoids on $X$
    \begin{center}
        \begin{tikzcd}
	{\eta_{\log}^*\mathscr{M}} &&& {\mathscr{N}} \\
	{\overline{\eta_{\log}^*\mathscr{M}}} & {\eta^{-1}\overline{\mathscr{M}}} & {\eta^{-1}\eta_*\overline{\mathscr{N}}} & {\overline{\mathscr{N}}}
	\arrow[from=1-1, to=1-4]
	\arrow[from=1-1, to=2-1]
	\arrow[from=1-4, to=2-4]
	\arrow[Rightarrow, no head, from=2-1, to=2-2]
	\arrow["\cong", from=2-2, to=2-3]
	\arrow[from=2-3, to=2-4]
\end{tikzcd}
    \end{center}
    If $\mathscr{N}$ is defined in the Zariski topology, the upper horizontal map is an isomorphism so, by commutativity of the square above, the morphism $\eta^{-1}\eta_*\overline{\mathscr{N}} \to \overline{\mathscr{N}}$ has to be an isomorphism. Conversly, if $\eta^{-1}\eta_*\overline{\mathscr{N}} \to \overline{\mathscr{N}}$ is an isomorphism, then the lower horizontal map is an isomorphism and since the log structures $\mathscr{N}$ and $\eta_{\log}^*\mathscr{M}$ are $u$-integral, \cite[I., Proposition 4.1.2.]{Og18} implies that the upper horizontal map has to be an isomorphism as well, which implies that $\mathscr{N}$ is defined in the Zariski topology.
\end{preuve}

\begin{corollaire} \label{criterion_zariski_log_str_constant_sheaf}
   Let $\mathscr{N}$ be a $u$-integral log structure on $X_{\text{\'{e}t}}$. \\
   If there exists a Zariski sheaf of monoids $\mathscr{P}$ on $X$ such that $\eta^{-1}\mathscr{P} \cong \overline{\mathscr{N}}$ and such that the map $\mathscr{P} \cong \eta_*\eta^{-1}\mathscr{P}$ is an isomorphism, then $\mathscr{N}$ is defined in the Zariski topology. In particular, if $\mathscr{P}$ is a constant sheaf such that $\eta^{-1}\mathscr{P} \cong \overline{\mathscr{N}}$, then $\mathscr{N}$ is defined in the Zariski topology.
\end{corollaire}

\begin{preuve}
    This is a direct consequence of the corollary above. The addendum follows from the fact that if $\mathscr{F}$ is a constant Zariski sheaf, then the canonical morphism of sheaves $\mathscr{F} \to \eta^*\eta^{-1}\mathscr{F}$ is an isomorphism.
\end{preuve}

\subsubsection{Idealized logarithmic schemes.}

Many definitions from the theory of ideals in commutative rings can be naturally adapted to the setting of monoids, leading to analogous concepts such as \textit{ideals} of a monoid, \textit{prime ideals}, and the \textit{spectrum} of a monoid. This gives rise to geometric objects known as \textit{monoschemes} --- the interested reader can read \cite[I., 1.4.]{Og18}. However, unlike in classical scheme theory, closed immersions do not correspond to surjective morphisms of structure sheaves (or, equivalently, to (sheaf of) ideals of the source). For this reason, it is sometimes important to explicitly keep track of the (sheaf of) ideal(s) involved, which motivates the definition of an \emph{idealized} logarithmic scheme.

\begin{definition} \label{definition_idealized_log_scheme}
    An \emph{idealized logarithmic scheme} is a triple $(X, \mathscr{M}_X, \mathscr{K}_X)$ where $(X, \mathscr{M}_X)$ is a logarithmic scheme and $\mathscr{K}_X$ is a sheaf of $\mathscr{M}_X$-ideals such that $\alpha_X(\mathscr{K}_X)=\lbrace 0 \rbrace$.
\end{definition}

\begin{definition}
    A morphism $f : X \to Y$ between two idealized logarithmic schemes is a morphism of logarithmic schemes such that $f^\flat$ maps $f^{-1}\mathscr{K}_Y$ to $\mathscr{K}_X$. We say that $f$ is \emph{ideally strict} if $\mathscr{K}_X$ is generated by $f^{-1}\mathscr{K}_Y$.
\end{definition}

We now introduce the idealized variant of the monoidal logarithmic scheme.

\begin{paragraphe}
    Let $(M, I)$ be an idealized monoid such that either $I$ is a proper ideal of $M$ or $M=0$\footnote{In \cite{Og18}, this is called an \emph{acceptable} idealized monoid.}. If $R$ is a ring, $R[I]$ is an ideal of $R[M]$ and we consider the $R$-algebra
    \begin{center}
        $R[M,I]:=R[M]/R[I]$.
    \end{center}
    This $R$-algebra is called the \emph{monoid algebra of the idealized monoid $(M,I)$}. We denote by $A_{M,I}$ the idealized logarithmic scheme whose underlying scheme is $\mathrm{Spec}(\mathbb{Z}[M,I])$ and whose logarithmic structure is the logarithmic structure associated to the prelogarithmic structure defined by the morphism of monoids $M \to \mathbb{Z}[M,I]$ and whose sheaf of ideals of monoids is that generated by $I$ in $\mathscr{M}_{A_P}$.
\end{paragraphe}

There is also an idealized version of the notion of chart.

\begin{definition}
    Let $X$ be an idealized logarithmic scheme.\\
    A \emph{chart} for $X$ is a pair $(\beta : \underline{P}_X \to \mathscr{M}_X,K)$ where $\beta$ is a chart for $\mathscr{M}_X$ and $K \subseteq P$ is an ideal such that $\mathscr{K}$ is equal to the sheafification of the presheaf that maps an open subset $U \subseteq X$ to the ideal of $\Gamma(U, \mathscr{M}_X)$ generated by $\beta(K)$.
\end{definition}

\begin{remarque}
    As in the non-idealized case, this can be reformulated more geometrically saying that a chart for an idealized logarithmic scheme $X$ is a morphism of idealized logarithmic schemes $X \to A_{P,K}$ that is strict and ideally strict.    
\end{remarque}

\subsubsection{The cotangent sheaf of a logarithmic scheme and logarithmic smoothness.}

The definition of smoothness for logarithmic schemes is similar to Grothendieck's definition of smoothness as ability to lift morphisms along infinitesimal thickenings.

\begin{definition} \label{definition_log_smoothness}
    Let $f : X \to Y$ be a morphism between fine logarithmic schemes. \\
    We say that $f$ is \emph{log smooth} (resp. \emph{log \'{e}tale}, resp. \emph{log unramified}) if $\underline{f}$ is locally of finite presentation and if, for any commutative diagram of fine logarithmic schemes
    \begin{center}
        \begin{tikzcd}
	{T'} & X \\
	T & Y
	\arrow["{h'}", from=1-1, to=1-2]
	\arrow["i"', from=1-1, to=2-1]
	\arrow["f", from=1-2, to=2-2]
	\arrow["g", dotted, from=2-1, to=1-2]
	\arrow["h"', from=2-1, to=2-2]
\end{tikzcd}
    \end{center}
    where $i : T' \to T$ is a strict thickening of order one (i.e., a closed immersion with ideal sheaf of square zero and that is strict as a morphism of logarithmic schemes), locally on $T$, there exists a (resp. there exists a unique, resp. there exists at most one) morphism of log schemes $g : T \to X$ making the diagram above commutative \footnote{This second property is called \emph{formal log smoothness}, resp. \emph{formal log \'{e}taleness}, resp. \emph{formal log unramification}.}.
\end{definition}

We prove a cancellation lemma for log \'{e}tale morphisms that will be useful in the sequel.

\begin{lemme} \label{cancellation_log_etale}
    Let $X \xrightarrow{f} Y \xrightarrow{g} Z$ be morphisms of logarithmic schemes. \\
    If $g$ is formally log \'{e}tale, then the following two statements are equivalent:
    \begin{enumerate}
        \item $f$ is formally log smooth, resp.  formally log \'{e}tale.
        \item $g \circ f$ is formally log smooth, resp.  formally log \'{e}tale.
    \end{enumerate}
    Moreover, if $g$ is only formally log unramified, $2 \Rightarrow 1$ remains true.
\end{lemme}

\begin{preuve}
    The implication $2 \Rightarrow 1$ direcly follows from the fact that a formally log \'{e}tale morphism is formally log smooth and that formally log smooth morphisms are stable under composition. Let us prove the converse: we consider a commutative square
    \begin{equation} \label{square_infinitesimal_lifting}
        \begin{tikzcd}
	{T'} & X \\
	T & Y
	\arrow["{h}", from=1-1, to=1-2]
	\arrow["i"', hook, from=1-1, to=2-1]
	\arrow["f", from=1-2, to=2-2]
	\arrow["{h'}"', from=2-1, to=2-2]
\end{tikzcd}
    \end{equation}
    where $i$ is a strict thickening of order one. It can be completed to a commutative diagram
    \begin{center}
        \begin{tikzcd}
	{T'} & X & X \\
	T & Y & Z
	\arrow["h", from=1-1, to=1-2]
	\arrow["i"', hook, from=1-1, to=2-1]
	\arrow[Rightarrow, no head, from=1-2, to=1-3]
	\arrow["f", from=1-2, to=2-2]
	\arrow["{g \circ f}", from=1-3, to=2-3]
	\arrow["{\tilde{h}}"{pos=0.3}, shift left, dotted, from=2-1, to=1-3]
	\arrow["{h'}"', from=2-1, to=2-2]
	\arrow["g", from=2-2, to=2-3]
\end{tikzcd}
    \end{center}
    where $\tilde{h}$ is obtained by formal log smoothness of $g \circ f$ and namely satisfies $g \circ f \circ \tilde{h}=g \circ h'$. We now want to prove that the square
    \begin{center}
        \begin{tikzcd}
	{T'} & X \\
	T & Y
	\arrow["h", from=1-1, to=1-2]
	\arrow["i"', hook, from=1-1, to=2-1]
	\arrow["f", from=1-2, to=2-2]
	\arrow["{\tilde{h}}", from=2-1, to=1-2]
	\arrow["{h'}"', from=2-1, to=2-2]
\end{tikzcd}
    \end{center}
    is commutative, the commutation of the upper triangle is already known by construction of the lifting $\tilde{h}$ of $h$, we thus only have to prove that $f \circ \tilde{h} = h'$. This follows from the fact that $f \circ \tilde{h}$ and $h'$ are two liftings of $f \circ h$ to $T$ and since $g$ is formally log \'{e}tale, it implies that $f \circ \tilde{h} = h'$; hence, $f$ is formally log smooth. In fact, to prove that $f \circ \tilde{h} = h'$, we just need that $g$ is formally log unramified so this proves the addemdum as well. The proof in the formally log \'{e}tale case is completely identical, noting that two liftings of $h$ as in the diagram (\ref{square_infinitesimal_lifting}) induces two liftings of $h$ in the same diagram except that $f$ is replaced with $g \circ f$.
\end{preuve}

We say some words about the notion of smoothness of \emph{idealized} logarithmic schemes, also known as ideal logarithmic smoothness.

\begin{paragraphe}
    If $f : X \to Y$ is a morphism of idealized logarithmic schemes, we say that $f$ is \emph{ideally log smooth} (resp. \emph{ideally log \'{e}tale}, resp. \emph{ideally log unramified}) there exists a lift $g$ of $f$ to $T$ as in the definition \ref{definition_log_smoothness} but where $i : T' \hookrightarrow T$ is also ideally strict. In particular, Lemma \ref{cancellation_log_etale} still holds in the idealized setting. Note that if $f$ is log smooth, then $f$ is ideally log smooth.
\end{paragraphe}

We now state the \emph{toroidal criterion for logarithmic smoothness} due to K. Kato (\cite[Theorem 3.5.]{KK88}) that essentially says that log smoothness \footnote{In characteristic zero.} is modelled on smoothness between toroidal varieties.

\begin{theoreme} \label{kato_toroidal_criterion}
    Let $f : X \to Y$ be a morphism between fine logarithmic schemes. If there is a chart $\underline{Q}_Y \xrightarrow{q_Y} \mathscr{M}_Y$ where $Q$ a finitely generated monoid, the following statements are equivalent:
    \begin{enumerate}
        \item The morphism $f$ is log smooth, resp. log \'{e}tale.
        \item Locally on $X$, the chart for $\mathscr{M}_Y$ can be extended to a chart $(\underline{P}_X \xrightarrow{q_X} \mathscr{M}_X, \underline{Q}_Y \xrightarrow{q_Y} \mathscr{M}_Y, Q \xrightarrow{a} P)$ for $f$ such that
        \begin{itemize}
            \item The kernel and the torsion part of the cokernel (resp. the kernel and the cokernel) of $a^{\mathrm{gp}} : Q^{\mathrm{gp}} \to P^{\mathrm{gp}}$ are finite groups of order (seen as a constant function on $X$) invertible on $X$.
            \item The morphism of schemes $\underline{X} \to \underline{Y} \times_{\underline{A_Q}} \underline{A_P}$ induced by $f$ is smooth, resp. \'{e}tale.
        \end{itemize}
    \end{enumerate}
\end{theoreme}

\begin{remarque}
    The first condition ensures that the morphism of logarithmic schemes $A_P \to A_Q$ is log smooth, resp. log \'{e}tale.
\end{remarque}

We now prove a partial idealized version of Kato's toroidal criterion that relies on the classical one.

\begin{theoreme} \label{idealized_toroidal_criterion}
    Let $f : X \to Y$ be a morphism between idealized fine logarithmic schemes. If there is a chart $(Q, J) \to (\mathscr{M}_Y, \mathscr{K}_Y)$ where $Q$ is a finintely generated monoid and if, in addition, locally on $X$, the chart for $(\mathscr{M}_Y, \mathscr{K}_Y)$ can be extended to a chart 
        \begin{center}
            $\big( (P,I) \to (\mathscr{M}_X, \mathscr{K}_X), (Q, J) \to (\mathscr{M}_Y, \mathscr{K}_Y), (Q,J) \xrightarrow{a} (P,I) \big)$
        \end{center} 
        for $f$ such that
        \begin{itemize}
            \item The kernel and the torsion part of the cokernel (resp. the kernel and the cokernel) of $a^{\mathrm{gp}} : Q^{\mathrm{gp}} \to P^{\mathrm{gp}}$ are finite groups of order invertible on $X$
            \item The morphism of schemes $ \tilde{f} : \underline{X} \to \underline{Y} \times_{\underline{A_{Q,J}}} \underline{A_{P,I}}$ induced by $f$ is smooth, resp. \'{e}tale,
        \end{itemize}
        then the morphism $f$ is ideally log smooth, resp. ideally log \'{e}tale.
\end{theoreme}

\begin{preuve}
    As in the proof of the classical toroidal criterion, we have a commutative diagram
    \begin{center}
        \begin{tikzcd}
	X & {X \times_{A_{Q,J}} A_{P,I}} & {A_{P,I}} \\
	& Y & {A_{Q,J}}
	\arrow["{\tilde{f}}", from=1-1, to=1-2]
	\arrow["f"', from=1-1, to=2-2]
	\arrow[from=1-2, to=1-3]
	\arrow[from=1-2, to=2-2]
	\arrow["\lrcorner"{anchor=center, pos=0.125}, draw=none, from=1-2, to=2-3]
	\arrow["h", from=1-3, to=2-3]
	\arrow[from=2-2, to=2-3]
\end{tikzcd} 
    \end{center}
    where the square on the right is cartesian and $\tilde{f}$ is strict and ideally strict. If we prove that $h$ is ideally log smooth (resp. ideally log \'{e}tale), then we are done. Indeed, say that $h$ is ideally log smooth, for example; then, the vertical map on the left of the square is also ideally log smooth as the base change of an ideally log smooth map and since $\tilde{f}$ is strict and ideally strict and that the underlying scheme morphism is smooth, \cite[IV., Variant 3.1.22.]{Og18} and \cite[IV., Proposition 3.1.6.]{Og18} imply that $\tilde{f}$ is ideally log smooth. Therefore, $f$ is ideally log smooth as the composition of two ideally log smooth morphisms.

    Let us now use the first hypothesis on the chart for $f$ to prove that $h$ is ideally log smooth. By \cite[IV., Variant 3.1.21.]{Og18}, we have a commutative square of idealized logarithmic schemes
    \begin{center}
        \begin{tikzcd}
	{A_{P,I}} & {A_{P,\emptyset}} \\
	{A_{Q,J}} & {A_{Q, \emptyset}}
	\arrow["{i_P}", from=1-1, to=1-2]
	\arrow["h"', from=1-1, to=2-1]
	\arrow["{h_{\emptyset}}", from=1-2, to=2-2]
	\arrow["{i_Q}"', from=2-1, to=2-2]
\end{tikzcd}
    \end{center}
    where $i_P$ and $i_Q$ are ideally log \'{e}tale; note that $h_{\emptyset}$ is ideally strict since both on its target and its source, the sheaf of ideals of monoids is the empty one. By \cite[IV., Theorem 3.1.8.]{Og18}, the first hypothesis on the chart of $f$ implies that $h_{\emptyset}$ is log smooth and since it is ideally strict, it is also ideally log smooth. In particular, the commutativity of the square above implies that the composition $i_Q \circ h$ is ideally log smooth. Since $i_Q$ is ideally log \'{e}tale, the (idealized version of the) lemma \ref{cancellation_log_etale} implies that $h$ is ideally log smooth as well, which is exactly what we wanted to prove.
\end{preuve}

Let us treat the example of toric varieties.

\begin{paragraphe}
    Let $n \geq 1$ and let $\sigma$ be a cone in $\mathbb{R}^n$. We consider its associated toric variety $\mathrm{Spec}(\mathbb{Z}[P])$ where $P:=\sigma^\vee \cap \mathbb{Z}^n$. If we take $Q:=\lbrace 0 \rbrace \subseteq P$, then we have an isomorphism of schemes
    \begin{center}
        $\mathrm{Spec}(k) \times_{\mathrm{Spec}(\mathbb{Z}[Q])} \mathrm{Spec}(\mathbb{Z}[P]) \cong \mathrm{Spec}(k[P])$
    \end{center}
    which implies that toric varieties are log smooth over the trivial logarithmic point.
\end{paragraphe}

Let us continue with the example of the central fiber of a semistable degeneration.

\begin{paragraphe} \label{example_log_semistable_degeneration}
    Let $X$ be a normal crossing complex algebraic variety equipped with a logarithmic structure of semistable type. By definition, we have a morphism of sheaves $\underline{\mathbb{Z}}_X \to \mathscr{M}^{\mathrm{gp}}$ and therefore have a morphism of sheaves 
    \begin{center}
        $\underline{\mathbb{Z}}_X \times_{\mathscr{M}^{\mathrm{gp}}} \mathscr{M} \cong \underline{\mathbb{N}}_X \to \mathscr{M}$
    \end{center}
    which induces a morphism of logarithmic schemes $X \to S_0$ called a \emph{logarithmic semistable degeneration}. Let us prove that $f$ is log smooth. For $S_0$, we have the global chart $1 \in \mathbb{N} \mapsto 0 \in \mathbb{C}$ and \'{e}tale locally on $X$, we have seen that we had a chart $e_i \in \mathbb{N}^r \mapsto z_i \in \mathscr{O}_X$. In particular, we have
    \begin{center}
        $\mathrm{Spec}(\mathbb{C}) \times_{\mathrm{Spec}(\mathbb{Z}[t])} \mathrm{Spec}(\mathbb{Z}[t_1, \dots, t_r]) \cong \mathrm{Spec}(\mathbb{C} \otimes_{\mathbb{Z}[t]} \mathbb{Z}[t_1, \dots, t_r]) \cong \mathrm{Spec}(\mathbb{C}[t_1, \dots, t_r]/(t_1 \dots t_r))$
    \end{center}
    and the morphism of schemes in the second hypothesis of Kato's toroidal criterion is simply given \'{e}tale locally by
    \begin{center}
        $\mathrm{Spec}(\mathbb{C}[t_1, \dots, t_n]/(t_1 \dots t_r)) \to \mathrm{Spec}(\mathbb{C}[t_1, \dots, t_r]/(t_1 \dots t_r))$
    \end{center}
    which is indeed smooth and since smoothness is \'{e}tale local, $X \to \mathrm{Spec}(\mathbb{C}[t_1, \dots, t_r]/(t_1 \dots t_r))$ is smooth as well. Since we are in characteristic zero, the first hypothesis of Kato's toroidal criterion is automatically satisfied so $f$ is log smooth.
\end{paragraphe}

Let us now introduce the cotangent sheaf of a morphism between fine logarithmic scheme.

\begin{definition} \label{definition_log_cotangent_sheaf}
    Let $f : X \to Y$ be a morphism of logarithmic schemes. \\
    The \emph{cotangent sheaf of $X$ over $Y$} is the $\mathscr{O}_X$-module $\Omega^1_{X/Y}$ defined by
    \begin{center}
        $\Omega^1_{X/Y}:=\big( \Omega^1_{\underline{X}/\underline{Y}} \oplus ( \mathscr{O}_X \otimes_{\underline{\mathbb{Z}}_X} \mathscr{M}_X^{\mathrm{gp}} \big)/\mathscr{R}$
    \end{center}
    where $\mathscr{R}$ is the sub-$\mathscr{O}_X$-module generated by the following two types of sections:
    \begin{itemize}
        \item $(\mathrm{d}(\alpha_X(m)),-\alpha_X(m) \otimes m)$ for all local section $m$ of $\mathscr{M}_X$
        \item $(0, 1 \otimes f^\flat(n))$ for all local section $n$ of $f^{-1}\mathscr{M}_Y$.
    \end{itemize}
    It comes with two morphisms $\mathrm{d} : \mathscr{O}_X \to \Omega^1_{X/Y}$ and $\mathrm{dlog} : \mathscr{M}_X \to \Omega^1_{X/Y}$. For $p \geq 1$, we define $\Omega^p_{X/Y}:=\bigwedge^p \Omega^1_{X/Y}$.
\end{definition}

This construction can be thought as a generalization of the logarithmic de Rham complex as the following example shows.

\begin{paragraphe}
    Let $X$ be a normal algebraic variety over a field $k$ and let $D \subseteq X$ be a reduced divisor, we equip $X$ with the divisorial logarithmic structure associated to $D$. Then, the cotangent sheaf $\Omega^1_{X/\star}$ is equal to the sheaf of logarithmic de Rham differentials $\Omega^1_{\underline{X}}(\log D)$. If we have a morphism $X \to S_0$ that admits a chart whose monoid homomorphism is diagonal, we get the relative version.
\end{paragraphe}

\begin{paragraphe}
    Let $f : X \to S_0$ be a logarithmic semistable degeneration as in the paragraph \ref{example_log_semistable_degeneration}, we will determine the \'{e}tale-local structure of $\Omega^1_{X/S_0}$. As we have seen in that paragraph, the logarithmic structure of semistable type on $X$ has \'{e}tale-locally a chart given by $e_i \in \mathbb{N}^r \mapsto z_i \in \mathscr{O}_X$; in particular, this implies that the stalk of the ghost sheaf of $X$ is isomorphic to $\mathbb{N}^r$ so we have 
    \begin{center}
        $\displaystyle \Omega^1_{\underline{X}/\underline{S_0}} \oplus (\mathscr{O}_X \otimes \mathscr{M}_X^{\mathrm{gp}}) \cong \Omega^1_{\underline{X}/\underline{S_0}} \oplus \bigoplus_{i=0}^r e_i \mathscr{O}_X $
    \end{center}
    locally for the \'{e}tale topology. Let us write down the two type of relations in the definition of $\Omega^1_{X/S_0}$:
    \begin{itemize}
        \item The first type of relation implies that $\mathrm{d}z_i=z_ie_i$ for all $0 \leq i \leq r$. In other words, $e_i$ should be thought as the logarithmic derivative $\frac{\mathrm{d}z_i}{z_i}$.
        \item Since $f^\flat$ is given in a chart (at the level of ghost sheaves) by the diagonal homomorphism $\mathbb{N} \to \mathbb{N}^r$, there is only one relation of the second type, which is $\sum_{i=0}^r e_i =0$ --- or, with the other interpretation $\sum_{i=0}^r \frac{\mathrm{d}z_i}{z_i}=0$.
    \end{itemize}
    In particular, we see that the sheaf $\Omega^1_{X/S_0}$ is locally free of rank 2.
\end{paragraphe}

As in the classical case, the local freeness of the logarithmic cotangent sheaf in the example above is in fact just a consequence of the log smoothness of $X \to S_0$.

\begin{proposition} \label{cotangent_sheaf_locally_free_log_smooth}
    If $Y$ is locally noetherian and $X$ locally of finite type over $Y$, the sheaf $\Omega^1_{X/Y}$ is a coherent $\mathscr{O}_X$-module. If $X \to Y$ is in addition log smooth, $\Omega^1_{X/Y}$ is a locally free $\mathscr{O}_X$-module of finite rank.
\end{proposition}

We end this section with the computation of the dualizing sheaf of a logarithmic scheme with a log structure of semistable type over the standard logarithmic point. All definitions concerning duality theory can be found in \cite{Con00}.

\begin{proposition} \label{dualizing_sheaf_d_semistable variety}
    If $f : X \to S_0$ be a logarithmic scheme equipped with a logarithmic structure of semistable type, then $\omega_X = \Omega^{\dim(X)}_{X/S_0}$.
\end{proposition}

\begin{preuve}
    Firstly, the dualizing complex $\omega_Y^\bullet$ of $X$ is a \emph{sheaf} because $X$ is Gorenstein. The results now follows from \cite[Theorem 2.21]{Tsu99} since the morphism $X_0 \to S_0$ is vertical, which means that for all $x \in X$, the morphism of monoids $\overline{\mathscr{M}}_{S_0,f(x)} \to \overline{\mathscr{M}}_{X,x}$ is vertical in the sense of \cite[I., Definition 4.3.1.]{Og18}. The verticality of $f$ follows from the fact that this morphism is either an isomorphism or the diagonal map $\mathbb{N} \to \mathbb{N}^r$ which is vertical since the biggest proper faces of $\mathbb{N}^r$ are given by putting one of the summands of $\mathbb{N}^r$ to 0 and these faces do not contain the diagonal --- see \cite[I., Examples 1.4.8., 2.]{Og18}.
\end{preuve}

\subsubsection{The sheaf $\mathrm{LS}_X$.}

In this subsubsection, we recall the definition and some facts about the sheaf of sets $\mathrm{LS}_X$ associated to a normal crossing variety $X$. This sheaf encapsulates the local logarithmic structures of $X$ over the standard logarithmic point $S_0$ for which the structural morphism to $S_0$ is log smooth. All the proofs of the stated facts can be found in \cite[Section 5]{FRF21}.

\begin{definition} \label{definition_sheaf_LS}
    Let $X$ be a normal crossing complex variety equipped with a logarithmic structure. \\
    For any \'{e}tale map $U \to X$, we denote by $\Gamma(U \to X, \mathrm{LS}_X)$ the set of logarithmic semistable degenerations $U \to S_0$. This assignment defines an \'{e}tale sheaf of sets $\mathrm{LS}_X$ on $X$.
\end{definition}

In fact, the sheaf $\mathrm{LS}_X$ can be equipped with an action of $\mathscr{O}_X^\times$.

\begin{paragraphe}
    Let $f : U \to S_0$ be a local section of $\mathrm{LS}_X$. We have the following composition of morphisms of sheaves of monoids:
    \begin{center}
        $\underline{\mathbb{N}}_U \to \mathscr{M}_{S_0} \to \mathscr{M}_U$.
    \end{center}
    If we denote by $\rho \in \Gamma(U, \mathscr{M}_U)$ the image of 1 by that composition and if $u \in \Gamma(U, \mathscr{O}_U^\times)$, we define the element $u\cdot f \in \Gamma(U, \mathrm{LS}_X)$ as the unique log smooth morphism $U \to S_0$ such that the image of $1 \in \underline{\mathbb{N}}_U$ in $\Gamma(U, \mathscr{M}_U)$ is equal to $\rho \cdot \alpha_U^{-1}(u)$. This defines an action of the sheaf $\mathscr{O}_X^\times$ on the sheaf $\mathrm{LS}_X$. Since $X$ is a normal crossing variety, this action is transitive.
\end{paragraphe}

In fact, this sheaf $\mathrm{LS}_X$ has to do with the sheaf $\mathscr{T}^1_X:=\mathscr{E}xt^1_{\mathscr{O}_X}(\Omega^1_{\underline{X}}, \mathscr{O}_X)$ of local first order deformations of $X$ consider in the definition of d-semistability.

\begin{paragraphe} \label{facts_about_ls}
    We define a morphism of sheaves of sets $\mathrm{LS}_X \to \mathscr{T}^1_X$ as follows: if $f : U \to S_0$ is a local section of $\mathrm{LS}_X$, then, locally on $U$, a log smooth lifting $f_\varepsilon : U_\varepsilon \to S_{1,1}$ of $f$ to $S_{1,1}:=\mathrm{Spec}(\mathbb{N} \xrightarrow{1 \mapsto \varepsilon}\mathbb{C}[\varepsilon])$ exists uniquely and the conormal short exact sequence
    \begin{center}
        $0 \to \mathscr{O}_U \to i^*\Omega^1_{\underline{U_\varepsilon}} \to \Omega^1_{\underline{U}} \to 0$
    \end{center}
    associated to the closed immersion $i : U \hookrightarrow U_\varepsilon$ defines a section of the sheaf $\mathscr{T}^1_X$ on $U$. This morphism is in fact $\mathscr{O}_X^\times$-equivariant (\cite[Proposition 5.2.]{FRF21}) and is injective in the case of a normal crossing variety. One can prove that, in that case, its image coincide with the set of local generators of the $\mathscr{O}_{X_\mathrm{sing}}$-module $\mathscr{T}^1_X$ --- see \cite[Theorem 5.5]{FRF21}.
\end{paragraphe}

Let us write some implications of this in the example of the variety $X$ defined in the first section.

\begin{paragraphe} \label{global_sections_LS}
    Since the variety $X$ defined in the first section is d-semistable, we have $\mathscr{T}^1_X \cong \mathscr{O}_D$. By \ref{facts_about_ls}, the image of the morphism $\Gamma(X, \mathrm{LS}_X) \to \Gamma(X, \mathscr{T}^1_X) \cong \Gamma(D, \mathscr{O}_D)$ is isomorphic to $\mathbb{C}^\times$; this means that logarithmic structures of semistable type on $X$ are parameterized by $\mathbb{C}^\times$. One can write an explicit \'{e}tale-local chart for the global section of $\mathrm{LS}_X$ corresponding to $\lambda \in \mathbb{C}^\times$: it is given by the triple $(\mathbb{N}^2 \xrightarrow{a} \mathbb{C}[z_0,z_1,z_2]/(z_0z_1), \mathbb{N} \xrightarrow{b} \mathbb{C}, \mathbb{N} \xrightarrow{\delta} \mathbb{N}^2)$ where 
    \begin{itemize}
        \item $a$ is defined by $a(p,q):=\lambda^pz_0^pz_1^q$
        \item $b$ is defined by $b(1):=0$
        \item $\delta$ is defined by $\delta(1):=(1,1)$.
    \end{itemize}
\end{paragraphe}

\bibliography{bibliography}
\bibliographystyle{acm}

\end{document}